\documentclass[a4paper,11pt]{amsart}
\usepackage{multicol,ifthen}
\usepackage{amssymb,stmaryrd}
\usepackage[longnamesfirst]{natbib}
\usepackage{graphicx}
\usepackage[final]{pdfpages} 
\usepackage{caption}

\newcommand{\R}{\mathbb R}

\newcommand{\C}{\mathbb C}  
\newcommand{\M}{\mathbb M}

\newtheorem{thm}{Theorem}

\newcommand{\beq}{\begin{equation}}
\newcommand{\eeq}{\end{equation}}
\newcommand{\beqarr}{\begin{eqnarray}}
\newcommand{\eeqarr}{\end{eqnarray}}
\newcommand{\beqa}{\begin{eqnarray*}}
\newcommand{\eeqa}{\end{eqnarray*}}
\begin{document}

\thispagestyle{empty}
\renewcommand{\thefootnote}{\fnsymbol{footnote}}
\title[Weyl geometry in late 20th century physics]{The unexpected 
resurgence of Weyl  geometry in late 20-th century  physics}
\author[E. Scholz]{Erhard Scholz}
\date{2017  03/08}
\renewcommand{\thefootnote}{\arabic{footnote}}

\begin{abstract}
Weyl's original scale geometry of 1918 (``purely infinitesimal geometry'') was withdrawn by its author  from physical theorizing in the early 1920s. It had a comeback in the last third of the 20th century in different contexts: scalar tensor theories of gravity, foundations of gravity, foundations of quantum mechanics, elementary particle physics, and cosmology.  It seems  that Weyl geometry continues to offer  an open research potential for the foundations of  physics even after the turn to the new millennium.
\end{abstract}

\maketitle 

\setcounter{tocdepth}{2}
\tableofcontents

\section*{\small Introduction}
In the 1970s three groups of authors started, basically independent from each other, to reconsider Weyl's generalization of Riemannian geometry from 1918. Weyl  had proposed  the latter in the  perspective of building a geometrically unified theory of gravitation and electromagnetism. By the end of the 1920s,  after the successful reformulation of the underlying gauge idea in relativistic quantum physics,  most physicists including  Weyl himself  had  given up the idea of extending the geometry of spacetime by a ``localized'' scaling degree of freedom. It  was not to be expected that half a century later  researchers of the next generation would try again to give Weyl geometry  a new role in the changed context of  late 20th century physics. But some of them did.  A group of authors, in particular F. Ehlers, F. Pirani, and A. Schild,  used it as a conceptual framework for  clarifying  the foundations of gravity;  others explored extended gravity theories in the  generalized geometrical structure, and still others, like W. Drechsler and H. Tann in Munich,  investigated  connections between gravity and quantum physics. In several of these approaches  a  scalar field extending the gravitational structure played a crucial role. Although none of the attempts found an immediate broader response,  many of them led to follow up papers. In the result, different research perspectives exploring questions of recent physics from a Weyl geometric viewpoint emerged, but they remained  too heterogeneous for coalescing  to a  coherent literary tradition or even forming a common research community.

 The call for papers for the Mainz conference proceedings was   a splendid incentive for taking stock of the broader range of Weyl geometric investigations in physics, which took a new start in the last three decades of the 20th century.  Of course the following survey cannot be complete; it rather  has to be confined by specified boundaries. So this paper is restricted to the more classical parts of gravity with some, relative limited  outlooks at connections to quantum theory. 
{\em Not covered} in this survey is the whole range of Weyl geometric methods in Kaluza-Klein theories,  in supergravity, and in string theory.

In order to facilitate the  reading of the following survey, the paper starts with a very short  introduction to, or a reminder of, central features of Weyl geometry and gravity (section 1.1). Because a considerable amount of the following developments utilize a  scale covariant scalar field coupled to the Hilbert term similar to the one in Jordan--Brans--Dicke (JBD) gravity, the second part of the first section is devoted to a short glance at JBD theory from a Weyl geometric perspective  (section 1.2). The other sections give a partly historical, partly systematic survey  of the attempts for using Weyl geometric methods in recent physics.

In section 2 three different,  partially overlapping, approaches of the 1970s are described. The already mentioned paper of Ehlers, Pirani and Schild (EPS) on the foundations of gravity  and some follow up papers are dealt with in section 2.1. A completely different retake arose from proposals put forward  by a group of Japanese physicists, M. Omote, R. Utiyama et al. and independently by  P.A.M. Dirac. They investigated a scalar field coupling  to the  Hilbert term similar to JBD gravity, but in the  scale invariant approach of Weyl geometry. Dirac's and the Japanese physicists'  interpretations of the Weylian scale connection were not the same. They and their respective immediate successors had  different research contexts in mind, gravity, astrophysics, cosmology and electromagnetism in Dirac's case, nuclear and elementary particle physics in Utiyama's (section 2.2). Finally, although less noticed in the wider community, a specific road to Weyl geometric structures arose in the research on  gauge theories of gravity  arising from the Kibble-Sciama program of deriving gravitational structures (fields) from ``localizing'' symmetries in Minkowski space, often considered from a wider perspective  than that of  the Poincar\'e group. In this view Weyl geometry appeared as a special case of Cartan geometry, and Weyl geometric gravity ought generically be extended by a translational connection component, viz. torsion. It is a surprising fact that these three re-starts of Weyl geometric gravity, although arising from completely different backgrounds and pursuing different goals, were undertaken and published in the short interval  1971 -- 1974, exactly the time when the basics for the standard model of elementary particle physics were established (section 2.3). 

Before we come to the follow up investigations which made use of  these approaches in the standard model of elementary particle physics and/or in astrophysics and cosmology, we turn towards an even   more surprising recourse to Weyl geometry in attempts to  geometrize  quantum mechanics (QM) in the wake of the Bohmian heterodoxy (section 3). In order to make this kind of geometrization accessible to readers not versed in Bohmian quantum mechanics, the basic ideas necessary to understand the  geometrization proposals are shortly resumed in section 3.1. A survey of a peculiar road towards geometrizing  configurations spaces of QM by Weyl geometry, developed in the 1980s by  E. Santamato's and continued after the turn to the 2010s with his colleague F. De Martini, follows (section 3.2). A more fragile idea of a Bohm-type quantization procedure in cosmology leading to a Weyl geometric framework, proposed by  A. and F. Shohai and M. Golshani is the topic of section 3.3.

In  section 4 we turn towards different attempts at using Weyl geometric structures (mainly  scale invariance and the scale invariant affine connection)  and fields (Weylian scale connection and/or an additional scale covariant scalar field) in elementary particle physics. Three interrelated questions arise naturally if one wants to bring gravity closer to the physics of the standard model (SM): 
\begin{itemize}
\item[(i) ] Is it possible to bring conformal, or at least scale covariant generalizations of classical (Einsteinian) relativity 
 into a coherent common  frame with the standard model SM?\footnote{Such an attempt seemed  to be  supported experimentally by the phenomenon of (Bjorken) scaling in deep inelastic
 electron-proton scattering experiments. The latter indicated, at first glance, an active scaling symmetry of mass/energy in high energy physics; but it turned out to hold only approximatively and was of restricted range.}
\item[(ii)] Is  it  possible to embed classical relativity in a quantized theory of gravity or, the other way round, to derive classical relativity as an effective theory arising from a more fundamental quantum gravity theory at the classical level? 
\end{itemize}

The fact that all  the SM fields, with the only exception of the Higgs field, have conformally invariant Lagrangians, in the context of special relativity, i.e., Minkowski space, was considered among others by F. Englert and coauthors already in the mid-1970s. It cried out for investigations in a Weyl geometric perspective which then, of course, would invite  generalizing the spacetime environment of all SM fields, at least in their pseudo-classical form,\footnote{SM fields are here  called {\em pseudo-classical} if they are considered  before, or better abstracting from,  so-called second quantization. Mathematically they are classical fields (spinor fields or gauge connections), but the field components do not correspond to physically measurable quantities. Observationally relevant information can be extracted only after applying perturbative quantization methods.}
 to  Lorentzian or  Weylian  manifolds.  In this context the Weylian scale connection was identified by L. Smolin at the end of the 1970s  as a new, hypothetical, field which after quantization would lead to a  particle with mass close to the Planck scale. Roughly ten years later this particle was found again and called a ``Weylon'' by H. Cheng   (section 4.1). Again roughly a decade later the question of ``mass generation'' by breaking the scale symmetry in a Weyl geometric approach to  SM fields was studied at Munich by W. Drechsler and H. Tann (section 4.2). 

This question continued to attract the interest of researchers at least until the empirical detection of the Higgs boson in 2012.  In the last few years in particular H. Nishino and S. Rajpoot, but not only they, have studied the question of how the symmetry of the standard model may   be enhanced by a scale degree of freedom and may be broken by a peculiar interplay of an initially scale covariant scalar field and the ``Weylon''. All this was discussed at the pseudo-classical level  (section 4.3).
In the recent years some authors have  turned towards the difficult questions of Weyl scaling at the quantum level. A group of Italian authors, G. Codello, G. D'Orico, C Pagani, and R. Percacci brought forward new arguments with regard to the commonly shared view that scale symmetry is necessarily broken at the quantum level. They have proposed quantization procedures under which  scale invariance can be preserved under quantization. H. Ohanian has recently discussed the transition between a a scale invariant phase of fields close to the Planck scale to a lower energy regime with broken scale symmetry and Einstein gravity as effective field theory (section 4.4).

The rescaling allowed in Weyl geometry may change the geometrical picture underlying our usually assumed cosmological models. Scalar fields with conformal rescaling have been in use  for a long time in  ``early universe'' modelling (section 5.1). They invite  Weyl geometric investigations and  were dominated  for several decades by N. Rosen and M. Israelit, the first one an early  protagonist  of the Dirac approach to Weyl geometric gravity. In the last few decades also other authors jumped in with slightly different ideas (section 5.2f.). A  coherent tradition with a larger group of researchers in astrophysical and cosmological studies  has formed in Brazil around M. Novello. It  invokes a (weak) Weyl geometric framework and has defolded its research questions for more than two decades, more stable and with a wider group of contributors than any others line of research considered in this survey (section 5.4). But the question of  dark matter effects, if considered from the gravitational side, has to be measured at the successes of modified Newtonian dynamics, MOND. This remained outside the scope of the Brazilian school. First steps of reconstructing MOND-like phenomenology in a Weyl geometric approach to gravity, made recently at Wuppertal,  seem sufficiently striking  to include it here  (section 5.3).

A survey of a side-stream issue of recent research, as it is attempted
in this paper, cannot claim to tell a coherent, perhaps even success, story. It rather
has to collect views from necessarily heterogeneous perspectives and brings
them together in one panorama. In  this way it may invite for a look
backward and forward, in order to reflect on the development of methods and views in recent
mathematical and theoretical physics (section 6).

\section{\small Preliminaries: Weyl geometric gravity and Jordan-Brans-Dicke theory\label{section Preliminaries}}
\subsection{\small Weyl geometry and gravity\label{subsection Weyl geometry} }
\subsubsection{Basics of the geometrical framework}
Weyl geometry is a generalization of Riemannian geometry, arising from  two insights:  (i) The mathematical automorphisms of both, of Euclidean geometry  and  of special relativity, are the {\em similarities} (of Euclidean, or respectively of Lorentz signature) rather than the congruences. No unit of length is naturally given in Euclidean geometry, and likewise the basic structures of special relativity (inertial motion and causal structure) can be given without the use of clocks and rods.
(ii) The development of field theory and general relativity demands a conceptual implementation of this insight in a consequently {\em localized mode} (physics terminology).\footnote{In mathematical terminology, the implementation of a similarity structure happens at the {\em infinitesimal} level. A  discussion, given by Weyl later in his life, of the role of mathematical and physical autormorphisms can be found in  \citep{Weyl:similarity},  some aspects of it also in \citep[chap. III, sec. 14]{Weyl:PMNEnglish}.}
 In a more physical language  (i) and (ii) can be given the form of the postulate that fundamental field theories have  to be formulated  covariantly  under point dependent rescalings of the basic units of measurement, while the Lagrangian densities and the dynamical  laws (the ``natural laws'') are invariant under point dependent rescaling (see Dicke's postulate cited in section \ref{subsection JBD}). It remains an open question whether the resulting extension of the mathematical automorphism group of the theories may be of physical import, or whether it is purely mathematical refinement.

Based on these insights, Weyl developed what he called purely infinitesimal geometry ({\em reine Infinitesimalgeometrie}) building upon a conformal generalization of a  (pseudo-) Riemannian metric $g$ with coefficient matrix $(g_{\mu \nu}) $ with (point-dependent) rescaling $\tilde{g}(x)= \Omega  (x)^2 \, g(x)$ ($\Omega  $ a nowherere vanishing positive  function), and a scale (``length'') connection  given by a real valued differential form $\varphi = \varphi_{\mu } dx^{\mu}$ \citep{Weyl:GuE,Weyl:InfGeo}. If one rescales   the metric by $\Omega$   one has   to {\em  gauge transform} $\varphi$ by  $\tilde{\varphi} = {\varphi}-  d \log \Omega {} $. The {\em scale connection} $(\varphi_{\mu})$ expresses how to compare lengths of vectors (or other metrical quantities) at two infinitesimally close points, both measured in terms of  a representative $(g_{\mu \nu})$ of the conformal class. The typical symmetry of the geometry, at the infinitesimal level is thus the scale extended Poincar\'e group, sometimes called the {\em Weyl group} (although the same name is used in Lie group theory in a completely different sense).
In 1918 to roughly 1921/22 it seemed clear to  Weyl that this extension of Riemannian geometry can be used for unifiying gravity and electromagnetism; later he gave up this hope and considered his scale geometry as purely mathematical enterprise the most important features of which were transplanted to the $U(1)$-gauge theory of electromagnetism.\footnote{For more historical and philosophical details see, among others, 
\citep{Vizgin:UFT,Goenner:UFT,Ryckman:Relativity,Scholz:Connections,Scholz:DMV}.}

In hindsight, Weyl's generalization of Riemannian geometry may be embedded in E. Cartan's even  wider   program of geometries with infinitesimal symmetries. In the case of the scale extended Poincar\'e group one then arrives  at a Cartan-Weyl geometry with a translational Cartan connection and {\em torsion} as the typical extension of the structure.\footnote{For a modern presentation of Cartan geometry, including the Cartan-Weyl case, see, e.g., \citep[chap. 7]{Sharpe:Cartan_Spaces}; for the physical aspects of the extension studied since the 1970s  \citep[chap. 8]{Blagojevic/Hehl}.} 
With the exception of section \ref{subsection Cartan-Weyl} this paper will be restricted to the original form  of Weyl geometry without torsion; in large parts it even deals with the most simple case of an integrable scale connection. The reasons for this restriction will become apparent below. 

Metrical quantities in Weyl geometry are directly comparable only if they are  measured at the same point $p$ of the manifold. Quantities measured at different points $p \neq q$ of finite, i.e., non-infinitesimal, distance can be  compared metrically only after an integration of the scale connection along a path from $p$ to $q$.  Weyl realized that this structure is compatible with  a uniquely determined affine connection $\Gamma = (\Gamma ^{\mu }_{\nu \lambda } )$, the {\em affine connection of Weyl geometry}. If   the Levi-Civita connection of the Riemannian  part $g$  is denoted by ${}_g\Gamma^\mu _{\nu \lambda }$,  the Weylian affine connection is given by 
\begin{equation} \Gamma^{\mu }_{\nu \lambda } =   {}_g\Gamma^\mu _{\nu \lambda } + \delta ^{\mu }_{\nu } \varphi _{\lambda } +
\delta ^{\mu }_{\lambda } \varphi _{\nu } - g_{\nu \lambda } \varphi^{\mu }.\label{Levi-Civita}  
\end{equation} 
In the following the {\em covariant derivative} with regard to $\Gamma$ will be denoted as $ \nabla = \nabla_{\Gamma }$. Similarly the curvature expressions for the Riemann tensor, Ricci tensor and scalar curvature $ Riem,\, Ric,\, R$ will denote the Weyl geometric ones. The corresponding scale gauge dependent  Riemannian analogues derived from ${}_g\Gamma^\mu _{\nu \lambda }$ will be written as $_g\hspace{-0.2em}\nabla,\,  _g Riem,\, _g Ric,\, _g R$. The Weylian scalar curvature, e.g., is
\beq R = \, _g\hspace{-0em}R -(n-1)(n-2)\,\varphi_{\mu} \varphi^{\mu} - 2(n-1)\, _g\hspace{-0.2em}\nabla_{\mu} \, \varphi^{\mu} \, , \label{Weylian scalar curvature}
\eeq  
with $n$ the dimension of the manifold. 
A change of scale  neither  changes the connection (the left hand side of (\ref{Levi-Civita})) nor the covariant derivative;  only the composition from the underlying Riemannian part and the corresponding scale connection (right hand side) is shifted.

As every connection defines a unique curvature tensor, curvature concepts known from ``ordinary'' (Riemannian) differential geometry follow. The Riemann and Ricci  tensors, $Riem, Ric$, are  scale invariant by construction, although their expressions contain terms in $\varphi$. On the other hand,  the scalar curvature involves ``lifting'' of indices by the inverse metric and is thus scale covariant of weight $-2$ (see below). 

For vector and tensor fields (of dimensionful quantities) the appropriate scaling behaviour under change of the metrical scale has  to be taken into account. If a field, expressed by $X$ (leaving out indices) with regard to the metrical scale $g(x)$ transforms like
$ \tilde{X}= \Omega ^k X $  with regard to the scale choice $\tilde{g}(x)$ as above, $X$ is called a {\em scale covariant} field of {\em scale weight}, or {\em Weyl  weight} $w(X) := k$  (usually an integer or a fraction). It is the negative of the mass weights used in particle physics. 
 In general the covariant derivative, $\nabla X$,  of a scale covariant quantity $X$ is no longer scale covariant; but a scale covariance can be recovered. Adding a weight dependent term solves the problem. The {\em scale covariant derivative} $D$ of $X$ is  defined by $ DX := \nabla X + w(X)\varphi  \otimes X $,  in coordinate description 
\begin{equation} D_{\mu}X^{\nu} := \nabla_{\mu} X^{\nu} + w(X)\varphi_{\mu}  X ^{\nu}\; . \label{scale covariant derivative} \end{equation}

For example, the  derivative $\nabla g$ is not scale covariant, but $Dg$ is -- even with the result zero:
\beq Dg=\nabla g + 2  \varphi \otimes g =0\,  \label{metric compatibility}
\eeq
In Weyl geometry  the metric is thus  no longer  {\em constant} with regard to the  derivative $\nabla $  but {\em with regard to the scale covariant derivative $D$}. From the point of view of Riemannian geometry this appears as a ``non-metricity'' of the connection (in the literature often called ``semi-metricity''). From the Weyl geometric point of view  it is nothing but the {\em  metric compatibility condition} for $\Gamma$. 

Here it may suffice to have recalled these basic properties.   More  details on Weyl geometry can be found in Weyl's original papers  \citep{Weyl:GuE,Weyl:InfGeo},  those of his successors
\citep{Eddington:Relativity,Bergmann:Relativity,Dirac:1973} and  more recent 
literature.\footnote{Presentations of Weyl geometry can be found, among others, in 
\citep{Blagojevic:Gravitation,Israelit:1999Book},   \citep[appendix A]{Drechsler/Tann} and  \citep{Perlick:Diss} (difficult to access). For selected aspects see \citep{Codello_ea:2013} and  \citep[sec. 4]{Ohanian:2016}. Integrable Weyl geometry is presented in \citep{Dahia_ea:2008,Romero_ea:Weyl_frames,Almeida/Pucheu:2014,Quiros:2014a},  \citep[sec. 2.1]{Scholz:2011Annalen}.  Be aware of different conventions for the scale connection.   Expressions for Weyl geometric derivatives  and  curvature quantities are derived in  \citep{Gilkey_ea,Yuan/Huang:2013} and \citep[App.]{Miritzis:2004}. For  a more  mathematical  perspective consult  \citep{Folland:1970,%
Calderbank/Pedersen:1998,Gauduchon:1995,Higa:1993,Ornea:2001,Gilkey_ea}. \label{fn lit WG} }

\subsubsection{Weyl geometric gravity}
Weyl's generalization of Riemannian geometry arose with the perspective of  generalizing Einstein gravity, which would allow a geometrical unification of gravity and electromagnetism \citep{Vizgin:UFT,Goenner:UFT}. Any meaningful Lagrangian in this framework underlies the constraint of scale symmetry. Because of   $w(\sqrt{|g|})=4$, while the Weyl geometric scalar curvature $R$ is of weight $w(R)=-2$, Weyl could not work with the Hilbert-Einstein term but considered  quadratic expressions in the curvature terms for a generalization of the gravitational Lagrangian, e.g.\footnote{As the ``most simple and natural'' expression $\alpha_2=0$ in \citep{Weyl:GuE} and $\alpha_1=0$ as the most simple example in \citep[4th ed., 5th ed.]{Weyl:RZM}.}
\beq \mathfrak{L_W}= L_W \sqrt{|g|}   \quad \mbox{with}\; L_W = \alpha_1 R^{\mu}_{\; \nu \lambda \kappa}  R_{\mu}^{\; \nu  \lambda \kappa} + \alpha_2 \, R^2 \, . \label{Weyls grav Lagrangian}
\eeq 
Of course,  he added a term in the scale curvature  $f=d\varphi$  ($f_{\mu \nu} = \partial_{\mu} \varphi_{\nu}  -   \partial_{\nu} \varphi_{\mu}$) looking like the  Maxwell action:
\beq \mathfrak{L_f}= L_f \sqrt{|g|}   \quad \mbox{with}\; L_f = - \frac{1}{4}f_{\mu \nu}f^{\mu \nu} \label{fsquare Lagrangian}
\eeq 
Only much later -- in fact, about half a century later -- other gravitational Lagrangians with Weyl geometric scale symmetry started to be considered. They arose from the idea of a coupling between gravity, here the Weyl geometric scalar curvature $R$, and a scalar field  with ``correct'' complementary weight (subsection \ref{subsection Dirac-Utiyama}).

 In the period covered here we encounter two modes of Weyl geometric gravity. One is farther away from Einstein gravity and uses  square curvature Lagrangians, sometimes called {\em Weyl gravity} (in the strong sense);  the other is  closer to  Einstein gravity and works with a modified Hilbert term coupled to a scalar field, in the physics literature it is often called {\em Weyl geometric scalar tensor theory} (WST). The latter goes back to independent  proposals by M. Omote and R. Utiyama on the one hand and P.A.M. Dirac on the other for making use of a scalar field modification of the Hilbert term, analogous  to Jordan-Brans-Dicke theory (see section \ref{subsection Dirac-Utiyama}). Here the gravitational structure is characterized by an equivalence class of triples $(g, \varphi, \phi)$, with $g=g_{\mu\nu}dx^{\mu}dx^{\nu}$ the {\em Riemannian component} of the  Weylian  metric,  $\varphi=\varphi_{\mu}dx^{\mu}$ its  scale connection,  and $\phi$ an additional scalar field.\footnote{In a way, this may be called a  geometrical ``tensor-vector scalar'' theory {\em sui generis}, in which all components have  geometrical meaning.} 
The equivalence is given by combined rescaling transformations  $g \mapsto \tilde{g} = \Omega^2 g= e^{2\omega} g, \; \varphi \mapsto \tilde{\varphi}= \varphi - d  \omega ,\; \phi\mapsto \tilde{\phi}= e^{-\omega}\phi$. Because of scaling freedom, a  Weylian metric with nowhere vanishing scalar curvature can be gauged to  $R\doteq const$. Here, and elsewhere in this paper, $\doteq$ denotes an equality which holds only in a certain gauge  specified by the context. Weyl considered this as the ``natural'' gauge,   we  prefer to call   it the {\em Weyl gauge}.  

If the  scale connection is  an exact form, 
\beq
\varphi=-dw \, , \label{integrable scale connection}
\eeq
  with a scalar potential  $w$ scale transforming by $w\mapsto \tilde{w}=w+ \omega$, we  work in an {\em integrable Weyl geometric scalar tensor} theory (IWST). Then the gravitational structure  reduces to the Riemannian component of the metric plus, at face value, two scalar fields $(g,\phi=e^v,w)$ with equivalence under rescaling. As $v\mapsto\tilde{v}=v-  \omega$, the sum $v+w$  is a {\em scale invariant} scalar field of the gravitational structure and the {\em only crucial} one. Because of the scale gauge freedom $\phi$, respectively $v$ or $w$, can be given any chosen value, e.g. a  constant.\footnote{The dynamical consequences of this interdependence have been clarified by  \citep{Israelit:1999Book,Israelit:1999matter_creation}, see section \ref{subsection Israelit}. }  In the integrable case, two scale gauges are of particular importance in addition to  Weyl gauge: The {\em Riemann gauge} in which the scale connection is ``integrated away'' (for    $\omega={-w}$), then  $\tilde{\varphi }\doteq 0$. The other one  is the  {\em scalar field gauge} in which the scalar field is scaled to a constant, $\tilde{\phi}(x) \doteq \phi_o = const$ (for $\omega=v$). If the  value  of $\phi_o $ is specified such that it hooks up to Einstein gravity, $\phi_o=(8\pi G)^{-1}$ (up to a hierarchy factor if need be), it is  called {\em Einstein gauge}.\footnote{Obviously the Einstein gauge exists also in the non-integrable case.} 

For a vanishing  scale invariant sum, 
\beq
v+w=0\, ,  \label{v+w=0}
\eeq
the scalar field $\phi$ is essentially the potential for the scale connection, more precisely, 
\beq  \varphi = d v=  d \ln \phi \, . \label{triviality condition}
\eeq 
Then and only then,  Einstein gauge and Riemann gauge coincide and  IWST reduces to Einstein gravity. The Palatini approach varying the metric and the affine connection of a Lagrange density $\mathfrak{L}=\phi^2 R \sqrt{|g|}$ independently  enforces the constraint $v+w=0$  in addition to the integrability of the scale connection. This implies a  reduction of a Palatini-IWST to Einstein gravity.  The latter is then only  re-written in scale covariant form, but without any modification of the dynamics.\footnote{Cf. sections \ref{subsection trivial}, \ref{subsection Brazilian approach}. }
 If one considers IWST from the point of view of the metric-affine  scheme, one better   uses   variational constraints like in \citep{Cotsakis/Miritzis:Palatini} rather than the Palatini approach. Then the condition (\ref{v+w=0}) is not enforced and the scalar field, respectively the integrable scale connection  (\ref{triviality condition}), express an additional dynamical degree of freedom.

\subsection{\small Jordan-Brans-Dicke (JBD) gravity\label{subsection JBD}}
\subsubsection{Basics of JBD theory}
At the turn to the 1950s  {\em Pascual Jordan}  (Hamburg) and, a decade later, {\em Carl Brans} and {\em Robert Dicke} (Princeton) proposed a generalization of Einstein gravity by considering  a varying gravitational parameter. The motivations at Hamburg and at Princeton were different, but there was a wide overlap of the ensuing theory, here abbreviated by JBD. Jordan  started from an action principle \citep[p. 140]{Jordan:Schwerkraft}
\beq \mathcal{L_{J}}(\chi, g) = (  \chi R -  \frac{\xi}{\chi }\partial ^{\mu} \chi\, \partial _{\mu} \chi  ) \sqrt{|det\, g|} \, , \label{Lagrangian JBD} \eeq
with a   parameter $\xi$ and a real scalar field $\chi$ functioning as a kind  of  spacetime dependent (reciprocal)  gravitational ``constant'', $R$ here of course the Riemannian scalar curvature of the metric $g$ \citep[2nd. ed., 163, (3)]{Jordan:Schwerkraft}.\footnote{Warning: One has to check carefully the sign convention used in the definition of Riemann and scalar curvature. Jordan, e.g., used  sign inverted definitions of the curvature terms with respect to those  used here and in much of the present literature \citep[40]{Jordan:Schwerkraft}. In Fujii/Maeda's notation (see below) this would correspond to $\epsilon =-1$ and thus to a ``ghost'' field.}
 A Lagrange term $L_m$ for classical matter   could be foreseen
 (e.g., \citep[equ. (6)]{Brans/Dicke}). The hypothesis of a  ``varying gravitational constant'' had been brought up already more than  a decade earlier by P.A.M. Dirac, when he speculated about ``large numbers'' relations in physics.\footnote{For Dirac's role in this story see \citep{Kragh:VaryingGrav}, for a larger view at JBD theory \citep{Brans:Tenacious,Brans:JBD}.} 
  Pauli  reminded Jordan  that his ``extended gravity''  allowed for a class of conformal  transformations  which not only affect the metric but also the scalar field, 
\beq   \tilde{g}_{\mu \nu} = \Omega ^2 g_{\mu \nu} \, , \quad \tilde{\chi }= \Omega ^{-2} \chi \; . \label{Pauli conformal transformations}
\eeq
 Jordan  included this  generalization into the second edition of his book  
 \citep[2nd ed.,169]{Jordan:Schwerkraft}.\footnote{The conformal factor $\Omega $ was (unnecessarily) restricted by the condition  $\Omega ^2 = \chi ^{\gamma }$ for some constant $\gamma \in \R$.}

A few years later, Robert Dicke and Carl Brans   restarted the study of scale covariant scalar fields (including a  classical matter term $L_m$ in (\ref{Lagrangian JBD})  \citep[8]{Brans/Dicke}. Their motivation was to formulate a theory of gravity which took account of  Mach's  principle as understood by D.W. Sciama.\footnote{In the 1950s  Sciama had proposed to consider the possibility that the gravitational ``constant'' was related to the mass and the ``radius'' of the visible universe.}
For the two US physicists the main function of the scalar field was ``the determination of the local value of the gravitational constant'' \citep[929]{Brans/Dicke}. More clearly than in Jordan's work, the wave character of the dynamical equation of $\chi $ was emphasized by them \citep[equs. (9), (13)]{Brans/Dicke}. Moreover, they had a different view of the role of scale transformations. 

Their methodological goal  was a scale independent foundation of physical theories, with  a ``passive''   interpretation of scale transformations  in mind  \citep[927]{Brans/Dicke}, while Jordan and Pauli tended to think in  terms of  ``active'' scale transformations of material structures.  Dicke started
an  article dedicated  to  {\em transformations of units} in GRT \citep{Dicke:1962}   announcing as ``evident' the following principle:
\begin{quote}
It is evident that the particular values of the units of mass, length, and time employed are arbitrary and that the laws of physics must be invariant under a general coordinate-dependent transformation of units. \citep[2163]{Dicke:1962}
\end{quote}

That was very much in the spirit of Weyl's intentions of 1918, from which the latter  had disassociated himself with the shift of his gauge idea to quantum physics Weyl discussed this new view at different occasions in the 1940s,  e.g. in \citep[p. 165]{Weyl:similarity}.\footnote{Also in the  English edition of {\em Philosophy of Mathematics and Natural Sciences} Weyl  expressed this disassociation quite clearly, appealing to the constants of atomic physics which regulate the frequencies of spectral lines \citep[83]{Weyl:PMNEnglish}. But this was only one part of his perspective. In the appendix he  argued that for a deeper insight  it would be necessary to   understand how the  ``adaptation'' of the mass of the electron to the local field constellation is achieved \citep[288f.]{Weyl:PMNEnglish}. This was close to the intentions of his 1918 approach, although no longer a claim that the goal had been achieved. Einstein, in his later papers, agreed \citep[555f.]{Einstein/Schilpp}; see \citep{Lehmkuhl:2014}. }
It seems that  Dicke ``reinvented'' the idea of scale gauge invariance of the natural laws anew. He systematically discussed the scale transformations of physical quantities, based on the (quasi-axiomatic) principle of the invariance of the velocity of light $c$ and the Planck constant $\hbar$.  In particular, ``all three quantities, time, length, and reciprocal mass transform in the same way'' \citep[2164]{Dicke:1962}, i.e., 
\[ l' = \Omega \, l\, , \quad t' = \Omega  \, t \, , \quad m' = \Omega ^{-1}\, m \, .  \] 
In this sense, Weyl's scale gauge transformations  reappeared in the principles of Jordan-Brans-Dicke theory without being mentioned as such. It may be that at the time nobody  but Pauli was aware of  this close resemblance to Weyl's  theory.  Weyl's choice of (scale) gauge was translated by Dicke into the choice of a {\em frame} of measuring units, complementing the choice of a coordinate system.

In more recent papers the scalar field and the JBD parameter are written in slightly different form. With $\phi = \sqrt{2 \xi^{-1}\chi} $, scale weight   $w(\phi)= -1$,  and  $\xi=  \frac{\epsilon }{4 \omega }$  the Lagrangian (\ref{Lagrangian JBD}) turns into 
\beq  \mathcal{L}_{BD}= \left(\frac{1}{2} \xi  \phi^2 R - \frac{1}{2}\epsilon  \partial ^{\mu }  \phi \, \partial _{\mu}\phi  + L_{mat} \right) \sqrt{|det\, g|}     \, ,\label{Lagrangian JBD Fujii/Maeda} \eeq 
where  $\text{sig}\, g = (3,1) \sim (-+++)$ and generally $\epsilon=1$, while only in exceptional cases   $\epsilon = -1$ or  $0$ \cite[p. 5]{Fujii/Maeda}.\footnote{ $\epsilon =1 $ corresponds to a normal field having a positive energy, in other words, not a ``ghost''.  Fuji/Maeda add that $\epsilon =-1$  looks  unacceptable because it seems to indicate negative energy, but ``this need not be an immediate difficulty owing to the presence of the nonminimal coupling'' (ibid.).}   
In the following discussion this notation will be used as a  standard. 

A famous exception with $\epsilon=-1$ is the   special constellation of coefficients 
\beq L_{cc}= \phi^2 R + 6 \partial_{\mu}\phi \partial^{\mu}\phi  \, . \label{conformal coupling}
\eeq 
Then  the Lagrangian is invariant (up to an exact differential)  under conformal transformations \citep{Penrose:1965}. This case of   {\em conformal coupling} allowed to study  versions of gravity theory  ``in which scale invariance of matter is a consistency requirement on its coupling to gravitation'' \citep{Deser:1970}.  Deser  considered a conformally coupled scalar field as a paradigmatic example for matter and observed that the addition of a quadratic  term of the form $\frac{1}{2}\mu^2 \phi^2$, with $\mu$ a parameter of mass dimension, implies {\em breaking} of the {\em conformal symmetry}. In this case, the range  $\phi$ ``must clearly be cosmological in order not to lead to a clash with observation '' \citep[p. 252]{Deser:1970}. 

But in general, the three founding authors, Jordan, Brans and Dicke, considered it as evident that the ``conformal transformations''  (scale transformations) do {\em not} reduce the geometrical considerations to those of a purely conformal structure. They rather considered it as clear that  JBD theory  possesses a  covariant derivative $\nabla$, specified by the   reference metric $g$  underlying  (\ref{Lagrangian JBD}) from which Jordan  and Brans/Dicke started. Later this  scale was  called  {\em Jordan frame} (although Jordan was undecided, which scale might be  the ``natural''one). Because the Levi Civita connection of the Jordan frame metric determines the free fall trajectories of test particles,  many authors consider this one  as the ``physical frame'', the  other frames then appear as  mathematical auxiliary devices. On the other hand,  the JBD-field $\phi$ can be scaled to a constant. Then the gravitational part of the Lagrangian looks like the Hilbert term of Einstein gravity, while the remnants of the JBD scalar field appears in additional expressions of the Lagrange density.\footnote{See, e.g., \citep[chap. 3.6]{Capozziello/Faraoni}.}
 The resulting {\em Einstein frame} satisfies the Riemannian ``energy conservation'' condition for matter tensors. Since roughly the 1990s it has  found an increasing number of  supporters who now propose it as the proper frame for a  ``physical'' interpretation of JBD gravity. But no consensus in the JBD community has been achieved; the discussion  has remained undecided, to say the  least,     \citep{Faraoni/Nadeau,Quiros_ea:2013}.

\subsubsection{JBD in a Weyl geometric perspective}
The perspective of (integrable) Weyl geometry may help clarifying some aspects underlying this debate.  Let us denote the affine connection referred to by 
\beq \nabla := \,  _g\hspace{-0.2em}\nabla \, , \label{JBD affine connection}
\eeq
 where the r.h.s. expresses  the   Levi-Civita connection of $g$  in the JBD Lagrangian (\ref{Lagrangian JBD}). $\nabla$   is kept  unaffected, i.e. {\em invariant}, under scale transformations in JBD theory.  A structural view of Weyl geometry shows that the combination of a conformal structure $[g]$ of pseudo-Riemannian metrics $g$ and a specification of an invariant affine connection $\nabla$ with a compatibility condition, inbuilt here because of  (\ref{JBD affine connection}), determines a Weyl structure on a differentiable manifold $M$. In this way JBD gravity may be  embedded in the theoretical frame of Weyl geometry, independent of whether or not a single author knows.\footnote{See the  discussion in \citep{Quiros_ea:2013} and \citep{Scholz:2016Paving}.} Usually this is not being done (see, however, section \ref{subsection trivial}).

\section{\small Contributions to Weyl geometric gravity in the 1970s and 1980s\label{section WGG 1970s}}

\subsection{\small Ehlers/Pirani/Schild and subsequent work\label{subsection EPS}}
\subsubsection{An axiomatic approach to the foundations of gravity}
Weyl already discussed the relation between the physical concept of a {\em causal structure} and the mathematical concept  of a {\em conformal structure} on a differentiable manifold  \citep[4th. ed., appendix I]{Weyl:RZM}. He deemed it inadequate to think of an empirical determination of the metrical coefficients $g_{\mu \nu }$ by ``rods and clocks'' and looked for another empirical specification of a Weylian metric $(g, \varphi)$.
In a note added to a  letter to F. Klein\footnote{\citep{WeylanKlein:28Dez1920}}  (a little later published in  {\em G\"ottinger Nachrichten} as 
 \citep{Weyl:projektiv_konform})  he sketched an idea how this can be achieved. 
Assuming  his framework of the generalized ``purely infinitesimal'' geometry, Weyl showed that two of his generalized  metrics which have  identical  conformal structure and  the same projective geodesic path structure will coincide. This meant that, at least in the framework  of Weyl geometry, conformal and projective path structures  specify a  Weylian metric uniquely.\footnote{Weyl's note  \citep{Weyl:projektiv_konform} became better known by his calculation and discussion of projective and conformal curvature tensors, which followed.}

Weyl's argument on the combination of projective and conformal structure was taken up  and extended  by {\em  J\"urgen Ehlers, Felix Pirani and Alfred Schild} (EPS in the sequel), about the same time in which Dirac studied Weyl geometry in the context of scalar tensor theories
\citep{EPS}.  This paper  was written for a {\em Festschrift} in the honour of J.L. Synge.
Synge had become known for his proposal to base general relativity  on the behaviour of standard clocks ({\em chronometric approach}). From the foundational point of view,  clocks could appear  as a problematic choice, because they are  realized by complicated material systems. The question arose whether more basic signal structures of gravitational theory (light rays, particle trajectories) might do the job.

Using Hilbert's words, EPS ``laid the foundations  deeper'', combining  Weyl's idea of 1920/21 and the recently developed mathematical language and symbolic technology of differentiable manifolds with Hilbert's axiomatic method.\footnote{See \citep{Trautman:EPS}.} 
They started from three sets, $\mathcal{M}=\{p, q, \ldots\}, \, \mathcal{L}=\{L, N, \ldots\}, \,\mathcal{P}=\{P, Q, \ldots\} $, with $ \mathcal{L}, \mathcal{P} \subset \mathcal{M}$, and called the three sets respectively collections of {\em events, light rays} and {\em particles}. By postulates close to physical experimental concepts of light signal exchange between particles EPS formulated  different groups of axioms  in the Hilbertian style of foundations of geometry ($D_1, \ldots D_4, \, L_1,  L_2, \, P_1,  P_2, \, C$), which allowed them to introduce
 a $C^3$ differentiable structure on  $\mathcal{M}$ on which $ \mathcal{L} $ and $ \mathcal{P}$ then described smooth curves (axiom group $D$). Moreover,  
 a $C^2$ conformal structure was defined by  $\mathcal{L}$ (axioms $L$),
 and a differentiable projective path structure  by $\mathcal{P}$ (axioms $P$).
 
With a compatibility axiom $C$, basically postulating that light rays can be approximated arbitrarily well by particle trajectories, EPS could derive their main result.\\[-1.2em]
{\thm[Ehlers/Pirani/Schild 1972] A light ray structure $\mathcal{L}$ and a set of particle trajectories $\mathcal{P}$ defined on an event set 
$ \mathcal{M}$  which satisfy axioms $D, L, P, C$ endow $ \mathcal{M}$ with the structure of a ($C^3$-) differentiable manifold $M$ and a  ($C^2$-) Weylian metric  $[ (g,\varphi ) ]$. The latter is uniquely determined by the condition that its causal  and geodesic structures coincide with  $ \mathcal{L}$ and $ \mathcal{P}$ respectively.     }\\[-0.8em]

EPS posed the question, how a (pseudo-)Riemannian structure of classical (Einsteinian) relativity might arise from the Weylian one. A simple additional {\em Riemannian axiom},  postulating the vanishing of  the scale curvature, $d\varphi=0$,  could serve the purpose. Such a postulate did not seem nonsensical, as Weyl's interpretation of the  scale connection $\varphi$ as electromagnetic (e.m.) was obsolete anyhow and EPS did not adhere to it. But the authors did not exclude the possibility that a scale connection field $\varphi$ of nonvanishing scale curvature might play the role of a ``true'', although still unknown,  field.

\subsubsection{Subsequent work}
The paper of Ehlers, Pirani and Schild triggered a  line of investigations in the foundations of general relativity, sometimes called the {\em causal inertial approach} (Coleman/Kort\'e), sometimes subsumed under the more general search for a {\em constructive axiomatics} of GRT (Majer/Schmidt, Audretsch, L\"a\-m\-merzahl, Perlick and others). These investigations  turned towards a  basic conceptual analysis from the point of view of foundations of inertial geometry \citep{Coleman/Korte:inertial_conformal}, some even looking for  Desargues type characterization of free fall lines \citep{Pfister:Newtons_law}.
 How  a kind of   ``standard clocks'' can be introduced in the Weyl geometric setting  without taking refuge to atomic processes, by   just  using the observation of light rays and inertial trajectories, was studied by 
\citep{Perlick:Diss,Perlick:1987,Perlick:Observerfields}.
Another line of  follow up works explored the extension of the foundational argument of the causal inertial approach to quantum physics, where  particle trajectories  might  no longer appear acceptable as a foundational concept. 

This debate was opened by  \citep{Audretsch:1983}. It was soon continued by the collective work again of three authors \citep*{AGS}, cited in the sequel by AGS, and had  follow up studies, among them \citep{Audretsch/Laemmerzahl}. Audretsch  argued that the ``gap'' between Weylian and Riemannian geometry  can ``be closed if quantum theory as a theory of matter is made part of the total scheme'' \citep[2872]{Audretsch:1983}. He postulated  that quantum theory in the sense of Dirac or Klein-Gordon (K-G) fields on a Weylian manifold are  compatible with the latter's geometry,  
 if and only if the  WKB  (Wentzel-Kramers-Brillouin)  approximation  of  the  Dirac (or K-G) field leads  to streamlines which in the limit $\hbar \rightarrow 0 $ agree with geodesics  (Audretsch's {\em compatibility condition}).  

Working with scale covariant mass factors $m$ of Weyl weight $w(m)=-1$ 
Audretsch found that compatibility is possible only if the mass factor $m$ of the Dirac particle has vanishing covariant derivative, $\nabla_{\mu} m = 0$ in some gauge. He  observed that this implied  vanishing of Weyl's scale curvature
$  d \varphi = 0 $ \citep[equ. (6.14)]{Audretsch:1983} and concluded a bit rash:
\begin{quote}
{\em The consequence of the requirement is therefore that the Weyl space reduces to a Riemann space and the gap {\em [between Weylian and Riemannian geometry, ES]} described in Sec. I}  {\em is closed.} \citep[2881, emph. in original]{Audretsch:1983}
\end{quote}
 Because  $D_{\mu}m=0$ in any gauge if $\nabla_{\mu} m$ vanishes in Riemann gauge, Audretsch had only shown that  the limiting condition for streamlines of the WKB approximation of the Dirac field to classical geodesic trajectories implied {\em integrability} of the Weylian metric. The question whether this would also imply the choice of the Riemann gauge as ``physical''   was not  posed; it rather was imputed as  self-evident. 

In the AGS paper this question was taken up again  and stated carefully in the language of conformal fibre bundles for Dirac-  and for 
Klein-Gordon fields.  AGS showed that Audretsch's compatibility condition implies the possibility to  reduce the ``conformal'' group, here understood as  $\R^+ \times SO(1,3)$, to the orthogonal group. 
\vspace{-0.5em}
{\thm[Audretsch/G\"ahler/Straumann 1984] 
A Weylian manifold $(M,[(g,\varphi)])$ of Lorentzian signature  is locally  integrable, iff the  WKB 
 (Went\-zel-Kramers-Brillouin)  approximation  of  a  (locally defined) Dirac or Klein-Gordon field $\psi$ on $M$ leads  to streamlines which agree with geodesics in the limit $\hbar \rightarrow 0 $.
}

 The three authors formulated their consequence more carefully than Audretsch had done in his first  paper. They did not claim that their  investigation  had completely filled the gap to Riemannian geometry. 

All in all,   the three authors gave a more precise and mathematically modernized presentation of Audretsch's insight.  The gap between the Weylian and the Riemannian structure in the foundations of GRT was reduced but not completely closed. It  could seem natural to choose the Riemann gauge of the Weyl metric in order to reduce the structure group to $SO(3,1)$, but nothing compelled to do so.
The classical interpretation of geodesics as trajectories of mass points was foreign to the field theoretic context anyhow. It was now substituted by postulating  coherence between geodesic structure and the  flow-lines associated to pseudo-classical quantum fields.  

\subsection{\small Dirac's and Omote/Utiyama's  retake of Weyl geometry\label{subsection Dirac-Utiyama}}
\subsubsection{Dirac on  scale covariant ``varying'' gravity\label{subsection Dirac}}
In the 1970s {\em P.A.M. Dirac} introduced  Weyl geometry into the discourse of the  rising scalar-tensor theories. He was still  fascinated by the of interrelation of certain constellations of large numbers in physics, the ``large number hypothesis'' \citep{Dirac:1973,Dirac:LNH}.\footnote{Dirac presented his  proposal for a retake of Weyl geometry at the occasion of the symposium  honouring  his 70th birthday, 1972 at Trieste. This talk remained unpublished. According to \citep[p. 249 footnote]{Charap/Tait} the talk  was close to his 1973 publication.  For the broader historical context of this enterprise, the background in Dirac's reflection on large numbers in the 1920s,  and a surprising link to geophysics  see \citep{Kragh:VaryingGrav}.}
Largely following Eddington's notation and terminology \citep{Eddington:Relativity}, he introduced the readers to   Weyl geometry which was no longer generally known in the younger generation of physicists. He  then introduced a scalar-tensor theory of gravitation coupled in an  ``oldfashioned'', i.e. outdated, way to  electromagnetism. Like Weyl in 1918, he identified the potential of the electromagnetic field $F_{\mu \nu}$ with the Weylian scale curvature $f= d\varphi$
\[  F_{\mu \nu} = f_{\mu \nu} \; . \] 
In the sequel I  call this the {\em electromagnetic    (em) dogma} .
On the other hand,  he replaced  Weyl's original  gravity  Ansatz in the Lagrangian (using square curvature terms) by a JBD-type  Lagrangian using a 
real scalar field $\beta$ of weight $w(\beta)=-1$. He  added a biquadratic scale invariant potential  term,
\beq    \mathcal{L}_{Dir\, o} =-  \beta^2 R + k D^{\lambda}\beta\, D_{\lambda}\beta + c \beta^4 + \frac{1}{4} f_{\mu \nu}  f^{\mu \nu} \; ,     \label{scale invariant Dirac action}  \eeq
with   constant $k$. 
For $k=6$, the scale connection terms of the Lagrangian essentially cancel (i.e., they reduce to boundary terms and thus are  variationally negligible) like in (\ref{conformal coupling}). So Dirac wrote   the  Lagrangian   in the form
\beq       \mathcal{L}_{Dir\, 1} =  -\beta^2\,  _g\hspace{-0.1em} {R} + 6 \partial ^{\lambda}\beta\, \partial _{\lambda}\beta + c \beta^4 + \frac{1}{4} f_{\mu \nu}  f^{\mu \nu} \; ,                    \label{Dirac action} \eeq 
 known to be conformally  invariant \citep{Penrose:1965}, if  $_g\hspace{-0.1em}{R} $ denotes   the {\em sign inverted} Riemannian scalar curvature with respect to the generally accepted convention, while in (\ref{scale invariant Dirac action}) ${R}$ is the sign inverted scalar curvature of the  Weylian metric.\footnote{The qualifications ``sign inverted''  and  ``generally accepted''refers to the  sign  convention which agrees with the coordinate free definition $Riem(Y,Z)\,X=\nabla_Y \nabla_Z X - \nabla_Z \nabla_Y X - \nabla_{[Y,Z]}X$. It is  preferred in the mathematical literature including \citep[5th ed., 131]{Weyl:RZM} and also used in the majority of the more recent physics books. 
The ``sign inverted''  convention in some of the physics literature goes  back  to Einstein, e.g. \citep[801]{Einstein:GRT1916}, who in turn may have followed Ricci and Levi-Civita. It was continued in much of the physics literature of the first half of the 20-th century,   \citep[\S\, 37]{Eddington:Relativity},  \citep{Pauli:1921}  up to the influential  \citep[equ. (2.1.3)]{Weinberg:Cosmology_1972}.  Weyl, on the other hand, used the above convention long before the coordinate free definition of the Riemann tensor was available.  It seems to be  dominant in the more recent literature on GRT, although  Rindler speaks of a 50 \% distribution among the two conventions \citep[219]{Rindler:Relativity}.   \label{fn curvature conventions}}

Dirac derived dynamical equations and Noether identities for diffeomorphisms and scale transformations.
He distinguished Riemann gauge, $\varphi =0$, which existed, of course, only for vanishing e.m. field $f_{\mu \nu}=0$, from ``Einstein gauge'' (gravitational parameter constant, $\beta =1$) and ``atomic gauge'' (Weyl's natural gauge). He  warned that 
``all three gauges are liable to be different''  \citep[411]{Dirac:1973}. 

In a  discussion at the end of his article, why one should believe in the proposed ``drastic revision of our ideas of space and time'', Dirac announced a part of his research agenda, which was {\em independent} of the large number hypothesis:
\begin{quote}
There is one strong reason in support of the theory. It appears as one of the fundamental principles of Nature that the equations expressing basic laws should be invariant  under the widest possible group of transformations \ldots  The passage to Weyl's geometry is a further step in the direction of widening the group of transformations underlying physical laws [in addition to general coordinate transformations, E.S.]. 
One now has to consider transformations of gauge as well as transformations of curvilinear coordinates and one has to take on's physical laws to be invariant under all these transformations, which imposes stringent conditions on them.   \citep[418]{Dirac:1973}
\end{quote}

So far, Dirac's explanations agreed with the view of C. Brans and R. Dicke. 
He  followed a tendency of the time  for probing possible extensions of  the symmetries (automorphisms) of fundamental physics. 
In distinction to Pauli's insistence on a preferred scale, taken over  into the general discourse of JBD theory, Dirac  argued that at least three different gauges, Riemann, Einstein, and ``atomic'' gauge, indicated by atomic clocks, had to be considered in different theoretical or observational contexts. He saw no chance of a single preferred gauge; but sometimes the ``atomic'' gauge was assumed to be identical with Weyl gauge.\footnote{Weyl had argued that the atomic clocks somehow adapt to the local field constellation via the Weylian scalar curvature.}

\subsubsection{Some remarks on Dirac's followers\label{subsection Dirac followers}}
  Dirac's proposal for reconsidering Weyl geometry in a modified theory of gravity was taken up by field theorists and a few astronomers.
An immediate and often quoted paper by {\em Vittorio Canuto} and coauthors gave a broader and more detailed introduction to Dirac's view of Weyl geometry in gravity and field theory \citep{Canuto_ea}. The opening remark of the paper   motivated  the renewed interest in Weyl geometry with  actual developments in high energy physics:
\begin{quote}
In recent years, owing to the scaling behavior exhibited in high-energy particle scattering experiments there has been considerable interest in manifestly scale-invariant theories.
\citep[1643]{Canuto_ea}
\end{quote}
With the remark on ``considerable interest in manifestly scale-invariant theories'' in high energy physics the authors  referred  to Bjorken scaling and, in particular, the seminal paper  \citep{Callan/Coleman/Jackiw}.
But the authors were carefull not to claim field theoretic reality for Dirac's scalar function  $\beta$  \citep[1645]{Canuto_ea}. They rather  developed model consequences for the approach in several directions:  cosmology, including  ``LNH (large number hypothesis as a gauge condition'', modification of the Schwarzschild solution in the Dirac framework, consequences for planetary motion, and stellar structure. At the end the authors indicated certain heuristic links to gauge fields in high energy physics of the late 1970s.

Canuto was interested in exploring Dirac's idea that, perhaps, the gravitational units of measurement, expressed by a locally dependent parameter of gravity (in place of  a constant),  and a frequency change of gravitational clocks, like the period of planets revolving a star, might differ from the atomic units. This would imply a violation of the strong equivalence principle.  In careful evaluations of the astrophysical data available at the beginning of the 1980s he and his coauthor {\em Itzhak Goldman} concluded that a tiny difference might still be possible \citep{Canuto/Goldman:1983}.

For some years Dirac's approach attracted also some interest from astronomers at the Geneva observatory, {\em Pierre Bouvier, 
 Andr\`e Maeder} and coworkers.\footnote{For Canuto and Maeder compare \citep[pp. 126ff.]{Kragh:varying_c}}  
In November 1977, only a few months after the publication of \citep{Canuto_ea}, the two Geneva astronomers submitted a theoretical vindication of ``Weyl's geometry as a framework for gravitation'' to the journal {\em Astrophysics and Space Science} \citep{Bouvier/Maeder:WeylGeometry}. This paper was  meant as a background for a larger research program.  Maeder intended to ``build some new mechanics'' on Dirac-Weyl geometrical gravity. He conjectured that the determination of gravitating mass in gravitationally bound large systems (clusters, super clusters) on the basis of the virial theorem was  affected by the ``new mechanics'' that the {\em missing mass} identified observationally around the middle of the 1970s by astronomers and astrophysicists, might vanish (ibid, 341f.).\footnote{For an illuminating historical reports  on the rise of dark matter see  \citep{Sanders:DarkMatter}. 
From a methodological point of view Maeder's  hypothesis was not so far away from the later, more pragmatic and more successful approach of  modified Newtonian dynamics, MOND, by  {\em Mordechai Milgrom}. }
 First empirical investigations on the Coma cluster seemed to support Maeder's conjecture \citep{Maeder:CosmologyII,Maeder:LowVelocity}. But during the next  years the evidence in favour of his conjecture dissolved.  So the first attempt to bring Dirac's theory to bear in observational cosmology  faded out at the turn to the 1980s. We come back to this issue in section \ref{subsection dark matter}.

But theory development  has an open  horizon.   Dirac's program continued to be pursued during the following decades on the theoretical level among others by {\em Nathan Rosen} working during this time at the Technion Haifa and the University of Beer Sheva and {\em Mark Israelit},  who immigrated to Israel in 1971 and acquired his PhD at Haifa in 1975.  In their continuation of the Dirac program, Rosen and Israelit  sticked as far as possible  to the {\em  e.m.}  dogma
  for a  non-conformally coupled scalar field, $k\neq 6$, but with a  light massive (Proca-type) photon. But already in his  1982 paper Rosen discussed the possibility   of interpreting $\varphi_{\mu}$ as the potential of a  new, hypothetical,  heavy  massive boson field (see below).
 During the  1990s he and Israelit shifted to the last interpretation as the preferred physical view of the Weylian scale connection. 

Rosen extended Dirac's  approach in several respects. He added a scale invariant mass term $\mathfrak{L}_m$ to  the Lagrangian, studied  the dynamical equations, the corresponding currents, the Noether relations, and revisited  the question of different gauges   \citep{Rosen:1982}. Although he recognized the importance of the scale covariant derivatives corresponding to our (\ref{scale covariant derivative}) for giving the Lagrange density a scale invariant form, he {\em did  not } write the dynamical equations  {\em scale invariantly}. The  left hand side of the Einstein equation, e.g., appeared with the Einstein tensor of the Riemannian component of the Weyl  metric,   $_g\hspace{-0.05em}G =\, _g\hspace{-0.1em}Riem- \frac{_g\hspace{-0.08em} R}{2} g$ in the notation of our section \ref{subsection Weyl geometry}, rather than with the   respective (scale invariant) Weyl geometric tensors $G= Riem - \frac{R}{2}g$. Similarly the right hand side expressions  for the energy-momentum of mass and the scalar field  were neither scale covariant nor scale invariant. All terms of the dynamical equations were stripped down to their Riemannian cores. This deprived the Weyl geometric framework of much of its conceptual strenght, even though the equalities  were  valid in every scale gauge \citep[equ. (121)]{Rosen:1982}. This remained so in all of  his and Israelit's work. A  scale co/invariant form  for the dynamical equations  was introduced only a decade later in the work of Hung Cheng and  Drechsler/Tann  (section \ref{subsection Drechsler/Tann}), .

 Rosen also  posed the question how  Dirac's  ``atomic gauge'', in the sense of Weyl gauge, might be made consistent with a non-integrable Weyl geometric structure in order to remove the old problem which Einstein had raised in 1918 as an objection against Weyl's generalized geometry. He tried to back the ``atomic gauge'' by introducing what he called a ``standard vector'' (field). For any timelike vector field $u$ of Weyl weight $w(u)=-1$  the norm   $|u|=g(u,u)^{\frac{1}{2}}$ is  scale invariant ($w|u|=2-1-1=0$).  If $|u|$ is scale covariantly constant, i.e. $D|u| =0$ ($D$ the scale covariant derivative),  Rosen called it  a {\em standard vector} field  and considered the hypothesis that atoms carry a  ``standard vector'' field with them  \citep[p. 220f.]{Rosen:1982}. But he was cautious enough  not to declare this hypothesis as a definitive  solution of the measurement problem in Weyl geometric gravity.

He  found that the Noether relations due to the diffeomorphism invariance of the Lagrangian imply the equations of motion for  matter, while the Noether relations induced by its scale invariance show that the scalar field equation is a consequence of the Einstein equation and the generalized ``electrodynamical'' (i.e, scale curvature) equation  \citep[p. 230]{Rosen:1982}. Moreover, studying Dirac's Lagrangian (\ref{scale invariant Dirac action}) with general coefficient $k$, he realized that for the case of non-conformal coupling the scale curvature equation acquires the form of a generalized Proca equation
\beq \nabla_{\nu} f^{\mu \nu} + m^2 \varphi^{\mu}=0  \label{Rosens Proca equ}
\eeq 
with  $(m\, [ c \hbar])^2=\frac{1}{2}(6-k)$ \citep[233]{Rosen:1982}. This was consistent with Smolin's observation regarding the Weylian scale connection (cf. section \ref{subsection SM 1970s}), which Rosen apparently did not know. He concluded that in the case $k\neq 6$ two physical interpretations for the scale connection were possible:  $\varphi_{\mu}$ might represent an electromagnetic field with massive photons of very small mass, or  a ``meson'' field  extremely weak interacting with ordinary matter.\footnote{Rosen's ``meson'' was a hypothetical massive fundamental boson, no bound state of quarks like the ones of the SM.}
 He added:
\begin{quote}
These mesons could conceivably
accumulate at the center of galaxies and galaxy clusters and could provided (sic!)
the ``missing mass'' that is needed to give a closed universe. \citep[p. 234]{Rosen:1982}
\end{quote}
Rosen thus considered an early ``dark matter'' hypothesis  for the Weyl field, at a time when the conditions for the present understanding of dark matter in galaxies and structure formation was just forming \citep{Sanders:DarkMatter}. He mainly related it  to  the missing mass for cosmological models of positive spatial curvature and alluded at best implicitly to  Zwicky's early observations  of a  mass problem in galaxy clusters.   

 For   cosmological investigations Rosen also considered a vanishing scale curvature. That led to an integrable Weyl geometry with the logarithm of the Dirac scalar field as the potential of the scale connection like in our equ. (\ref{triviality condition})  \citep[equ. (136)]{Rosen:1982}. Although this implies a dynamically trivial extension of Einstein gravity, Rosen found it interesting to discuss scaling effects from a geometrical point of view. 
For a  Robertson-Walker type metric $\overline{g}_{\mu \nu}$ of the form 
\beq   d\overline{s}^2 = dt^2 - \frac{a(t)^2}{a_o^2} dl^2  \,  \label{Robertson Walker metric}  \eeq 
he introduced 
 \[ g_{\mu \nu} = \frac{a_o^2}{a(t)^2} \overline{g}_{\mu \nu} \, .  \]
as the  {\em cosmic gauge}. After an appropriate  reparametrization of  the time coordinate  this led to a   static Riemannian metric for  the model (\ref{Robertson Walker metric}),
\[ g_{\mu \nu}: \quad ds^2 =  dT^2 - dl^2\, , \]
and  $\beta =\frac{a(t)}{a_o}$ \citep[234ff.]{Rosen:1982}. 
He showed that   the cosmological redshift $z$ of a light signal, emitted at time $T_o$  and received at $T_1$, remains invariant under rescaling. In the ``cosmic gauge''  it appears no longer due to a spatial expansion of the geometry, but to the scalar field $\beta$, with $z+1 = \frac{\beta{T_1}}{\beta(T_o)}$, or equivalently, what Rosen did not mention, due to the Weylian scale connection in the ``cosmic gauge'' (cf. subsection \ref{subsection trivial}).

\subsubsection{Omote, Utiyama and the Japanese group\label{subsection Utiyama}}
Already in 1971 and unnoticed by Dirac, a Lagrangian field theory of gravity with a scale covariant scalar field coupling to the Hilbert term like in JBD theory, but
now explicitly formulated in the framework of Weyl geometry had been formulated by {\em M. Omote}, Tokyo, \citep{Omote:1971}.\footnote{This was more than a year before the Trieste symposium at which Dirac talked about his ideas. Apparently the paper remained unknown to Dirac. A second paper by Omote followed after Dirac's publication and after Utiyama had jumped in \citep{Omote:1974}.}
 A little later, and more or less 
at the time of  Dirac's retake of Weyl geometry   {\em Ryoyu Utiyama} (Toyonaka/Osaka)  headed toward a similar goal, although referring to A. Bregman's paper discussed in section \ref{subsection Cartan-Weyl}  and with a main interest in  elementary particle physics.\footnote{In his first paper of 1973 Omote was not mentioned; it was taken up, however, in the references of \citep{,Utiyama:1975a}.}
 Different to Dirac, he left the {\em em} dogma behind and tried to understand the (nontrivial) Weylian scale connection as a new fundamental field.  In a series of papers he 
 ventured toward  its bosonic interpretation   \citep{Utiyama:1973,Utiyama:1975a,Utiyama:1975II} and   presented his results at the Seventh International Conference on Gravitation and Relativity (Tel Aviv, June 1974). Utiyama emphasized that  a Brans-Dicke field $\phi$ of weight $-1$, imported to Weyl geometry, could serve as a  kind of {\em measure field}  (Utiyama's terminology) with respect to which gauge invariant measurable quantities could be expressed starting from any gauge \citep{Utiyama:1973,Utiyama:1975II}.

The import of the scalar field into a Weyl geometric structure would  let it appear natural that $\phi$ is accompanied by a Weylian scale connection $\varphi$ with non-vanishing curvature  (``Weyl's gauge field''). So Utiyama proposed to explore the ordinary Yang-Mills Lagrangian term for a Weylian scale connection
\beq   \mathcal{L}_{\varphi}= - \varepsilon  \frac{1}{4}  f_{\mu \nu } f^{\mu \nu } \sqrt{|det\, g|} \quad (\mbox{here with}\; \varepsilon =1) \eeq
\citep[(2.4)]{Utiyama:1975II}.\footnote{Dirac included a similar scale curvature term in his Lagrangian, but did  not study its consequences.} 
He studied conditions under which ``Weyl's gauge field''  admitted plane wave solutions, and came to the conclusion that it would be  ``tachyonic'', i.e. a  field which allowed superluminal propagation of perturbations.  In Utiyama's  view the   ``boson'' had therefore to be confined to the interior of matter particles. Nevertheless he thought that this ``unusual field $\varphi_{\mu }$ might play some role in establishing a model of a stable elementary particle'' \citep[2089]{Utiyama:1973}.  

Utiyama's results were not  generally accepted. {\em  Kenji Hayashi} and {\em Taichiro Kugo}, two younger colleagues from Tokyo resp. Kyoto, reanalyzed his calculations and argued that, with slight adaptations of the other parameters, the sign $\varepsilon $ could just as well be switched. Then an ordinary, at least non-tachyonic, field would result \citep[340f.]{Hayashi/Kugo:1979}.\footnote{Apparently  Hayashi/Kugo used different signature conventions from Utiyama, which resulted in another sign flip in $\varepsilon $. }
Even then the scale connection would still have strange physical properties. The two physicists showed, after a careful introduction of Weyl geometric spinor fields and their Lagrangians (using scale covariant derivatives), that the  scale connection terms canceled. As they considered only the kinetic term of fermionic Lagrangians, no Yukawa term, in their approach neither the scalar $\phi$-field nor the scale connection $\varphi$  coupled effectively to spinor fields.\footnote{ $\phi$ and  $\varphi$ not even coupled to the electromagnetic field, as Hayashi and Shirafuri showed in another paper the same year.}

At the very moment that a Weylian scale connection $\varphi$ was  interpreted as a ``physical'' field beyond electromagnetism,  it started to puzzle its investigators and, at first, posed more riddles than it was able to solve. It  seemed not to  couple  to matter fields at all (Hayashi/Kugo),  looked either  ``tachyonic'' (Utiyama) or, as we shall see below (Smolin, Nieh, Hung Cheng), appeared to be of Planck mass, far beyond anything observable.

\subsection{\small   Cartan-Weyl geometric  approaches\label{subsection Cartan-Weyl}}
\subsubsection{The Cartan geometric approach to gravity}
Another input to Weyl geometric gravity  came  from the research tradition  started by  {\em Dennis Sciama} and {\em Thomas Kibble}  who   developed a  theory of gravity  by ``localizing''   the symmetries of Minkowski space, i.e., the Poincar\'e group \citep{Sciama:1962,Kibble:1961}.
They  treated the ``external'' symmetries of spacetime  similar   to  the ``internal'' ones
investigated in the   gauge theories over Minkowski space in elementary particle physics  (isospin, $SU_2$, later $SU_3$ and generalizations) which arose from the works of Yang/Mills and Utiyama.   Without explicit recourse to Cartan they reproduced  basic structures of Cartan geometry in field theoretic terms  written in classical tensor calculus.  The   dynamical nature of the infinitesimal translations component of   the ``localized'' Poincar\'e group found its expression in the asymmetry of the linear connection, viz.   torsion. In the sequel different Lagrangians were investigated, and more general groups, in particular the scale extended Poincar\'e group or the affine group, were studied.  In this  way a broad field of  {\em    gauge theories of gravitation} arose \citep{Blagojevic/Hehl}.

 During the 1970s  several  authors introduced Cartan geometric  methods into this research program, particularly prolific among them {\em Andr\'e Trautman}, Warszaw,  and {\em Friedrich Hehl},  Cologne.\footnote{\citep{Hehl:1970Habil,Trautman:1973,Hehl_ea:1976}. For ex-post surveys see  \citep{Trautman:2006,Hehl:Dennis} and the  rich reader     \citep{Blagojevic/Hehl}. } 
They showed that Cartan geometry offered a tailor-made geometric framework for  infinitesimalizing  (``localizing'' in the language of physicists) the symmetries and the  currents known from Minkowski space and special relativity. About the same time also the first publications studying the scale extended Poincar\'e group, often called the {\em Weyl group}, 
\beq 
 \quad \mathfrak{W} = \R^n \rtimes (SO(1,n-1) \times \R^+) \, ,
\eeq 
appeared. 
In the global view $(\R^+, \cdot)$ operates as the  {\em dilation} group on the translations and, in case of a global view, on the underlying Minkowski space $\R^{(1,3)} \cong  \mathfrak{W}/ (SO(1,3) \times \R^+)$ in the case $n=4$. Under localization, or equivalently in the corresponding Cartan space, the  infinitesimal groups (Lie algebras)  are related in such a way that $so(1,3) \oplus \R$ operates on the infinitesimal translations, $\R^4$.  $\R^4$ is ``soldered'' point dependently to the tangent spaces of the underlying differentiable space $M$ by  specifying a  tetrad field or more generally a {\em frame field}, i.e., a family of  bases of the tangent spaces.  In more recent mathematical terms  this corresponds to the choice of a  {\em Cartan gauge} in a Cartan space modelled after $\mathfrak{W}/ (SO(1,3) \times \R^+)$, respectively the corresponding Lie algebras \citep{Sharpe:Cartan_Spaces}. What appears in the global view as an operation on the space itself was thus reshaped,  in the infinitesimalized situation, as a mere change of a Cartan gauge.  
Weyl's intentions of his 1918 geometry and the ideas of   Dicke and Dirac regarding unit scaling were  well expressed in this approach, and at the same time extended by introducing translational curvature, {\em torsion} in Cartan's terminology.

\subsubsection{Alexander Bregmann,} \ldots
 at that time working at Kyoto, inferred  from  \citep{Omote:1971}  that  localized rescaling  could be  separated from Weyl's geometrical interpretation of the infinitesimal length transport. He argued that the point-dependent scale transformations  could be treated ``analogous to the introduction of a space-time dependence into the constant parameters of Isospin or Poincar\'e transformations''. The global scale dimensions  $d$ of a physical field $X$ could then be  taken over as ``Weyl weight'' (Bregman's terminology) of $X$ to the localized theory   \citep[p. 668]{Bregman:1973}.  He first  developed a  Kibble-like approach to gravity built upon the  Poincar\'e group with  tetrad fields $h^{\mu}_a, \; (a=0,\ldots,3$  indexing the tetrads, $\mu = 0, \ldots 3$ the coordinates). He introduced    a  covariant derivative in terms of tetrad coordinates allowing for torsion, and a spin connection expressed  by coefficient systems of the form $A^{m n}_{\mu}$ with regard to generators $S_{mn}$ of the Lorentz group.\footnote{Notations  have been slightly adapted.} 
Then he went on ``to accommodate''  the Weylian scale  transformations to the tetrad calculus, in particular rescaling the tetrads with weight $-1$ (ibid. pp. 675ff.)
\beq  h^{\mu}_{\;\;  a} \mapsto   h'^{\mu}_{\;\; a} =   \Omega^{-1} h^{\mu}_{\;\; a} \label{Bregman tetrad rescaling}
\eeq 
This expressed an operation of the scale group on the tetrads, not on the tangent vectors which remained unaffected by rescaling.\footnote{For $g_{\mu \nu} = h_{\mu}^{\;\;a} h_{\nu a}$ the convention (\ref{Bregman tetrad rescaling}) boils down to $g_{\mu \nu} \mapsto g'_{\mu \nu} = \Omega^2 g_{\mu \nu}$.} 
That  made it necessary to extend the spin connection by a component in the Lie algebra $R$ of the scale group, i.e., a Weylian scale connection  $\varphi = \varphi_{\mu}dx^{\mu}$,\footnote{Bregman, like many other of our authors, used a sign inverted convention for the scale connection form. } 
\beq \hat{A}^{m n}_{\mu} = A^{m n}_{\mu} - (h_{\mu}^{\;\;m} h^{\nu  n} - h_{\mu}^{\;\;n} h^{\nu m})          \varphi_{\nu} \qquad \mbox{\citep[equ. (3.6)]{Bregman:1973}.}
\eeq 
Bregman  remarked that the  modified spin connection represented by $\hat{A}^{m n}_{\mu}$ is  ``Weyl invariant'' and used it to define an associated  scale covariant derivative  $\hat{D}_{k}=h^{\mu}_{\;\;k}\hat{D}_{\mu}$ with the   property $\hat{D}_{\lambda}g_{\mu \nu}=-\varphi_{\lambda}g_{\mu \nu}$, typical for a Weylian metric like in our equ. (\ref{metric compatibility}). The corresponding  linear connection $\hat{\Gamma}^{\lambda}_{\mu \nu} = \Gamma^{\lambda}_{\mu \nu} + \delta^{\lambda}_{\mu}\varphi_{\nu} + \delta^{\lambda}_{\nu}\varphi_{\mu}-g_{\mu \nu} \varphi^{\lambda} $ generalized the Weylian affine connection but was no longer symmetric; it rather included the scale invariant  {\em torsion} tensor
\beq T^{\lambda}_{\mu \nu}=\hat{\Gamma}^{\lambda}_{\nu \mu}-\hat{\Gamma}^{\lambda}_{\mu \nu}={\Gamma}^{\lambda}_{\nu \mu}-{\Gamma}^{\lambda}_{\mu \nu} \, .\label{Bregman torsion}
\eeq

In retrospect we can see in Bregman's paper  a symbolism for working in  a Cartan space modelled after the homogeneous space $\mathfrak{W}/(SO(1,3)\times \R^+$,  later  called a {\em Cartan-Weyl space} (or the other way round).\footnote{Cf. \citep{Sharpe:Cartan_Spaces}.} 
This terminology was not Bregman's; he used Cartan geometric language rather parsimoniously, only with regard to the underlying Riemann-Cartan  structure, modelled after  $\mathfrak{P}/SO(1,3)$ with $\mathfrak{P}= \R^4 \rtimes (SO(1,3)$ the Poincar\'e group. He did not think in geometric terms about the extension of this structure   by   rescaling  the tetrads, the  physical (spinor, vector etc.) fields, the associated spin connection etc.
	
Bregman was more interested in showing how to form Weyl invariant Lagrangians $L$ from matter Lagrangians $L^M$ (his notation) of scale dimension $-4$ with regard to ``constant parameter scale transformations''. He noticed that  many Lagrange densities studied in field theory are also invariant under all conformal transformations of the Minkowski space, including the special conformal ones  (``in particular this is generally true of theories whose quantized versions are renormalizable'') and added:
\begin{quote} In our case such a wider invariance of $L^M$ implies in turn that
the Poincare gauge invariant lagrangian $L^P$ is already Weyl invariant with $L=L^P$ \citep[p. 678]{Bregman:1973}.
\end{quote}
This sharp minded remark generalized  Pauli's observation  that a massless Dirac-spinor field is invariant under Weyl transformations without assuming a coupling to a Weylian scale connection \citep{Pauli:1940}. 

Finally Bregman  gave a short discussion of an integrable scale connection with potential $\sigma$, $\varphi_{\mu}=- \partial_{\mu} \sigma$. He considered $\sigma$ as an  ``independent dynamical variable''  which  is  ``connected to the translation or spin gauge fields'' only through the field equations (ibid. p. 687). This approach facilitated the building of Weyl invariant Lagrange densities. As an example he presented a Lagrangian, which was similar to Omote's  Lagrangian (and Dirac's   not yet published one) including  an additional torsion term  \citep[equ. (5.2)]{Bregman:1973}.\footnote{The torsion term $\frac{2}{3} T^{\alpha}_{\mu \alpha} \partial^{\mu}\sigma^2$ in Bregman's equation is not scale invariant for itself, but his whole Lagrangian density is. }
  All in all,  this  was a remarkable paper which seems to have been underestimated in the following development. 
	
	\subsubsection{Charap/Tait}
	About a year later {\em John Charap} and {\em W. Tait},  London, presented a ``gauge theory of the Weyl group'' building upon the papers  \citep{Yang/Mills:1954,Utiyama:1956,Kibble:1961,Dirac:1973}, while  Bregman's  paper  remained unnoticed by them. They introduced the Weyl group as the ``simplest possible non-trivial enlargement of the Poincar\'e group'' \citep[p. 250]{Charap/Tait}, where by ``non-trivial''  they apparently hinted at  the semidirect product operation of $\R^+$ on the translations. Like Bregman before them they explored ``the consequences of demanding for a theory of matter fields that it be invariant under the transformations of the Weyl group'' (ibid.). 
	
They started by studying the infinitesimal Weyl transformations on  Minkowski space endowed with a globally Weyl invariant Lagrangian $L(\chi, \chi')$ depending on  a couple of fields $\chi$ and their first derivatives (indicated by $ \chi'$). They  derived the Noether relations with regard to translations, rotations and dilations without mentioning Noether.\footnote{This was characteristic for the time. Over several decades the knowledge about the  invariance properties of Lagrangian field theories and the know-how of dealing with it spread  with marginal or no reference at all to Noether's seminal paper \citep{Noether:1918}, cf. \citep{Kosmann:Noether}.} If the Euler-Lagrange equations for all fields are satisfied (``on shell'')  conservation laws for expressions corresponding to the symmetries follow (Noether's first theorem).  Most of them could  easily be identified with well known physical quantities and were  called   {\em  canonical currents}, : the { canonical energy momentum} current $_c\hspace{-0.15em}T^{\mu}_{\nu}$, the { canonical  angular momentum} current $M^{\mu}_{\nu \lambda}$  and additionally a  {\em canonical dilation current} $_c\hspace{-0.1em}\Delta^{\mu}$. The latter evaded an immediate physical interpretation. But in analogy  to the angular momentum, which can be decomposed into an internal (spin) contribution of the fields and an external, orbital component, the dilation current could be decomposed into
\beq  _c\hspace{-0.1em}\Delta{\mu} = J^{\mu} + \, _c\hspace{-0.15em}T^{\mu}_{\nu} x^{\nu}\, , \qquad  J^{\mu} = \frac{\partial L}{\partial (\partial_{\mu}\chi)} w(\chi) \chi\, ,
\eeq 
 with an internal component $J^{\mu}$ (summation of the field components understood to be included) and an external one depending on the origin of the coordinate system in Minkowski space. The external component was, of course, due to the dilational operation of $\R^+$ on the underlying space. It will be interesting to see what became of these under ``localization'' of the symmetries.

In the next step Charap and Tait localized the approach following Kibble's path; i.e., they made the group operations point dependent by introducing a tetrad field $h^{\mu}_k$ and  {\em gauge fields} (geometrically spoken connections), the first one with values in the Lorentz algebra, given by coefficients $A^{ij}_{\mu}$ with respect to the  algebra generators $S_{ij}$, and another one with values   in the real numbers, given by   a system of $A_{\mu}$. This allowed to define the  scale covariant derivative 
\beq D_{\mu} \chi = \partial_{\mu}\chi + \frac{1}{2}A^{ij}_{\mu}S_{ij}\,\chi - A_{\mu} w(\chi)\, \chi  \qquad \mbox{\citep[equ. (3.7)]{Charap/Tait}.}
\eeq 
Like Bregman they considered  Lagrange densities $\mathfrak{L}$ constructed from   globally Weyl invariant ones, $L$, by substituting partial derivatives by   scale covariant derivatives  (``minimal coupling'' in the later physics idiom) and using the volume element $|h|^{-1}dx$ with $|h| =det(h^{\mu}_k)$. They analyzed   the resulting gauge field dynamics and considered their  sources, usually the ``right hand side'' of the equations, as ``modified `currents'''\footnote{$ \mathfrak{T}^k_{\mu}= \frac{\partial \mathfrak{L}}{\partial h^{\mu}_k}\, , \quad  \mathfrak{S}^{\mu}_{ij}= -2  \frac{\partial \mathfrak{L}}{\partial A^{i j}_{\mu}}\, , \quad  \mathfrak{D}^{\mu}= \frac{\partial \mathfrak{L}}{\partial A_{\mu}}$ \citep[p. 256, notation slightly changed]{Charap/Tait}. Compare the {\em dynamical currents} by Hehl et al. below (\ref{hypermomentum current}). For the dynamical matter energy current variational and partial derivatives usually coincide.  Charap and Tait may have (wrongly) generalized this property to the other currents.}

Among the localized {\em  canonical}  (Noether) currents only the one referring to {\em energy-momentum } was explicitly mentioned by them.   In analogy to the global case they defined
\beq \,   _c\hspace{-0.1em}\mathfrak{T}^{\lambda}_{\mu} =  \frac{\partial \mathfrak{L}}{\partial (\partial_{\lambda} \chi )} \partial_{\mu} \chi - \delta^{\lambda}_{\mu}\mathfrak{L}
\eeq  
and commented upon its difference  from the {\em dynamical energy momentum}, arising from variational derivation of the matter Lagrangian with respect to the gravitational field. The latter written with respect to coordinate basis as a ``world tensor'', $\mathfrak{T}^{\lambda}_{\mu}$,  can be derived from the canonical energy momentum $ _c\hspace{-0.1em}\mathfrak{T}^{\lambda}_{\mu}$ by adding terms in  dynamical spin and dilation quantities.\footnote{$\mathfrak{T}^{\lambda}_{\mu}= _c\hspace{-0.12em}\mathfrak{T}^{\lambda}_{\mu}  - \frac{1}{2}\mathfrak{S}^{\lambda}_{ij}A^{ij}_{\mu} -   \mathfrak{D}^{\lambda}A_{\mu}$ \citep[equ. (3.22)]{Charap/Tait}.}
  But neither  the corresponding conservation theorem nor the other Noether currents were mentioned, not even the canonical dilation current which was round the corner,
\beq 
_c\hspace{-0.1em}\Delta^{\mu}  = - \frac{\partial \mathfrak{L}}{\partial (\partial_{\mu} \chi )} w(\chi)\chi\, . \label{dilational Noether current}
\eeq 
 In a separate passage on the  ``geometrical interpretation'' of their theory they  explained how  Weyl's scale geometry of 1918 was taken up  by their approach and how it was generalized by including torsion.  Cartan geometry was neither mentioned nor used.

\subsubsection{Hehl et al. at the Kiel conference}
In the next decade several papers on gravity in a Cartan-Weyl   environment appeared; among them \citep{Kasuya:1975,Obukhov:1982,Nieh:1982}.\footnote{All of them  are mentioned in  \citep[chap. 8]{Blagojevic/Hehl}.}  Not all of these contained new insights, nor were their statements always reliable; but they indicate a slow broadening of an  interest in Weyl geometric gravity with torsion.  This interest acquired   a wider  context with the rise of  Cartan geometric {\em metric-affine} studies of gravity. 
\label{metric-affine}
A general discussion of this research would go far beyond the limits of this survey.\footnote{In particular F. Hehl and his varying  coauthors were active in this field, \citep{Hehl_ea:1976short,Hehl_ea:1988Weyl_group,Hehl_ea:1988Kiel,Hehl_ea:1988Skaleninvarianz,Hehl_ea:1989local_scale,Hehl_ea:1995Report}. More papers by other authors, while only a few of those just mentioned, appear  in \citep[chap. 9]{Blagojevic/Hehl}.}
But at least one of the papers of this field has  to be commented upon here  \citep{Hehl_ea:1988Kiel}. This paper was presented by {\em Friedrich Hehl} and coauthors at the Weyl centenary conference at Kiel in 1985. It discussed a view of  Cartan geometric metric-affine  of gravity with  particular emphasis on the   Weyl group $\mathfrak{W}\subset \R^n  \rtimes GL(n,\R)$ and   may be taken as paradigmatic for how Weyl geometric aspects were dealt with in this research program. 

In modernized language, it worked in a Cartan geometry modelled after the Klein space $\mathfrak{A}/GL(n,\R)$ with the  affine group $\mathfrak{A}= \R^n  \rtimes GL(n,\R)$, respectively the corresponding Lie algebras. Additionally an  independently given (Lorentzian) metric $g$ was assumed. The local description involved again an $n$-frame field $h^{\mu}_a\, ({ a, \mu = 0, \ldots, n-1})$  and its dual system of forms $h_{\mu}^a$ characterizing a Cartan gauge or the translational connection. Moreover, a  connection with values in the Lie algebra of the isotropy group $gl(n,\R)$  (generalizing the rotations of Cartan-Riemann geometry) was given by the coefficient system $A_{\;a\, \mu}^{b}$,  and a metric of Lorentzian signature by  $g_{\mu \nu} = h_{\mu}^a  h_{\nu}^b\, g_{a b}$.  Let the corresponding covariant derivative be denoted by $D_a$, respectively $D_{\mu}$ if transcribed to coordinate indices. Of course, in general the derivative of the metric does  not vanish, $D_{\lambda}g_{\mu \nu} =- Q_{\lambda \mu \nu}$, with  $Q_{\lambda \mu \nu}$  usually called  the {\em non-metricity} of the derivative, which is symmetric in its last two indices. 
With respect to  coordinate bases  the linear connection is $\Gamma_{\mu \nu}^{\lambda}= h_{\mu}^a h^{\lambda}_b  A_{\;a\, \nu}^{b}$.

In the case of a symmetric $\Gamma$  it  can be decomposed,    $\Gamma= _g\hspace{-0.12em}\Gamma+_Q\hspace{-0.12em}\Gamma$,  into its  Riemannian (Levi-Civita) component $_g\hspace{-0.12em}\Gamma$ with respect to $g$ and a component due to the non-metricity $_Q\hspace{-0.12em}\Gamma$.\footnote{$  _Q\hspace{-0.12em}\Gamma^{\lambda}_{\mu \nu} = \frac{1}{2}(Q_{\mu \nu}^{\; \; \; \lambda} -  Q_{\; \; \mu \nu}^{\lambda}+ Q_{\nu\; \;  \mu}^{\; \; \lambda})$}
Also $Q_{\lambda \mu \nu}$ can  be  decomposed into a traceless part $/\hspace{-0.65em}Q_{\lambda \mu \nu}$, with regard to the last two indices, and and its trace $q_{\lambda}$.  Hehl and coauthors 
introduced what they now introduced as the ``celebrated Weyl vector'',  the trace part $q_{\lambda}$ of the non-metricity:  
\beq 
q_{\lambda}=\frac{1}{ n}Q_{\lambda \nu}^{\;\;\; \nu}\,  \quad \qquad \mbox{\citep[p. 252]{{Hehl_ea:1988Kiel}}}
\eeq
In our notation  $ \varphi_{\lambda} = \frac{1}{2}q_{\lambda} $ is a {\em part of  the connection} and can also be expressed by   $\varphi_{\lambda}=\frac{1}{n}A^a_{\; a\, \lambda} $. 
In the case of vanishing  $/\hspace{-0.65em}Q$  the ``Weyl vector''  satisfies the metric compatibility condition of Weyl geometry, our equ. (\ref{metric compatibility}). Then also the linear connection $\Gamma$  coincides with the Weylian invariant affine connection. This allowed to embed Weyl geometry in the wider framework of metric-affine geometry.

In this framework the  data $(h^{\mu}_a, A_{\;\;a\, \mu}^{b},g_{a b})$ were considered as   dynamically independent  field components of an extended theory of gravity \citep[sec. 6,  notation changed]{Hehl_ea:1988Kiel}. 
For  a set of  matter fields $\Psi$ and a matter Lagrangian $\mathfrak{L}_m$ minimally coupled to the affine-metric structure $\mathfrak{L}_m(g_{ab},h^{\mu}_a,\Psi, D_a \Psi )$ the authors defined the {\em dynamical}  matter {\em hypermomentum current}  of their generalized theory  by the variational derivative with regard to the full isotropy connection,
\beq  \Delta_{b}^{\;\; a\, \mu} = |h| \frac{\delta \mathfrak{L}_m}{\delta A_{\;a\, \mu}^{b}} \, , \qquad \mbox{with} \quad |h|=det(h^{\mu}_a) \, , \label{hypermomentum current}
\eeq 
 and  decomposed it into its rotational, dilational and shear  components \citep[p. 274f.]{Hehl_ea:1988Kiel}.  

In  the context of the talk, the dynamical {\em dilational  current} 
\beq
\Delta^{\mu} :=\Delta_{a}^{\; a\, \mu}= |h|  \frac{\delta \mathfrak{L}_m}{\delta A_{\;a\, \mu}^{a}} =\frac{|h|}{n} \frac{\delta \mathfrak{L}_m}{\delta \varphi_{\mu}} \label{Hehl ea dilational matter current}
\eeq
 was  of particular importance. The authors conceded  that the shear current was  ``remote from physical experience''; but   for the dilational current  they saw   ``supporting evidence''  in Bjorken scaling of deep inelastic scattering experiments and, on the theoretical level, in certain models of  supergravity \citep[pp. 244, 275]{Hehl_ea:1988Kiel}. This was a central  point for their talk. Aready in their abstract they  announced:
\begin{quote} 
In the light of modern developments in particle physics, this coupling of the Weyl vector to the material dilation current is an unalterable part in any viable theory of a general-relativistic type, which comprises a Weylian piece \citep[p. 241]{Hehl_ea:1988Kiel}.
\end{quote}
This  thesis stood somehow in contrast to the observation of Bregman (who was cited by our authors) that for  matter fields with conformally invariant Lagrangians there was no need for  assuming their  coupling to  the scale connection (the ``Weyl  vector''). But the authors gave an example of a scalar field with non-vanishing dilational current (see below).

Hehl and his coauthors analyzed the dynamics of the metric affine theory on a quite general level with, at first, no particular  Lagrangian specified. They described the form of the three dynamical equations corresponding to the decomposition of the general linear group into shear, rotational, and dilational components and argued that only two of them were dynamically independent because of the interdependencies due to the  relations of the second Noether theorem \citep[sec. 8]{Hehl_ea:1988Kiel}. A short discussion of specific Lagrangians followed, among them some consequences of a scale covariant  ``primordial scalar field, the so-called {\em dilaton field} $\sigma(x)$''  \citep[p. 282, emph. in original]{Hehl_ea:1988Kiel} of Weyl weight $w(\sigma)=-\frac{n-2}{2}$ and  the usual scale invariant quadratic kinetic term 
\beq 
\mathfrak{L}_{\sigma} = \frac{|h|}{2} D_{\mu}\sigma  D^{\mu} \sigma \, .
 \eeq 
 As  a  ``primordial''  field $\sigma$ was considered to be a part of the matter sector, and because of the  scale covariant derivative in $\mathfrak{L}_{\sigma}$ it  contributed to the  matter dilational current with
\beq 
\Delta^{\mu}(\sigma) = \frac{n-2}{2n}\sigma D^{\mu}\sigma \qquad \mbox{\citep[equ. (10.11)]{Hehl_ea:1988Kiel}.}
\eeq

But the  author team warned that one should not expect  easy empirical repercussions in laboratory experiments:
\begin{quote}
Local scale invariance of fundamental interactions is expected to be valid only approximately in the high energy limit of Bjorken scaling or exactly at the onset of the big bang (ibid, p. 285).
\end{quote}
They  assumed a breaking of scale symmetry down to the Poincar\'e group ``after a very short time lag'' (to the big bang) and proposed  a quartic potential for $\sigma$ with a symmetry breaking quadratic term similar to the Higgs potential. 
In the end, their gravitational Lagrangian {\em  boiled down to Einstein gravity} with cosmological constant ``plus some supplementary terms known from Poincar\'e gauge theory'' (ibid, p. 242). By the ``supplementary terms'' they apparently referred to the torsion-spin coupling which arises in Einstein-Cartan gravity.  It becomes a serious modification of Einstein gravity only at extremely high mass densities.\footnote{According to later estimates  the torsion-spin coupling of Einstein-Cartan gravity becomes important only close to  $10^{38}$ times the density of a neutron star, 
 which signifies energy densities at the hypothetical grand unification scale of elementary particle interactions \citep[p. 194]{Trautman:2006}, \citep[p. 108]{Blagojevic/Hehl}.}

In the framework of Hehl et al., Weyl geometric modifications of gravity were  to  be expected only under even more extreme condition than  for torsion. In  their  view, Weyl geometric effects seemed to be banned to a speculative realm close to the ``big bang'', one of the great adventure playgrounds of late 20th century physics. In the sections \ref{section SM},  \ref{section cosmology} we shall see that this need not necessarily be so, if  other perspectives are taken into account. But before we  turn to these researches, we have to pay attention to another road towards reviving Weyl's scale geometry. It had  different roots from those of  the Omote-Dirac-Utiyama and the Cartan geometric approaches  discussed in this section and arose from an attempt to geometrize the dynamics of non-relativistic quantum mechanics.

\section{\small   Weyl's scale connection a geometrical clue to quantum mechanics?\label{section geometrical quantum mechanics}}
\subsection{\small Bohmian mechanics as a background\label{subsection Bohm}}
\subsubsection{Bohm's ``causal'' approach to QM}
In the early 1950s {\em David Bohm} proposed an alternative approach to  non-relativistic quantum mechanics (QM) with an often discussed heterodox ``causal'' interpretation of the latter.\footnote{\citep{Bohm:1952a,Bohm:1952b}}
 His core idea was to reintroduce exact particle trajectories into the description of quantum systems, which  were guided by  a pilot wave evolving according to the Schr\"odinger equation of ordinary quantum mechanics. In this move he took up earlier ideas of {\em Louis de Broglie} on the dualism of wave and particle aspects in QM, which had been critically debated in the late 1920s.  This approach was mathematically close to a hydrodynamic picture of the Schr\"odinger equation, considered by {\em  Ernst Madelung} in 1926.  Madelung noticed that his ``hydrodynamical'' current was subject to a non-classical term which could be interpreted  as a kind of force function due to the `` `inner forces' of  the continuum'' \citep[p. 323]{Madelung:1926}. Bohm  extended these older ideas, among others, by an analysis of the measuring process.\footnote{Bohm realized the kinship of his approach to the earlier proposals of de Broglie only after he had finished his manuscript \citep[p. 167]{Bohm:1952a}. In a  footnote added in proof he also referred to Madelung's ``similar'' approach of 1926, adding the remark ``\ldots but like de Broglie he did not carry this interpretation to a logical conclusion'' (ibid.).} 
 He  could thus  avoid  to stipulate  a ``collapse'' of the wave function, which was usually assumed  for extracting real valued measuring values from the observables given by the Hermitian operators of QM \citep[chap. 11]{Bacciagaluppi/Valentini}.
Bohm wanted to challenge the mainstream  (``Copenhagen'')  interpretation of QM  which he  accepted as  consistent but  as unsatisfactory from a foundational and natural philosophic point of view. His goal was  to find an {\em alternative interpretation} of QM which did not affect the dynamics, at least not ``in the domain of dimensions of the order of $10^{-13}$ cm''. It ought to permit
\begin{quote} \dots    to conceive
of each individual system as being in a precisely definable state, whose changes with time are determined by definite laws, analogous to (but not identical with)
the classical equations of motion \citep[p. 167]{Bohm:1952a}.
\end{quote}

Bohm started from the observation that to any  Schr\"odinger equation for a  wave function $\psi(x)=a(x)e^{\frac{i}{\hbar}S(x)}$ ($x \in \R^3$) governed,  e.g. in the case of  a single particle of mass $m$ in an external potential $V(x)$, by the equation 
\beq  i \hbar \frac{\partial \psi}{\partial t} = - \frac{\hbar^2}{2m}\, \nabla^2 \psi + V(x) \psi  \, ,\label{Schroedinger equ}
\eeq 
one can associate a pair of coupled differential equations for the  phase $S(x)$ of $\psi(x)$ and  the probability density $\rho(x) = a(x)^2$:
\beqarr
 \qquad \qquad \frac{\partial \rho}{\partial t} + \nabla \cdot \left( \rho \frac{\nabla S}{m} \right)  &=&   0 \, , \label{Bohm continuity equ I} \\
\frac{\partial S}{\partial t} + H(x,t) &=& 0  \, ,     \label{Bohm Hamilton-Jacobi equ} \\
 \mbox{where}  \qquad H(x,t) =   \frac{(\nabla S)^2}{2m}   &+&  V(x) - \frac{\hbar^2}{4 m}\left( \frac{\nabla^2 \rho}{\rho} - \frac{1}{2}\frac{(\nabla \rho)^2}{\rho^2} \right) \,  \nonumber 
\eeqarr
\citep[equs. (6), (7)]{Bohm:1952a}. 
Equation (\ref{Bohm Hamilton-Jacobi equ}) has the form of a {\em Hamilton-Jacobi} equation for a point particle  with the principal function $S$ and  conjugate momenta $p_k = \partial_k S $, equivalently the velocity $v=\frac{\nabla S}{m}$. The total potential $ \tilde{V}(x) = V(x) + U(x) $ deviates  from the classical $V(x)$ by
\beq  U(x)=  - \frac{\hbar^2}{4 m}\left( \frac{\nabla^2 \rho}{\rho} - \frac{1}{2}\frac{(\nabla \rho)^2}{\rho^2} \right)  = - \frac{\hbar^2}{2m} \frac{\nabla^2 a}{a}  \, . \label{Bohm quantum potential}
\eeq 

Bohm considered $U(x)$ as a kind of {\em quantum potential} added to the classical one.
The trajectories of the Hamilton-Jacobi system  have velocities   $v=\frac{\nabla S}{m}$ normal to the level surfaces of constant values of $S$. Thus   (\ref{Bohm continuity equ I})  
acquires the form of a continuity equation  for an ensemble of point particles following the family of trajectories with the density $\rho$,
\beq
 \frac{\partial \rho}{\partial t} + \nabla \cdot \left( \rho v \right)  =   0  \,  . \label{Bohm continuity equ}
\eeq
 He argued that this might be ``the nucleus of an alternative interpretation
for Schroedinger's equation'' \citep[p. 170]{Bohm:1952a}. A pair of equations similar to  (\ref{Bohm Hamilton-Jacobi equ},  \ref{Bohm continuity equ}) had  already been investigated by Madelung   with a hydrodynamical interpretation.   But Madelung was more cautious. He did not consider his equations as  more fundamental than  Schr\"odinger's wave mechanics, but rather  as a ``model representation'' from which one could derive the essential features of the latter \citep{Madelung:1926}.\footnote{``Die hydrodynamischen Gleichungen sind also gleichwertig
mit denen von Schr\"odinger und liefern alles, was jene geben, d. h. sie
sind hinreichend, um die wesentlichen Momente der Quantentheorie der
Atome modellm\"a\ss{}ig darzustellen'' \citep[p. 325]{Madelung:1926}.}

In Bohm's alternative interpretation of the equations  a  quantum particle could seem to be no longer subject to the Heisenberg indeterminacy, because it appeared as though it follows  a specified  trajectory of the system  (\ref{Bohm Hamilton-Jacobi equ}); but the  Heisenberg uncertainty was implicitly preserved because only a probability  satisfying (\ref{Bohm continuity equ})    could be given for a trajectory passing specified regions in the level surfaces. Via  the quantum potential (\ref{Bohm quantum potential})  a wave  function satisfying  (\ref{Schroedinger equ}) would operate as a   non-local guiding structure, a ``pilot wave'' in a terminology not used by Bohm, for the motion of the quantum particle.  This was quite close to de Broglie's theory of the 1920s \citep{Bacciagaluppi/Valentini}.

 Although Bohm   extended de Broglie's and Madelung's view by an analysis of the measuring process, his proposal did not receive  immediate positive response in the quantum physics community \citep{Myvold:Bohm}. Only in the longer run  different authors took  it up and pursued programs along his lines, although sometimes with a different outlook on the underlying  ontology and enriched by new mathematical ideas.
 Independent of differing views on  ontology or mathematical techniques, they belong to  common  family of {\em  de Broglie -- Madelung -- Bohm} ({\em dBMB}) approaches.\footnote{\citep{Passon:2015,Passon:Bohm,Duerr_Ea:Compendium}, I thank O. Passon for his helpful explanations of Bohmian mechanics. For relativistic generalizations, see, e.g., \citep{Nicolic:2005}.}

Important for our context was the stepwise extension of the Bohmian approach to relativistic quantum mechanics, in particular for {\em Klein--Gordon}   particles given by a complex field of spin zero, $\psi(x)=a(x)e^{i\frac{S}{\hbar}}$. It had values in the complex numbers
  like a Schr\"odinger wave function,  but lived on the Minkowski space  with $x=(x^0, \ldots, x^3)$. Moreover, $\psi(x)$ demanded a  more intricate interpretation than Born's probability rule.
In case of an electromagnetic interaction with potential $A_{\mu}$, the wave field of a  Klein-Gordon particle of mass $m$ and charge $e$ satisfies the dynamical equation 
\beq
\left(\frac{\hbar}{i}\partial_{\mu} -\frac{e}{c}A_{\mu} \right) \left(\frac{\hbar}{i}\partial^{\mu} -\frac{e}{c}A^{\mu} \right) \psi = (mc)^2 \psi \label{Klein-Gordon equ}
\eeq 
in the signature $(+---)$ of the Minkowski space ($ \partial_o= c^{-1}\partial_t$). Here the Bohmian method  lead to the Hamilton-Jacobi and continuity equations 
\beqarr
\left( \partial_{\mu}S - -\frac{e}{c}A_{\mu}\right)\, \left( \partial^{\mu}S -\frac{e}{c}A^{\mu} \right) &=& m^2 c^2+ \hbar^2\, Q   \, , \label{Hamilton-Jacobi Klein-Gordon}\\
\partial_{\mu}( a^2\partial^{\mu}(S -\frac{e}{c}A_{\mu})= 0 \, , \label{continuity equ Klein-Gordon}
\eeqarr
with a  ``quantum term'' similar to (\ref{Bohm quantum potential}).\footnote{Cf. \citep[p. 554 ]{Nicolic:2005} for vanishing {\em em} potential. }
\beq 
 Q= \frac{\partial_{\mu}\partial^{\mu} a}{a}\, . \label{quantum potential Klein-Gordon simple}
\eeq 

\subsubsection{A geometrization idea by de Broglie\label{subsection de Broglie}}
In his later years de Broglie himself joined in the renewed research program. Among others, he  pondered about  a connection between the  Jacobi flow of the Klein-Gordon field and general relativity  \citep[pp. 118ff.]{deBroglie:1960}. He defended a hypothesis according to which  quantum particles are  constituted by  extremely dense tiny regions of the field, governed by unknown non-linear equations, while  in the exterior of these  regions the known linear equations of quantum mechanics hold. He called  these regions ``singular'' and investigated whether the motion of such  ``singular'' regions may follow geodesics similar to the motion of  singular regions in general relativity, which had been  studied by Einstein and Grommer  in the 1920s.\footnote{According to de Broglie it was J.-P. Vigier who made him awar of a parallel between his hypothesis and the work of Einstein and Grommer \citep[p. 92]{deBroglie:1960}.}
 In this context he considered the right hand side of (\ref{Hamilton-Jacobi Klein-Gordon}) as a kind of ``variable rest mass'' which had to be  calculated in the ``immediate vicinity of the particle'' \citep[p. 116]{deBroglie:1960}:
\beq \mathfrak{M}_o= \sqrt{m_o^2 + \frac{\hbar^2}{c^2} \frac{\partial{\nu}\partial^{\nu}a}{a}}
\eeq 
Some authors would later  call $\mathfrak{M}_o$ the  ``quantum mass'' of a  Klein-Gordon field (see subsection \ref{subsection Shojai ea}).

De Broglie considered the trajectories  $x^{\mu}(s)$ of the Hamilton-Jacobi flow of a particle ``{\em in the absence of electromagnetic and gravitational fields}'' (ibid. p. 119, emph. in original) with  4-velocity,   $u^{\mu}=\frac{d}{ds}x^{\mu}$, normalized to $u_{\mu}u^{\mu}=1$,
\beq 
u^{\mu} = (\mathfrak{M}_o)^{-1} \partial^{\mu} S \, .
\eeq
and found that they satisfy the geodesic equations of a  metric $g_{\mu \nu}$ arising from the Minkowski metric $\eta_{\mu \nu}$ by conformal rescaling
\beq 
g_{\mu \nu} = \frac{\mathfrak{M}_o^2}{m_o^2}\,\eta_{\mu \nu} \label{de Broglie quantum metric}
\eeq
He concluded:
\begin{quote}
Thus, even if the particle is not subjected to any gravitational or electromagnetic field, its possible trajectories (\ldots) are the same as if space-time possessed non-Euclidean metrics defined by [$g_{\mu\nu}$]  \citep[p. 120]{deBroglie:1960}.
\end{quote}

This was an interesting geometrization argument and de Broglie did not remain the only one to ponder on  a connection between the dBMB quantum mechanics and general relativity. 
Here we are mainly interested in later authors who tried to make progress by attempting  a geometrization  of QM in the framework of  Weyl geometry.

\subsection{\small Santamato's proposal for geometrizing  quantum mechanics\label{subsection Santamato}}
\subsubsection{Two phases of work on the program\label{Santamato two phases}}
In the 1980s {\em Enrico Santamato}, Napoli, proposed a new approach to quantum mechanics \citep{Santamato:1984a,Santamato:1984b,Santamato:1985}. It was based on studying   weak random processes of  ensembles of point particles  moving in a Weyl geometrically modified configuration space. 
He compared his approach with that of   {\em Madelung-Bohm}  and   the stochastic program of {\em Feyn\`es-Nelson}.\footnote{For E. Nelson's program  to re-derive the  quantum dynamics from classical stochastic processes and classical probability see \citep{Bacciagaluppi:Nelson}.}
While the latter dealt with stochastic (Brownian) processes,  Santamato's approach  was closer to the view of  Madelung and Bohm  because it assumed only random initial conditions, with classical trajectories  given in Hamilton-Jacobi form (this explains the attribute ``weak'' above). One can read Bohm's particle trajectories as deviating  from those  expected in Newtonian mechanics by some ``quantum force''. Santamato found this an intriguing idea but deplored the latter's ``mysterious nature'' which   ``prevents carrying out a natural and acceptable theory along this line''. He  hoped to find a {\em rational explanation}  for the effects of the ``quantum force'' by  {\em geometry} with a modified affine connection of the  system's configuration space. Then the deviation from classical mechanics would appear as the outcome of  ``fundamental properties of space''  \citep[p. 216]{Santamato:1984a}, which has to be understood in the sense of {\em configuration} space, as we may add. 

In his first paper paper Santamato started from a configuration space with coordinates $(q^1, \ldots, q^n)$  endowed with a Euclidean metric. More generally,  his approach allowed for a general positive definite metric $g_{ij}$,  and later even a metric of indefinite signature,  for dealing with general coordinates of $n$-particle systems and perhaps, in a further extension, with spin. The Lagrangian of the system, and the corresponding Hamilton-Jacobi equation, contained  the metric explicitly or implicitly.  This Euclidean, or more generally Riemannian, basic structure was complemented by a Weylian scale connection. Santamato's central idea was that the modification of the Hamilton-Jacobi equation induced by a properly determined {\em scale connection} can be used to  express the {\em  quantum modification}  of the classical Hamiltonian like in the  Madelung-Bohm approach. 
  Then the quantum aspects of the systems would be geometrized in terms of the Weyl geometry, surely a striking and even beautiful idea, if it works. 

Santamato thus headed towards a  new program of   {\em geometrical quantization sui generis}. It   had nothing to do with the better known geometric quantization program initiated more than a decade earlier by  J.-N. Souriau, B. Kostant and others, which was already well under way in the 1980s \citep{Souriau:1966,Kostant:1970,Simms:1978}. In the latter  geometrical methods underlying  the canonical quantization were studied. Starting from a {\em symplectic  phase space} manifold of a classical system, the   observables were ``prequantized'' in a Hermitian line bundle, and  finally the Hilbert space representation of quantum mechanics was constructed on this basis.\footnote{See, e.g., \citep{Woodhouse:GeoQuant} or \citep[chaps. 22/23]{Hall:QuantumTheory}.
A classical monograph on the symplectic approach to {\em classical} mechanics is \citep{Abraham/Marsden};  but it does not discuss spinning particles.  In the 1980s the symplectic approach was already used as a starting platform for (pre-)quantization to which proper quantization procedures  could then hook up, see e.g. \citep{Sniatycki}.  Souriau was an early advocate of  this program. In his  book he   discussed relativistic particles with spin \citep[\S 14]{Souriau:1970}.} 
Santamato's geometrization was  built upon a different
 structure, Weyl geometry rather than symplectic geometry, and he had rather different goals.

Like other proposals in the dBMB  (de Broglie-Madelung-Bohm) family, Santamato's program did not find immediate positive response.  In the following decades he 
  shifted the center of his research to  nonlinear optics of liquid crystals and to quantum optics, even with a strong empirical component, and stopped publishing on the foundational topic.
 Perhaps a critical paper by 
 {\em Carlos Castro Perelman}, a younger colleague who knew the program nearly from its beginnings, contributed to the extended period of  interruption? Castro discussed ``a series of technical points'' which seemed important for Santamato's program from the physical point of view \citep[p. 872]{Castro:1992}.\footnote{Carlos Castro later added his mother's name Perelman to his second name. Under this name he is mentioned in the acknowledgements of \citep{Santamato:1984b}. At that time he was research assistant at the University of Texas, Austin, where he acquired his Ph.D. in 1991.} 
Among the problems he mentioned were several of a more foundational than of purely technical  import: {\em (i)} the problem of specifying the random intial data for the ensemble of particle paths, {\em (ii)} the Hilbert space interpretation of the theory, {\em (iii)}  the relationship of Santamato's approach  to the Feynman path integral quantization. He also criticized the lack  of a rigorous hypothesis in the choice of the particle's Lagrangian,  the non-definite character of the probability density in the case of a Klein-Gordon particle  (which could appear if  the  foliation with respect to the principal function $S$ is not timelike),   the un-understood  dependence of the particle's effective mass on the Weylian scalar curvature (in the configuration space), and some other more technical points.

After the turn to the new century/millennium Santamato came back to foundational questions in  close cooperation with his  colleague {\em Franceso De Martini} from the University of Rome. Both had  cooperated in quantum optics already  for many years. In the 2010s they turned  to geometrical quantization  in a series of joint publications and continued the program started by Santamato  three decades earlier. They  showed how to deal with spinor fields in this framework, in particular with the Dirac equation \citep{DeMartini/Santamato:2013} and discussed the famous {\em Einstein-Podolsky-Rosen} (EPR) non-locality question \citep{DeMartini/Santamato:2014a,DeMartini/Santamato:2014b}. Moreover, they analyzed the helicity of elementary particles and showed that the {\em  spin-statistics relationship} of relativistic quantum mechanics can be derived in their framework without invoking arguments from quantum field theory \citep{DeMartini/Santamato:2015a,DeMartini/Santamato:2015b,DeMartini/Santamato:2016a}.  In this new series of papers Minkowski space formed the starting point for the construction  of the configuration spaces which could be extended by internal degrees of freedom. Moreover,  a  transition from point dynamical Lagrangians as the dominant view to a dynamically equivalent  description in terms of  scale invariant  field  theoretic Lagrangians in two  scalar fields enlarged the perspective (see below).   

It is unnecessary to go into details of these often quite technical articles; here I concentrate on the basic question of the geometrization program.\footnote{Santamato and De Martini propagated the geometrical side of their research under different labels: at first they talked about ``affine quantum mechanics (AQM)'' in ``conformal differential geometry'' \citep{DeMartini/Santamato:2013},  then they shifted to ``conformal quantum geometrodynamics (CQG)'' \citep{DeMartini/Santamato:2014a,DeMartini/Santamato:2014b}. }
Here we want to  see how  Santamato's  
intriguing idea of introducing a Weyl geometric structure on the configuration space, in order to model the Bohmian effects of quantum systems, works.

\subsubsection{The geometrization of the configuration space\label{subsection configuration space}}
Summing up, Santamato's idea was to consider dynamical systems with finite degrees of freedom, parametrized by a configuration space $V$ with parameters $q^1, \ldots , q^n$, endowed with a pseudo-Riemannian metric $g_{ij}$ which could be of  any signature \citep{Santamato:1984a}. In the case of a non-relativistic  $k$-particle system without inner degrees of freedom it could be  the product of Euclidean 3-metrics, for relativistic particles in Minkowski space with metric $\eta  = diag(-1,1,1,1)$  it was of signature $(3k,3)$ \citep{Santamato:1984b}. 

In the case of a relativistic  1-particle system with spin the product of the Minkowski space $\M$ and the Lorentz group served as configuration space, $V=\M \times SO(3,1)$, where the second factor  parametrizes   ``hidden'' rotational degrees of freedom of the particle  \citep[p. 634]{DeMartini/Santamato:2013}. 
By an astute choice of coordinates  $(q^1, \ldots q^{10}) =(x^{\mu}, \theta^{\alpha})$ (with $\alpha = 1, \ldots, 6$) in  $V$, with generalized ``Euler angles'' $ \theta^{\alpha}$ for parametrizing $ SO(3,1) $, 
the authors  introduced a metric  $(g_{ij})$  by a  block matrix composed of the Minkowski metric $\eta_{\mu \nu}$ and  a ``metric of the parameter space of the Lorentz group'' $g_{\alpha \beta}$ with signature $(+++---)$.   A frame given by $e^{\mu}_a$ can be charactized by the Lorentz transformation $\theta$  which transforms the standard basis into the given one, which may now be written as $e^{\mu}_a(\theta)$.
The metric on the Lorentz group component was  derived from  the group operation on the frames, by measuring  the Minkowski squared norm induced by infinitesimal rotations of the Euler angles (summation over all frame vectors):
\beqarr g_{\alpha \beta}(\theta) &=& - a^2\, \eta_{\mu\rho}\eta_{\nu\sigma}\,\omega^{\mu\nu}_{\alpha}(\theta)\omega^{\rho\sigma}_{\beta}(\theta) \\
\mbox{with}\qquad \omega^{\mu\nu}_{\alpha}(\theta) &=& g^{\rho\nu} e^a_{\rho}(\theta)\, \frac{\partial}{\partial {\theta^{\alpha}}} e^{\mu}_a(\theta) \nonumber
\eeqarr 
The  factor  $a^2$ was  not  mentioned  at this place and   only  made explicit by the authors   in passing elsewhere \citep[p. 3313]{DeMartini/Santamato:2014a}.  Here $a$  expressed the {\em gyromagnetic radius} of a relativistic top $a=\sqrt{6} \frac{\hbar}{mc}$ and was 
  important, because due to it the geometry would ``know'' about the mass of the spinning particle. For the metric they found a  constant Riemannian scalar curvature $_gR=\frac{6}{a^2}= \frac{(m c)^2}{\hbar^2}$ induced from the Lorentz component \citep[p. 3313]{DeMartini/Santamato:2014a}.

This was a surprising {\em Riemannian geometrization of the configuration space} of the, up to here, non-quantum, relativistic top by a  non-definite metric $g_{ij}$ of signature $(3+3, 1+3)$ with constant scalar curvature.

\subsubsection{Santamato's random processes in the 1980s}
At first, the particle motion in the configuration space had to be analyzed. Santamato characterized it as a (weak) random process described by an ensemble of  trajectories $q^i(t,\omega)$ with $\omega$ ``the sample tag''\footnote{That is, $\omega\in \Omega$, the sample space of a probability triple $(\Omega,\mathfrak{F}, P)$, where $\mathfrak{F}\subset \mathfrak{P}(\Omega)$ are the random events and $P$ is a probability measure on $\Omega$.}  
and a  well defined and normalized probability density $\overline{\rho}(q,t)$ satisfying the 
 continuity equation \citep[p. 217]{Santamato:1984a}
\beq \partial_t \,  \overline{\rho} + \partial_i \,\left(   \overline{\rho} \, v^i\right) = 0  \, . \label{continuity equ. I}
\eeq 

He gave a peculiar derivation for the velocity field $v^i$ of his random process associated to a given  Lagrangian $L(q,\dot{q},t)$. After  shifting the   Lagrangian to $L^{\ast} = L + \frac{d}{dt}S$ for some sufficiently differentiable function $S(q,t),$\footnote{$L^{\ast}$ has  has the same Euler-Lagrange equations as  the original $L$.}
 he analyzed the {\em averaged action functional} 
\beq I(t_o,t_1)= E\left( \int_{t_o}^{t_1} L^{\ast}(q(t,\omega), \dot{q}(t,\omega),t) dt \right) \,  \label{averaged action functional}
\eeq 
with $E(\ldots)$ the expectation value. He looked for the minimum of $I$ under variation of $v^i=\dot{q}^i$, with respect to all random motions obeying a flow equation and satisfying given  initial data. As a  necessary condition for the existence of such a minimum it turned out that $S$ has to solve  the {\em Hamilton-Jacobi} equation 
\citep[app. A]{Santamato:1984a}
\beq \partial_t S + H(q, \nabla S,t) = 0 \, ,\label{Santamatos Hamilton-Jacobi equ}
\eeq  
with $H(q,p,t)$ the classical Hamiltionian corresponding to $L(q,\dot{q},t)$. Then the minimizing velocity field of (\ref{averaged action functional}) is the corresponding Hamilton-Jacobi flow. 

  Santamato could  hope that the wave equations of QM might be derivable from his random processes if a  classical Lagrangian (if there is any) was modified in a convincing way. ``Convincing'' would mean for him a change of the Lagrangian by geometrical terms, where the geometry is  influenced by the particle's (random) motion.
\begin{quote}
Geometry is not prescribed; rather it is determined by physical reality. In
turn, geometry acts as a ``guidance field'' for matter. \citep[p. 216]{Santamato:1984a}
\end{quote}
 He argued that such a ``feedback mechanism between geometry of space and particle motion'' was ``quite analogous'' to general relativity and might lead to ``a theory that is
physically indistinguishable from traditional quantum mechanics'' (ibid.).

At this point Santamato complemented the originally 
 Euclidean, or more generally Riemannian, basic structure of the configuration space by a {\em Weylian scale connection}. He called it a ``vector transplantation law'' and denoted  it by $\phi_k$,  corresponding to our $-\varphi_k$. In the case of a non-Euclidean Riemannian component of the metric $g_{ij}$ the continuity equation for the adapted probability density $\rho={|g|}^{-\frac{1}{2}}\,\overline{\rho}$ turns into the covariant equation:
\beq \partial_t \rho + \, _g\hspace{-0.2em}\nabla\hspace{-0.1em}_i (\rho v^i) = 0 \label{continuity equ II}
\eeq 

A classical Lagragian $L_c(q,\dot{q},t)$  on the original non-relativistic configuration space was then  modified on the Weylianized space according to \\[0.2em] {\em Santamato's 1-st postulate} \citep[equ. (8)]{Santamato:1984a}:
\beq
L(q,\dot{q},t) = L_c(q,\dot{q},t) + \gamma \frac{\hbar^2}{2 m} R(q,t)\, , \qquad \mbox{with} \quad \gamma=\frac{n-2}{4(n-1)} \, , \label{Santamatos 1st postulate}
\eeq 
 where $R(q,t)$ denotes the complete Weylian scalar curvature. 
With a  sign inverted convention for the scalar curvature (\ref{Weylian scalar curvature}),\footnote{Cf. fn \ref{fn curvature conventions}.} 
Santamato   wrote it as 
\beq R = \, _gR +(n-1)(n-2)\,\phi_i \phi^i - 2(n-1)\, _g\hspace{-0.2em}\nabla\hspace{-0.1em}_i \, \phi^i  \, .\label{Santamatos scalar curvature I}
\eeq
The  term in $R$  enters (\ref{Santamatos 1st postulate}) like an add on to the potential. The Hamilton-Jacobi equation of the random flow (\ref{Santamatos Hamilton-Jacobi equ}) thus becomes
\beq 
 \partial_t S + H_c(q, \nabla S,t)- \gamma \frac{\hbar^2}{m}R = 0 \, .\label{Santamatos Hamilton-Jacobi equ II}
\eeq

According to our author's program, $R$ should depend on the random process and was assumed to be time dependent. 
Therefore the $\phi_k$ cannot be arbitrarily given but ought to be determined by the probability density of the matter flow in  the configuration space.  Santamato  applied his averaged least action principle (\ref{averaged action functional}) another time and  evaluated it with (\ref{Santamatos 1st postulate})  for vanishing $L_c$, i.e., for $R(q,t)$ alone, with the encouraging result (ibid. equ. (19))
\beq
\phi_i = - (n-2)^{-1}\, \partial_i \ln \rho \, . \label{Santamatos scale connection}
\eeq 
Then the  scalar curvature   $ R = \, _gR +\,_{\varphi}R$ turned out to contain \citep[equ. (20)]{Santamato:1984a}
\beq \,  _{\varphi}R
 = \frac{1}{\gamma \sqrt{\rho}}\;\left(   _g\hspace{-0.2em}\nabla\hspace{-0.1em}_i \, \partial^i \sqrt{\rho} \right) \, .\label{Santamatos scalar curvature II}
\eeq
This form stood in  striking accord with  Bohm's quantum potential (\ref{Bohm quantum potential}).    Santamato jumped without hesitation from the recognition of the formal agreement to a  realistic conclusion: 
\begin{quote}
\ldots according to Eq. [(\ref{Santamatos scale connection})], the geometric properties of
space (\ldots) are indeed affected by the presence of the particle
itself. In turn, this alteration of geometry of space
acts on the particle through the quantum force
$f_i =\gamma \frac{\hbar^2}{m}\, \partial_i\, R$, which, according to Eq. [(\ref{Santamatos scalar curvature I})], depends
on the gauge vector and its first and second derivatives. \citep[219, equ. numbers adapted]{Santamato:1984a}
\end{quote}
This was a strong statement. It suggested a close {\em kinship of Santamato's  modification} of geometry to the one in the {\em general theory of relativity} (GR), although  his modification did   not refer to the spacetime manifold of GR, the ``extensive medium of the world'' as Weyl liked to formulate, but to the configuration space of a dynamical system. 

 In his next paper Santamato  derived the Klein-Gordon equation  (\ref{Klein-Gordon equ}) in the same way starting from a random process. He used  a  configuration  space arising from Minkowski space $\M$ by superimposing   a Weylian scale connection (\ref{Santamatos scalar curvature I}). Including  electromagnetic terms his Lagrangian for the    {\em relativistic ensemble} was 
\beq
L(x,\dot{x}) = \left(1+\gamma \frac{\hbar^2}{(mc)^2}R(x) \right)^{\frac{1}{2}}\, |x| + \frac{e}{mc^2} A_\mu\,  \dot{x}^\mu \, ,  \label{Lagrangian Klein Gordon ensemble}
\eeq
with $R=\, _{\varphi}R$ the Weylian scalar curvature ($\,_gR=0$). Taking into account equ. (\ref{Santamatos scalar curvature II}) it followed that  a  complex function  $\psi = \sqrt{\rho}\, e^{i S}$ constructed as usual (up to a factor $\hbar^{-1}$ in the exponent) from the flow quantities ``obeys the Klein-Gordon equation''\citep[p. 2479]{Santamato:1984b}. This was no small achievement; but Santamato did not continue his research along these lines for many years.

\subsubsection{A  look at the second phase in  cooperation with De Martini}
After a long interruption Santamato, now in joint work with  his colleague De Martini, gave a new derivation for the  Weyl geometric approach to the foundations of quantum mechanics. Moroever, in the new series of papers we find a  much clearer emphasis on the  underlying scale co-/invariant structure. 
The paper \citep{DeMartini/Santamato:2014a} started from a  field theoretic Lagrangian in a metric-affine approach (cf. pp. \pageref{metric-affine}ff.). It involved two scalar fields   $\rho, \sigma$ with weights $w(\rho)=-2,\, w(\sigma)=0$ under conformal rescaling  and a  scalar curvature term $R$ defined with regard to a metric  $g_{ij}$ and an independently defined affine (torsion free) connection $\Gamma_{ij}^k$. $\sigma$  now took  over the role of the former Hamilton-Jacobi principal function $S$ \citep[equ. (1)]{DeMartini/Santamato:2014a}.
\beq
\mathfrak{L}=\rho(\partial_{\mu}\sigma\partial^{\mu} \sigma + \gamma \hbar^2 R)\sqrt{|g|}\, . \label{Lagrangian ``geometrodynamics''}
\eeq
Variation with regard to the scalar fields leads to the dynamical equations:
\beqarr
\partial_{\mu}\sigma\,\partial^{\mu} \sigma +  \gamma \hbar^2 R &=& 0  \label{equ 1 ``geometrodynamics''} \\
\partial_{\mu}\left(\sqrt{|g|}\, \rho\,\partial^{\mu} \sigma\right)    &=& 0  \qquad \longleftrightarrow \qquad\,
_g\hspace{-0.2em}\nabla_{\mu} \left( \rho\, \partial^{\mu}  \sigma \right) = 0 \label{equ 2 ``geometrodynamics''}
\eeqarr
(\ref{equ 1 ``geometrodynamics''}) has the same form as the Hamilton-Jacobi equation of an uncharged Klein-Gordon field (\ref{Hamilton-Jacobi Klein-Gordon}), where the scalar curvature, up to sign, takes the place of the ``quantum potential''. (\ref{equ 2 ``geometrodynamics''}) may be read as  a continuity equation for a flow with density $\rho$ and velocity given by $\partial^{\mu} \sigma$ (if timelike).
By variation with regard to the affine connection,\footnote{Compare subsection \ref{section Palatini}.}
the authors concluded that the  affine connection has the Weyl geometric form  (\ref{Levi-Civita}) with the scale connection as in (\ref{Santamatos scale connection}), just like the one  Santamato had derived in the 1980s from his average action
 principle.\footnote{A sign error in  the  formula of the Weyl geometric affine connection \citep[equ. (4)]{DeMartini/Santamato:2014a} notwithstanding. }

The authors did not consider a variation of the metric because they had a  de Broglie--Madelung--Bohm context in mind in which the Riemannian metric of the configuration space was determined by the Lagrangian of a classical system.
They  immediately turned to it by a {\em mechanical interpretation} of their scalar field theory.  In the relativistic case the equations (\ref{equ 1 ``geometrodynamics''}), (\ref{equ 2 ``geometrodynamics''}) can be derived just as well as the Hamilton-Jacobi and continuity equations of a variational problem $\delta \int L\, d\tau = 0$ with 
\beq L_r= \sqrt{-\gamma \hbar^2 \, R(q) g_{\mu \nu}\, \dot{q}^{\mu}  \dot{q}^{\nu} }\label{Lagrangian relativistic system Santamato}
\eeq
This fits well to the program of geometrizing a configuration space with Riemannian metric related to a classical process, which is amended by  a Weylian scale connection standing  in  ``backreaction'' with a solution  pair  $(\sigma, \rho)$ of  (\ref{equ 1 ``geometrodynamics''}), (\ref{equ 2 ``geometrodynamics''}). 

With $\gamma= \frac{n-2}{4(n-1)} $ like above, and $n=4$, the Weylian component of the scalar curvature is in fact\footnote{\citep[p. 3310]{DeMartini/Santamato:2014a}}
\beqarr \, _{\varphi}R &=&  \frac{1}{4 \gamma}\,\left(  2 \rho^{-1}\,_g\hspace{-0.2em}\nabla_i \partial^i \rho -  \rho^{-2}\, \partial_i\rho\, \partial^i\rho    \right)
      = \gamma^{-1} \frac{\,_g\hspace{-0.2em}\nabla_i \partial^i \sqrt{\rho}}{ \sqrt{\rho}}  \\
			&=&  \gamma^{-1}\, \frac{2m}{\hbar^2 }  U =  \gamma^{-1}  Q \, ,
	\eeqarr 
where $U$ and $Q$ are the additional terms (``quantum potentials'') (\ref{Bohm quantum potential}), (\ref{quantum potential Klein-Gordon simple}) on the right hand side of the Hamilton-Jacobi equations of a Schr\"odinger, respectively a Klein-Gordon particle. It has to be understood that  (\ref{Lagrangian relativistic system Santamato}) holds only for  relativistic particles (Klein-Gordon and Dirac), while for the non-relativistic case of a Schr\"odinger particle  with $R=\,  _{\varphi}R $ the Lagrangian is (\ref{Santamatos 1st postulate}).

For investigating {\em relativistic spinning} particles De Martini and Santamato
considered  a point dynamics with internal degrees of freedom  in the configuration space $V=\M \times SO(3,1)$ described in subsection \ref{subsection configuration space}.  The Hamilton-Jacobi equation of a process  governed by the Lagrangian (\ref{Lagrangian relativistic system Santamato}) plus an electromagnetic term $L_{em}$
 is given by \citep[equ. (7)]{DeMartini/Santamato:2013}
\beq (\partial_{\mu}S - \frac{e}{c}A_{\mu})(\partial^{\mu} S- \frac{e}{c}A^{\mu}) + \hbar^2 \gamma R = 0 \, , \label{Hamilton Jacobi equ conformal top}
\eeq
 where $S$ satisfies  the divergence equation
\beq  D_{\mu} (\partial^{\mu} S - \frac{e}{c}A^{\mu}) =0 
\eeq 
with the scale covariant derivative, here $D_{\mu} = \nabla_{\mu}  -2 \phi_{\mu}$ in our weight convention with $w(g_{\mu \nu})=2$, and with $\nabla_{\mu}$ the Weyl geometric covariant derivative. For a current defined by
$ \mathfrak{j}^{\mu} = \chi^{-(n-2)}\sqrt{|g|}\,(\partial^{\mu}S - \frac{e}{c}A^{\mu})$  
this boils down to an ordinary continuity equation 
\beq
\partial_{\mu} \mathfrak{j}^{\mu}= 0 \, . \label{continuity equ top}
\eeq
  Obviously $ \mathfrak{j}^{\mu}$ is scale invariant. 

The  transition to a  complex wave function depending on all coordinates $q$ of the configuration space
\beq
\psi(q) = \sqrt{\rho}\,e^{\frac{i}{\hbar}S} \,  \label{wave function top}
\eeq 
 transforms the equs. (\ref{Hamilton Jacobi equ conformal top}), (\ref{continuity equ top}) into the {\em linear differential equation} of second order
\beq
(\hat{p}^{\mu} - - \frac{e}{c}A^{\mu})(\hat{p}_{\mu} - - \frac{e}{c}A_{\mu})\,\psi + \hbar^2 \gamma \, _gR \; \psi = 0 \, ,\label{pre-Dirac equ top}
\eeq
where  $\hat{p}$ denotes the differential operator with $\hat{p}_{\mu} = - i \hbar \partial_{\mu}$. It has the form of the Klein-Gordon equation (\ref{Klein-Gordon equ})
with a mass factor which contains only  the {\em Riemannian part} of the scalar curvature, $ \hbar^2 \gamma\,\, _gR = \hbar^2 \gamma\,\frac{6}{a^2}= m^2 c^2 $. The  {\em Weylian component} $_\varphi R$ is controlled via (\ref{Santamatos scale connection}) by the density of the quantum flow. 
The authors commented
\begin{quote}
This is a striking result as it demonstrates that the Hamilton-
Jacobi equation, applied to a general dynamical problem can be transformed into a
linear eigenvalue equation, the foremost ingredient of the formal structure of quantum
mechanics and of the Hilbert space theory. \citep[p. 636]{DeMartini/Santamato:2013}
\end{quote}
The first step towards a reconstruction of the Hilbert space quantization of the relativistic top was achieved.

In the next step the authors analyzed the decomposition of a solution $\psi$ of (\ref{pre-Dirac equ top})  into components $\psi_{u,v}$ lying in  finite dimensional representations  of $SO(3,1)$ of type $D^{(u,v)}$ with $2u, 2v \in \{0, 1, 2, \dots  \} $. Then $\psi_{u,v}(q)$  can be factorized into functions of the spatial variable $x$  with values in the representation space of $D^{(u,v)}$, and  $\theta$-dependent representation matrices operating on the latter; in  spinor notation similar to van der Waerden's symbolism:\footnote{For van der Waerdens spinor symbolism see \citep{Schneider:Diss}.}
\beq \psi_{u,v}(q) = D^{(u,v)}(\Lambda(\theta))_{\sigma'}^{\sigma}\psi_{\sigma}^{\sigma'}(x) +  D^{(v,u)}(\Lambda(\theta))_{\dot{\sigma'}}^{\dot{\sigma}}\psi_{\dot{\sigma}}^{\dot{\sigma'}}(x) \,
\eeq
For the particular choice $u=v=\frac{1}{2}$ this leads to a pair of 2-component spinor field on Minkowski space, equivalent to  a  4-component Dirac field 
$\Psi(x)= \left(  \begin{array}{c }
                     \psi_{\sigma}^{\sigma'}(x)        \\
												\psi_{\dot{\sigma}}^{\dot{\sigma'}}(x)  
\end{array}
\right) $ in the Weyl representation. Then the equation (\ref{pre-Dirac equ top}) acquires a form which, after neglecting an extremely small term in the electromagnetic field strength,\footnote{This term, $ \frac{e^2 a^2}{c^2}(H^2-E^2)$, is  comparable with the linear term in the field strenghts only under the condition of very large field strenghts, $E\sim 10^{18}\, V m^{-1}, \, H \sim 10^9\, T$. ``To have an idea how large is this field, an electron at rest is accelerated by such field up to $10^9\, GeV$ in a linear accelerator 1 m long'' \citep[p. 3315]{DeMartini/Santamato:2014a}.}
  was identified by the authors as  the {\em squared Dirac equation}
\beqarr
\mathfrak{D}_+ \mathfrak{D}_- \Psi &=& \mathfrak{D}_- \mathfrak{D}_+ \Psi = 0 \, ,\label{squared Dirac equ}  \\
\mbox{where} \qquad      \mathfrak{D}_\pm &=&   \gamma^{\mu}(p_{\mu} -\frac{e}{c} A_{\mu})  \pm m         \nonumber
\eeqarr
with the Dirac matrices $ \gamma^{\mu}$ ($\mu = 0, \ldots, 3$) \citep[p. 639]{DeMartini/Santamato:2013}.\footnote{It remains unclear to me (E.S.) how the  representation matrices of the ``Euler angles'' of configuration space are suppressed, while the change of coordinate frames in Minkowski space gets represented on the spinor fields.}

Solutions of (\ref{squared Dirac equ}) can be decomposed into a superposition of the linear Dirac equation with positive and with negative mass. The authors proposed that the negative mass contributions have to be ``disregarded as unphysical'' \citep[p. 641]{DeMartini/Santamato:2013}. Even without trying to assess this proposal, it is clear that  by this  model of relativistic spinning particles Santamato and De Martini had achieved a surprising step forward for the geometrization program of the dBMB  approach started in the 
1980s.

They did not stop here, but went on by investigating  the nonlocality of EPR systems in their approach. Their considerations led to a justification of  the spin-statistic relation which usually is derived by quantum field methods \citep{DeMartini/Santamato:2014a,DeMartini/Santamato:2015a}. In order  not to blow our survey  these derivations, although central for the content  of their papers, have to be shunted here.

\subsection{\small An attempt at bridge building to gravity\label{subsection Shojai ea}}
We still have to review  attempts  at  connecting the dBMB approach to gravity with a specific reference to Weyl geometry. Different authors tried to do  so. The main thrust in this direction  was developed independently of the two Italian authors by {\em Fatimah Shojai, Ali Shojai} and {\em Mehdi Golshani} working at Tehran.  Another,  to my taste slightly more bizarre, step in this direction was made by {\em Giorgio Papini} and {\em Robert Wood} at the occasion of a symposium honouring J.-P. Vigier \citep{Wood/Papini:1997}. Some years earlier they had tried to fix a defect of  Dirac's 1972 proposal to revive Weyl's original interpretation the scale connection as the electromagnetic potential, resulting from the non-integrability of the scale connection.\footnote{Among others, this had led to the measurement problem by atomic clocks.}
  Papini and Wood proposed to solve this problem by considering ``bubbles'' in the environment of atoms, in which the scale symmetry is broken, while it holds in the large, outside the ``bubble'' \citep{Wood/Papini:1992}. For the Vigier symposium they recycled their idea by establishing a connection to a  dBMB approach governing the dynamics in the bubble,  similar to de Broglie's proposal.

At the end of the 1990s F.  and A. Shojai, sometimes coauthored by Golshani, started with investigations of their own, in which they hoped to be able to use a Bohmian approach for a peculiar way of quantizing a part of the gravitational structure \citep{Shojai/Shojai/Golshani:1998a,Shojai/Shojai/Golshani:1998b,Shojai/Shojai/Golshani:1998c,Shojai/Golshani:1998,Shojai/Shojai:2000}. To do so they used methods from scalar-tensor theories of gravity. A specific emphasis of conformal ideas brought their approach close to Weyl geometry. During a sojourn at the {\em Max Planck Institute for Gravitational Physics} at Potsdam they laid this connection open  and proclaimed it as the correct framework of their approach \citep{Shojai/Shojai:2003}. We want to see what that meant.

In the 1980s {\em Jayant  Narlikar}  and {\em Thanu Padmanabhan} had started to study a simplified version of  quantum gravity which was invariant under conformal changes of the metric, $g_{\mu \nu} \mapsto   \Omega^2 g_{\mu \nu}$. They  proposed to {\em quantize  only the factor} $\Omega$, viz. the scale degree of freedom of the metric. This had the great advantage of keeping the conformal structure unaffected by the quantization and circumvented the infamous obstacle of a fuzzy causal structure, which other approaches towards quantum gravity encountered. On this basis Narlikar and Padmanabhan  calculated semiclassical approximations for cosmological solutions of the Einstein equation \citep{Narlikar/Padmanabhan:1983,Padmanabhan:1989}.\footnote{In the physics literature, so also in the paper by Shojai and Golshani,   $\Omega$ is often talked about as the  ``conformal'' degree of freedom of the metric, or even the ``conformal structure''.The latter  is clearly mistaken, the first one at least  misleading. Therefore I avoid this terminology in favour of  {\em scaling degree of freedom}. }

In one of their early joint papers F. Shojai and M. Golshani took  this  idea up. In contrast to Narlikar and Padmanabhan  they attempted a  Bohmian path towards quantizing the scale factor \citep{Shojai/Golshani:1998}. This was quite daring  because Bohmian quantum mechanics had been developed for systems of finite degrees of freedom only.  Shojai and Golshani, however, invoked the idea of de Broglie to re-interpret  the ``quantum mass'' $\mathfrak{M}_o^2 = m^2 + \frac{\hbar^2}{c^2}Q$    of a Klein-Gordon system (\ref{Hamilton-Jacobi Klein-Gordon})  as a conformal modification of the Minkowski metric using  $\Omega^2 = \frac{\mathfrak{M}_o^2}{m^2}=1 + \frac{\hbar^2}{m^2 c^2} \frac{\nabla_{\mu}\partial^{\mu}\sqrt{\rho}}{\sqrt{a}} $. They considered this rescaling factor  as a representative for the quantum degrees of freedom of a globally defined Klein-Gordon field.   

Another problem was that $\mathfrak{M}_o^2$  could become negative. The Shojais and Golshani solved  it    by  passing over to  the exponential \citep[equ. (12)]{Shojai/Golshani:1998}\,
\beq
\mathfrak{M}^2 = m^2 e^{ \frac{\hbar^2}{m^2 c^2} \frac{\nabla_{\mu}\partial^{\mu}\sqrt{\rho}}{\sqrt{a}}} \, . \label{Shojais' quantum mass}
\eeq
The linear approximation coincides with  $\mathfrak{M}_o^2$, and the  conformal factor became
\beq
\Omega^2 = \frac{\mathfrak{M}^2}{m^2}=e^{ \frac{\hbar^2}{m^2 c^2} \frac{\nabla_{\mu}\partial^{\mu}\sqrt{\rho}}{\sqrt{\rho}}} \, . \label{Shojais conformal factor}
\eeq 

They started from a Lagrangian with Einstein-Hilbert term and a matter Lagrangian in terms of a Hamilton-Jacob function $S$ and flow density $\rho$,
\beq \mathfrak{L}_m=\frac{\hbar^2}{m}\left(\frac{\rho}{\hbar^2}\partial_{\mu}S \, \partial^{\mu}S - \frac{m^2}{\hbar^2} \rho \right)\sqrt{|g|} \, 
\eeq 
characteristic for a classical Jacobi-Hamilton system.  Santamato's viewpoint (which was apparently unknown to the Tehran authors)  had been  that the 
 introduction of the quantum potential turned the corresponding dynamical system into the Hamilton-Jacobi form of a Klein-Gordon system (\ref{Hamilton-Jacobi Klein-Gordon}), (\ref{continuity equ Klein-Gordon}). The Shojais proceeded differently.  
Sustained by de Broglie's argumentation they argued:
\begin{quote}
\ldots the de Broglie remark leads to the conclusion that the introduction of the quantum potential which contains the quantal behaviors of the particles is {\em equivalent} to the introduction of a conformal factor $\Omega^2 = \frac{\mathfrak{M}^2}{m^2}$ in the metric \citep[p. 683, emphasis E.S.]{Shojai/Golshani:1998}.
\end{quote}

This was a  puzzling statement. De Broglie had considered a geometrization for a single particle  {\em in the absence of electromagnetic and gravitational fields} (subsection \ref{subsection de Broglie}). It remained unclear whether the argument could be transferred to the case of gravitational fields and in which sense such an  ``equivalence'' was to be understood.  

In  Santamato's  geometrization of a  dBMB Hamilton-Jacobi system     the  ``prepotential'' of a Weylian scale connection on the configuration space  (\ref{Santamatos scale connection}) was  $\ln \rho$, up to a constant factor. It  leads to the Weylian curvature expression  (\ref{Santamatos scalar curvature II}) equivalent to a Bohmian ``quantum potential'' in the dynamical equation (\ref{Santamatos Hamilton-Jacobi equ II}).  De Broglie  and with him the Shojais used a different geometrizaton idea. Their ``prepotential'' of the Weylian scale connection was the scale factor  $\Omega$ between a classical metric and the  metric describing a  quantum system.  In the work of our authors it was  the exponential expression  (\ref{Shojais conformal factor}).  Following de Broglie, one  had to consider geodesic flows with the implicit  constraint of orthogonal initial conditions to a  level surface of a related Hamilton-Jacobi principal function $S$.
The  modification (\ref{Shojais' quantum mass}) of the usual ``quantum mass'' formula implies  that we cannot expect  equivalence in the  literal sense. 
Even if one wants to read the argument as a motivation for a new type of dBMB-like quantization procedure, following the de Broglie paradigm, a  justification for the attempted generalization from de Broglie's  case (no gravitation) to the general case considered had to be given.

But our  authors did not hesitate to take this step as  a starting point for investigating  cosmological models in which  matter fields were given  in different versions of scalar tensor theories.\footnote{\citep{Shojai/Golshani:1998,Shojai/Shojai/Golshani:1998a,Shojai/Shojai/Golshani:1998b,Shojai/Shojai:2000,Shojai/Shojai:2001}.} 
In the result  a Klein-Gordon field appeared on  large scales, rather than as a descriptor of the motion of a single quantum particle. At some places it played 
 the role of a matter field \citep[p. 683]{Shojai/Golshani:1998}, \citep[p. 1762]{Shojai:2000}, at others that of a ``quantum gravity'' modification of the metric field \citep[p. 2728]{Shojai/Shojai/Golshani:1998a}. 
The Shojais  were convinced 
\begin{quote}
\ldots that the theory works for a particle as well as for a real ensemble of the particle under consideration and that it includes pure quantum gravity effects \citep[1763]{Shojai/Shojai:2000}.
\end{quote}
But it remained unclear what  ``quantum gravity'' would mean here. 

One of the papers dealt explicitly with  conformal transformations in scalar-tensor theories \citep{Shojai/Shojai/Golshani:1998a}.\footnote{The authors made a difference  between ``scale transformations'' and ``conformal transformations''. In their terminology the first operated only on the metric, while the latter rescaled all physical fields according to their weights.}
The three authors distinguished  between a  ``{\em background metric}'', in which they considered the quantum effects being encoded by the varying  ``quantum mass'' $\mathfrak{M}$, while in a ``{\em physical metric}''  $\mathfrak{M}$ was rescaled to a constant value $\overline{m}$. Then  ``some part of the curvature of space-time represent the quantum effects'' \citep[p. 2726]{Shojai/Shojai/Golshani:1998a}. Independent of the physical interpretation and reasonability of this and some other observations one might  wonder whether a reformulation in Weyl geometry could at least help to clarify the mathematical side of such statements.

This is what  A. and F. Shojai attempted 
in \citep{Shojai/Shojai:2003} and a following preprint \citep{Shojai/Shojai:2004}. In the meantime they had adopted Dirac's theory of 1972 (see section \ref{subsection Dirac}), but  did not follow Dirac's {\em em} dogma. They rather considered the scale connection as ``a part of the geometry of the space-time'', implicitly constrained in their context by the integrability condition.\footnote{The Dirac Lagrangian was stripped of the Yang-Mills term of the scale connection \citep[equ. 6]{Shojai/Shojai:2004}.}
But without much hesitation they declared  that   Dirac's scalar field $\beta$  in (\ref{scale invariant Dirac action})   ``represents the quantum mass field'' in the sense of their embryonic theory outlined above \citep[p. 7 preprint]{Shojai/Shojai:2003}. They did not discuss  how   the different Lagrangians for the Dirac field and their Klein-Gordon field could be related to each other. Only a rather opaque  perturbative argument was given as to why a solution of the $\beta$-scalar field equation may  be identified with an expression of the ``quantum mass'' type, $\beta \mapsto \mathfrak{M}$ \citep[p. 13f.]{Shojai/Shojai:2003}.\footnote{The claim of the possibility to identify a Dirac-type scalar field with a ``quantum mass''  field   remains, in my view, an unfounded  speculation; E.S.} 
On the other hand, this identification  allowed to clarify their discussion of different  frames  a little. They now considered ``different  conformal frames'' as ``identical pictures of the gravitational and quantum phenomena'' (ibid., p. 9).\footnote{The authors even conventionalized this idea as a ``conformal equivalence principle''\citep[p. 10]{Shojai/Shojai:2003},  \citep[p. 63]{Shojai/Shojai:2004}.}

In the light of such  open spots  A. and F. Shojai's  conclusion that Weyl geometry ``provides a unified geometrical framework for understanding the gravitational and quantum forces'' \citep[p. 10 preprint]{Shojai/Shojai:2003} was at least premature and reads like too grand a speculation. Not all readers had this impression. Their program  found at least one active successor, {\em R. Carroll} 
 \citep{Carroll:quantum-potential}.\footnote{{\em Not S.} Carroll, as it sometimes  erroneously appears in the bibligraphy of later papers.}
But the critical points of justification for the ``Tehran'' program seem not to be clarified in this work either.


\section{\small Scale covariance in the standard model of elementary particle physics\label{section SM}}
About  the middle of the 1970s the standard model of elementary particle physics (SM) started to become widely accepted as the key to the basic structures of matter  \citep[chap. 22]{Kragh:Quantum}, \citep{Pickering:Quarks}.
 Besides the point dependent (localized) {\em internal} symmetries of the electroweak forces, $SU(2)\times U(1)$ and the chromodynamic symmetry of the strong forces $SU(3)$ the new paradigm of gauge field theories worked with non-localized (``global'') {\em external} symmetries of special relativity, the Lorentz group. Characteristic for the paradigm was a global but only  nearly  respected scale invariance of the field Lagrangians,   broken only by the mass term of the Higgs field.  The Higgs field $\Phi$, a scalar field with values in an isospin $\frac{1}{2}$ representation of the electroweak group, was the clue for making electroweak symmetry of elementary particles consistent  with mass terms.  The latter was understood as a ``spontaneous breaking'' of the electroweak symmetry  and became to be known as the ``Higgs mechanism'' \citep{Borrelli:Higgs}.  In this section we look at some attempts for bridging the gap between the Higgs field and the  scalar field of gravity.

\subsection{\small Englert,  Smolin and  Cheng,  1970/80s\label{subsection SM 1970s}}
\subsubsection{A conformal approach}
One of the originators of this  theory (Higgs mechanism), {\em Fran\c{cois} Englert},\footnote{See, e.g.  \cite{Karaca:Higgs}.} tried to play a similar game of ``spontaneous symmetry breaking'' in gravity, here with a real valued scalar field with scale symmetry in the sense of conformal rescaling. 

 In a common paper written with  {\em Edgar Gunzig, C Truffin} and {\em P. Windey}, the authors established an explicit link to JBD gravity \citep{Englert/Gunzig:1975}. 
But in contrast to \citep{Deser:1970}, Englert and  coworkers   considered conformal gravity as part of the quantum field program. They assumed a ``dimensionless'', i.e. scale invariant,  Lagrangian for gravitation with a square curvature term of an affine connection  $\Gamma$ {\em not} bound to the metric,  $\mathcal{L}_{\text{grav}} = R^2 \, \sqrt{|det\, g|} $,  in addition to a Lagrangian matter term  \citep{Englert/Gunzig:1975}. In consequence, the authors  varied with respect to the metric $g$ and the connection $\Gamma$  independently. 

``To make contact with General Relativity'' (p. 74) the authors  assumed the scalar curvature as  expressed by a scalar function,  $R \sim \omega^2$ (they used the symbol  $\varphi$ instead of  $\omega$). The Euler-Lagrange equation of the affine connection  resulted  in a relation like  (\ref{Levi-Civita}) for the  Weyl geometric case, with an integrable integrable scale connection $\varphi = d \log \omega$ \citep[equs.(7), (8)]{Englert/Gunzig:1975}.
By such a specialization, their  approach looked as though it was  touching upon a Weyl structure. But this was not the point of view of the authors; they rather 
 proceeded as ``conformal'' as possible on their search for connecting paths between quantum field theory of scalar fields and  general relativity.

 After some tentative quantum considerations the authors came back to a  ``classical phenomenological description''  of their theory \citep[76]{Englert/Gunzig:1975}. For  this description they  introduced a scalar  field $\phi(x) =\lambda^{-1} e^{\lambda \sigma(x)}$ coupled to gravity like in our equ. (\ref{Lagrangian JBD Fujii/Maeda}),  with the necessary specification $\xi = \frac{1}{6}$ in order to secure  conformal symmetry \citep[equ. (16)]{Englert/Gunzig:1975}.  They considered $\sigma$ to be a ``dilation field'' (sic!) which represented     a ``Nambu-Goldstone boson'' coupling  to the mass terms. 

 After some turns and twists they summed up  that their original action principle
\begin{quote}
\ldots  matches all the results of General Relativity at a classical level, provided
mass originates in dynamical breakdown of symmetry. Thus, the fundamental finite component fields must
be massless and of the kind currently used in gauge field theories, but without scalar mesons \citep[76]{Englert/Gunzig:1975}.
\end{quote}

In one of the following papers Englert, now  with Truffin as only  coauthor,  studied the perturbative behaviour of his version of  conformal gravity  ($\xi = \frac{n-2}{4(n-1)}$) coupled to massless fermions and photons in $n \geq 4$ dimensions.\footnote{The motivation or considering $n \geq 4$ was the method of dimensional regularization for the quantization of the theory.} He came to the conclusion that anomalies arising in the calculations for non-conformal actions disappeared at the tree and 1-loop levels in their approach. The two authors took this as an indicator that  gravitation might perhaps arise in a ``natural way from spontaneous breakdown of conformal invariance''  \citep[426]{Englert_ea:1976}. 

\subsubsection{Smolin introduces  Weyl geometry}
Englert's e. a. paper was one of the early steps into the direction (i) of our introduction. 
Other authors followed and extended this view, some of them explicitly in a Weyl geometric setting, others continued to use the language of conformal geometry. The first strategy was  chosen by  {\em Lee Smolin} in his paper \citep{Smolin:1979}. In section 2 of the paper he gave an explicit and clear introduction to Weyl geometry.\footnote{In his bibliography he went back directly to \citep{Weyl:STM} and \citep{Weyl:GuE}; he  did not   quote any of the later literature on Weyl geometry. }
The  ``conformally metric gravitation'', as he called it, was built  upon a matter-free Lagrangian built from Weyl geometric curvature terms $R, \, Ric = (R_{\mu \nu })$,  $f=(f_{\mu \upsilon })$ for scale curvature and used a gravitational Lagrangian of order two.  In  a slight adaptation of notation using the  scale covariant Weylian derivatives $D$ it was \cite[equ. (13)]{Smolin:1979}: 
\beqarr  |det\, g|^{-\frac{1}{2}}  \mathcal{L}_{\text{grav}} &=&  -\frac{1}{2}c \,  \phi^2 R + \; [  - {e_1} R^{\mu \nu }R_{\mu \nu } -  {e_2}  R^2 ] \label{Smolin's Lagrangian}\\
& & + \quad \frac{1}{2}D^{\mu}\phi \,D_{\mu}\phi -\frac{1}{4g^2}f_{\mu \nu } f^{\mu \nu } - \lambda \phi^4 \nonumber
 \eeqarr 
 where  $c, e_1, e_2, g, \lambda $ are coupling coefficients.\footnote{Signs  have to be taken with caution. They may depend on conventions for defining the Riemann curvature, the Ricci contraction, and the signature. Smolin, e.g., used a different sign convention for $Riem$  to the one used in this survey. Signs given here are adapted to ${signature}\, g = (3,1)$. The Riemann tensor and Ricci contraction are those usually adopted in the   mathematical literature, see fn. \ref{fn curvature conventions}. } 

For  coefficients 
of the quadratic curvature terms (in square brackets) with  $e_2 = - \frac{1}{3}e_1 $, the latter  is variationally equivalent (equal up to divergence) to the  squared conformal curvature $C^2 = C_{\mu \nu \kappa \lambda } C^{\mu \nu \kappa \lambda }.$\footnote{This seems to have been widely  known. For an explicit statement see, e.g., \citep{Hehl_ea:1996quad_curv}.}

Smolin introduced the scalar field $\phi$ not only by formal reasons (``to write a conformally invariant Lagrangian with the required properties''), but with a physical  interpretations similar to those given by Englert e.a.,\footnote{\citep{Englert/Gunzig:1975} was not quoted by Smolin.}
namely ``as an order parameter to indicate the spontaneous breaking of the conformal invariance'' \citep[260]{Smolin:1979}. His Lagrangian used  a modified adaptation from JBD theory, ``with some additional couplings'' between the scale connection  $\varphi$ and the scalar field $\phi$.  Smolin  emphasized that ``these additional couplings go against the spirit of
 Brans-Dicke theory''  because from the Riemannian point of view they introduced a non-vanishing divergence of the non-gravitational fields.

For low energy considerations  Smolin dropped the square curvature term (square brackets in  (\ref{Smolin's Lagrangian})),  added an ``effective'' potential term of the scalar field $V_{\text{eff}}(\phi )$ and derived the equations of motion by varying with respect to  $g, \phi, \varphi$. Results were  Einstein equation,  scalar field equation, and Yang-Mills equation for the scale connection.

 Smolin's reduced Lagrangian contained  terms in the scale connection:\footnote{In scalar field gauge with $\phi \doteq \phi_o= F$, his  reduced  Lagrangian (square gravitational terms dropped) was  \citep[equ. (3.17)]{Smolin:1979} \[  |det\, g|^{-\frac{1}{2}}  \mathcal{L}_{grav} \doteq - \frac{1}{2}c \, F^2\,  _g\hspace{-1pt}R -\frac{1}{4g^2}f_{\mu \nu } f^{\mu \nu } + \frac{1}{8}(1+6c)F^2 \varphi_{\mu} \varphi^{\mu } - V_{\text{eff}}(F) \; . \] 
}
\beq  -\frac{1}{4g^2}f_{\mu \nu } f^{\mu \nu } + \frac{1}{8}(1+6c)F^2 \varphi_{\mu} \varphi^{\mu }
\eeq 
That looked like a mass term for the scale connection  $\varphi$, the  potential of the scale curvature field $f_{\mu \nu }$ called ``Weyl field'' by Smolin. By comparison with the Lagrangian of the Proca equation in electromagnetic theory,  Smolin   concluded that the ``Weyl field'' has  mass close to the Planck scale, given by
\beq  M^2_{\varphi}  =  \frac{1}{4}(1+6c)F^2 \, . \eeq 
He commented that in his Weyl geometric gravitation theory  ``general relativity couples to a massive vector field'' $\varphi$. The scalar field $\phi$, on the other hand, ``may be absorbed into the scalar parts'' of $g_{\mu \nu }$ and $\varphi_{\mu}$,\footnote{It is possible to choose the scale gauge such that $\phi$ becomes constant (scalar field gauge, see section \ref{subsection Weyl geometry}.}
 by a change of variables and ``remains massless'' \citep[263]{Smolin:1979}. In this way, 
Smolin brought Weyl geometic gravity closer to the field theoretic frame of particle physics. He did not discuss mass and interaction fields of the SM. Morover, the huge mass of the ``Weyl field'' must have appeared  irritating.

\subsubsection{Interlude}
At the time Smolin's paper appeared,  the program of  so-called {\em induced gravity}, entered an  active phase. Its central goal was to derive the action of conventional or modified Einstein gravity from an extended scheme of standard model type quantization. Among the authors involved in this program {\em Stephen Adler} and {\em Anthony Zee} stand  out.   We cannot go  into this story here.\footnote{ For a survey of the status of investigations in 1981 see  \citep{Adler:Report_1982}; but note in particular 
\citep{Zee:1982a,Zee:1983}. The topic of ``origin of spontaneous symmetry breaking'' by radiative correction was much older \citep{Borrelli:Higgs,Karaca:Higgs}. A famous paper was \citep{Coleman/Weinberg:1973}. 
In fact, Zee's first publication  on the subject preceded Smolin's. \citep{Zee:1979} was submitted in December 1978 and  published in February 1979; \citep{Smolin:1979} was submitted in June 1979.}

 Smolin's view that  the  structure of Weyl geometry might be suited  to bring classical gravity into a coherent frame with standard model physics did not find much immediate response. But it was ``rediscovered''  at least twice (plus  an independently developed conformal version). In 1987/88   Hung Cheng at the MIT,   and a  decade later  Wolfgang Drechsler and  Hanno Tann  at Munich,   arrived  at similar insights and established an explicit extension of Weyl geometric gravity to standard model (SM) fields  \citep{Cheng:1988,Drechsler/Tann,Drechsler:Higgs}. Simultaneous to  Cheng, the core of the  idea was  once more discovered  by  Mosh\'e Flato (Dijon)  and Ryszard R\c{a}cka  (during that time at Trieste), although  they formulated it in a strictly conformal framework without  Weyl structure  \citep{Flato/Raczka}. 
Neither  Cheng, nor Flato/R\c{a}cka or Drechsler/Tann seem to have known Smolin's proposal (at least Smolin was not cited by them), nor did they refer to the papers of each other.\footnote{Flato/R\c{a}cka's paper appeared as a preprint of the {\em Scuola Internazionale Superiore di Studi Avanzati}, Trieste, in 1987; the  paper itself  was submitted in December 1987 to {\em Physics Letters B} and published in July 1988.   Cheng's paper was submitted in February 1988,  published in November. Only a  decade later, in March 2009, Drechsler and Tann got acquainted with the other two papers. This indicates  that the Weyl geometric approach in field theory had not yet acquired the coherence of a research program with a stable  communication network.  }
All three approaches had their own achievements. Here we can give only give a short presentation of the main points of the work directly related to Weyl geometry.

\subsubsection{Hung Cheng and his  ``vector meson''}
{\em Hung Cheng} started out from a Weyl geometric background, apparently inherited from the papers of the Japanese  group of authors around Utiyama. The latter had taken up Weyl geometry in the early 1970s in a way not too different from Smolin's later approach (see section \ref{subsection Utiyama}).
Cheng   extended Utiyama's theory  explicitly to the electroweak sector of the SM. He replaced the complex scalar field $\phi$  by the {\em Higgs field} $\Phi $, again  of weight $-1$ but now  with values in an  isospin $\frac{1}{2}$  representation, and coupled it to the  Weyl  geometric scalar curvature  $R$ and postulated:\footnote{In the sequel the isospin extended scalar field will  be denoted by $\Phi $.}
\beqarr  \mathcal{L}_{R} &=&    \frac{1}{2} \beta \,  \Phi^{\ast} \Phi R  \, |det\, g|^{\frac{1}{2}}  \label{Hilbert term Weyl geometric} \\
 \mathcal{L}_{\Phi} &=&  \quad \frac{1}{2}\tilde{D}^{\mu} \Phi^{\ast}  \tilde{D}_{\mu}\Phi   \, |det\, g|^{\frac{1}{2}}  \; \label{phi term Weyl geometric}
\eeqarr
The  scale covariant derivatives were extended to a  localized  electroweak ({\em ew}) group $SU(2)\times U(1)$. With the usual denotation of  the standard model,   $W^j_{\mu} $ for the field components of the $su(2)$ part (with respect to the Pauli matrices $\sigma _j\; (j=0,1,2)$) and  $B_{\mu}$ for $u(1)_Y \cong \R$ and coupling coefficients $g, g'$ the derivative  read\footnote{Cheng added another coupling coefficient for the scale connection, which is here suppressed. }
\beq  \tilde{D}_{\mu}\Phi    = (\partial _{\mu} - \varphi_{\mu} + \frac{1}{2} i g W^j_{\mu} \sigma_j + \frac{1}{2} g'B_{\mu})\Phi \, . \label{ew covariant derivative}
\eeq 
 The sign of the kinetic term of the Higgs field (\ref{phi term Weyl geometric}) shows that Cheng  supposed $sig\; g =(+ - - - )$, which agrees with his high energy context, while  the sign of (\ref{Hilbert term Weyl geometric})  indicates   that he used the sign inverted  convention for curvature.\footnote{See fn \ref{fn curvature conventions}.}
He
 added Yang-Mills interaction Lagrangians for the {\em ew} interaction fields $F$ and $G$ of the potentials $W$ (values in $su_2$), respectively $B$ (values in  $u(1)_Y$),   and added a scalar curvature term in  $f=(f_{\mu \nu}) = d\varphi$
\beq   \mathcal{L}_{\text{YM}}= -\frac{1}{4}\left( f_{\mu \nu}  f^{\mu \nu}  + F_{\mu \nu}  F^{\mu \nu}  + G_{\mu \nu}  G^{\mu \nu}  \right) \, |det\, g|^{\frac{1}{2}} \, .  \label{Cheng's ew interaction Lagrangian} \eeq
Finally he introduced spin $\frac{1}{2}$ fermion fields $\psi$ with the  weight convention $w(\psi)=-\frac{3}{2}$, and a Lagrangian  $\mathcal{L}_{\psi}$ similar to the one formulated later by Drechsler,  discussed below   (\ref{Lagrangian Dirac field}).\footnote{The second term in (\ref{Lagrangian Dirac field}) is missing in Cheng's publication. That is probably not intended, but a misprint. Moreover he did not discuss scale weights for Dirac matrices in the tetrad approach.} 

Thus Cheng's general relativistic scalar field $\Phi$ resembled very much the Higgs field of the SM, which at that time was still a highly hypothetical object.  He 
 called  the scale connection, respectively its curvature, {\em Weyl's meson} field. Referring to Hayashi's e.a. observation  that  the scale connection does not influence the equation of motion of the spinor fields, he  concluded:\footnote{Remember that the $\varphi$ terms of scale covariant derivatives in the Lagrangian of spinor fields cancel.}
\begin{quote}
\ldots  Weyl's vector meson does not interact with leptons or quarks. Neither does it interact with other vector mesons. The only interaction the Weyl's meson has is that with the graviton. \citep[2183]{Cheng:1988}
\end{quote}

 Because of the tremendously high mass of ``Weyl's vector meson'' Cheng conjectured that even such a minute coupling might  be of some cosmological import. More precisely, he   wondered, ``whether Weyl's meson may account for at least part of the dark matter of the universe'' (ibid.). Similar conjectures were stated once and again over the next decades,  if theoretical entities were encountered which might represent massive particles without  experimental evidence. Weyl geometric field theory was not spared this fate.  

\subsubsection{Can gravity do what the Higgs does?  }
In the same year in which  Cheng's paper appeared, {\em Mosh\'e Flato} and {\em Ryszard  R\c{a}czka} sket\-ched  an approach  in which they  put gravity into a quantum physical perspective.\footnote{More than a decade earlier Flato had sketched a covariant (``curved space'') generalization of the Wightman axioms  \citep{Flato/Simon}, different from the one discussed by R. Wald in this volume. }
In our context, this paper matters because it introduced a  scale covariant Brans-Dicke like field in an isospin representation similar to Hung Cheng's, but in a strictly conformal framework  \citep{Flato/Raczka}. 

Six years later, R. R\c{a}czka took up the thread again, now in  cooperation  with {\em Marek 
Paw{\l}owski}.   In the meantime Paw{\l}owski 
 had joined the research program  by a  paper  in which he addressed the question whether  gravity ``can do what the Higgs does'' \citep{Pawlowski:1990}. In  a couple of  preprints\footnote{
 \citep{Pawlowski/Raczka:1994a,Pawlowski/Raczka:1995_0,Pawlowski/Raczka:1995} }
and two   refereed papers 
\citep{Pawlowski/Raczka:1994FoP,Pawlowski/Raczka:1995} the two physicists proposed   a ``Higgs free model for fundamental interactions'', as they  described it.  This proposal was formulated in a strictly conformal setting. Although it is very interesting in itself, we cannot discuss it here in more detail.

\subsection{\small Mass generation  and Weyl geometric gravity ``at Munich'', 1980/90s\label{subsection Drechsler/Tann} }
\subsubsection{1990: Drechsler and Tann}
 A view  closer to Cheng's establishing a  connection between gravity and electroweak  fields in the framework of  Weyl geometry   was developed  a decade later   by {\em Wolfgang Drechsler} and his PhD student {\em Hanno Tann} at Munich. Drechsler had been active for more than twenty years in differential geometric aspects of field theory.\footnote{For example \citep{Drechsler/Mayer:1977}.}  In  cooperation with {\em D. Hartley} he developed an  approach of his own  to Weyl geometric gravity evolving form investigations in Kaluza-Klein  theories \citep{Drechsler/Hartley}.
Tann  joined the activity a little later during his work on his PhD thesis \citep{Tann:Diss},  coming from a background interest in geometric properties of the de Broglie-Bohm interpretation of quantum mechanics (see section \ref{subsection Bohm}). In their  joint work \citep{Drechsler/Tann}, as well as in their separate publications \citep{Tann:Diss,Drechsler:Higgs}  Weyl geometric structures were used in a coherent way, clearer than in  most of the other  physical papers discussed up to now. 

Tann studied  a complex valued scalar field $\Phi$, Drechsler, and the common paper of both, investigated a  scalar field with values in an isospin $\frac{1}{2}$ representation of the {\em ew}  group (like Cheng) with gravitational Lagrangian
\beq  \mathcal{L}_{\text{grav}} =  \mathcal{L}_{R} +  \mathcal{L}_{R^2} \label{Drechsler/Tann Lagrangian}
\eeq
where $ \mathcal{L}_{R^2} = \tilde{\alpha } R^2  \sqrt{|det\, g|}  $ and $ \mathfrak{L}_{R}= \frac{1}{12} \Phi^{\ast}\Phi\, R $ \citep{Drechsler:Higgs}. 
A common form of their linear gravitational  Lagrangian with modified Hilbert term $ \mathfrak{L}_{R}$ and the kinetic term of the scalar field is 
\beq \mathfrak{L}_{R,\Phi}= \frac{\beta}{2} \Phi^{\ast}\Phi\, R + \frac{1}{2}\,(D_{\nu}\Phi)^{\ast}\, D^{\nu}\Phi \, , \qquad \beta = \frac{1}{6}\; ,  \label{Lagrangian D/T}
\eeq
 with $\Phi^{\ast}$  the adjoint  (often written as $\Phi^{\dag}$) which in the case of Tann reduces to complex conjugation (often $\overline{\Phi}$), $R$ the Weyl geometric scalar curvature, signature of $g$ $(1,3)\sim(+---)$ and $D_{\nu}$ the scale covariant derivation, in Drechsler's case extended to the electroweak bundle.\footnote{\citep[equ. (372)]{Tann:Diss}, \citep[equ. (2.29)]{Drechsler:Higgs}. Both authors used coefficients like in the case of conformal coupling in Riemannian geometry, $\beta=\frac{1}{6}$.
In the Weyl geometric framework  this was an unnecessary restriction, because scale covariance holds for any $\beta$. In addition, Tann wrote the modified Hilbert term with a negative sign,  because the used the sign inverted convention for the Riemann tensor, see fn. \ref{fn curvature conventions}. } 
In such a Lagrangian they tried to straddle  the gap between the gravitational scalar field and a Higgs-like scalar field of electroweak theory. 

Both authors arrived  at a scale covariant expression for  the (metrical) energy momentum tensor  of the scalar field ($w(T)=-2$)  including terms, here with  factors  $ \beta^{-1} $, which result from varying the scale invariant Hilbert-Einstein term:\footnote{ \citep[equ. (372)]{Tann:Diss}, \citep[equ. (2.46)]{Drechsler:Higgs}.}
 \beqarr   T_{\Phi} &=&    D_{(\mu } \Phi ^{\ast} D_{\nu )}\Phi - \beta^{-1} D_{(\mu } D_{\nu )} |\Phi |^2  
  \label{energy-momentum phi}  \\
& & \quad \quad  - g_{\mu \nu }
 \left( \frac{1}{2}D^{\lambda }\Phi ^{\ast}D_{\lambda} \Phi -   \beta^{-1} \, D^{\lambda} D_{\lambda}(\Phi ^{\ast} \Phi)  + V(\Phi) \right)
\;  \nonumber
\eeqarr
 Drechsler noticed that the  $\beta^{-1}$ terms are identical to those  introduced in \citep{Callan/Coleman/Jackiw} for  ``improving'' the  ``energy-momentum tensor''  of a scalar field by  quantum physical considerations.

 In their common paper, Drechsler and Tann introduced fermionic Dirac fields into the analysis of Weyl geometry  \citep{Drechsler/Tann}. Their  gravitational Lagrangian had the form (\ref{Drechsler/Tann Lagrangian}).\footnote{In the appendix Drechsler and Tann showed that the squared Weyl geometric conformal curvature $C^2 = C_{\lambda \mu \nu \rho }C^{\lambda \mu \nu \rho }$ arises from the conformal curvature of the Riemannian component $_g\hspace{-1pt}C^2$  by adding a scale curvature term: $C^2 = \, _g\hspace{-1pt}C^2 +\frac{3}{2} f_{\mu \nu }f^{\mu \nu }$ \citep[(A 54)]{Drechsler/Tann}. So one may wonder, why they did not replace the square term $ \mathcal{L}_{R^2}$ by the Weyl geometric conformal curvature term   $\mathcal{L}_{\text{conf}} = \tilde{\alpha } C^2  \sqrt{|det\, g|}$.  }
 For the development of a Weyl geometric theory of the Dirac  field, 
they  introduced an adapted Lagrangian
\beq  \mathcal{L}_{\psi} =  \frac{i}{2} \left (\psi^{\ast}  \gamma ^{\mu}D_{\mu}\psi - D^{\ast}_{\mu}\psi^{\ast}  \gamma ^{\mu} \psi \right) + \gamma  |\Phi|  \psi^{\ast}\psi \, \label{Lagrangian Dirac field} \eeq  
with (scale invariant) coupling constant $ \gamma $ and Dirac matrices $\gamma^{\mu}$ with symmetric product $\frac{1}{2} \{ \gamma^{\mu},\gamma^{\nu} \} = g^{\mu \nu} \mathbf{1} $ \citep[(3.8)]{Drechsler/Tann}. Here the covariant derivative had to be  lifted to the spinor bundle, It  included  an  $U(1)$ electromagnetic potential $A=(A_{\mu})$, 
\beq  D_{\mu} \psi = \left( \partial _{\mu} + i \tilde{\Gamma} _{\mu} + \frac{i q}{\hbar c} A_{\mu}   \right) \psi \; ,  \eeq 
$q$ electric charge of the fermion field, $w(\psi) =-\frac{3}{2}$, $ \tilde{\Gamma} $ spin connection lifted from the Weylian affine connection.\footnote{(\ref{Lagrangian Dirac field}) can  equivalently be written with a  Weylianized scale covariant derivative  $\overline{D}_{\mu}  = \left( \partial _{\mu} + i \tilde{\Gamma} _{\mu}  + w(\psi) \varphi_{\mu}  + \frac{i q}{\hbar c} A_{\mu}   \right)$. Because $\varphi_{\mu}$ is real, the scale connection terms  $ w(\psi) \varphi_{\mu}$  in the Lagrangian  cancel.} 
This amounted to a   (local) construction of a spin $\frac{1}{2}$ bunde. Assuming the underlying spacetime $M$ to be  spin, they worked in a Dirac spin bundle $\mathcal{D}$ over the Weylian manifold $(M, [(g,\varphi )]$. Its structure group was  $G= Spin(3,1) \times R^+ \times U(1) \cong Spin(3,1) \times \C^{\ast}$, where $ \C^{\ast}=\C \setminus {0}$.\footnote{One could then just as well consider a complex valued connection $z= (z_{\mu})$ with values $z_{\mu}= \varphi_{\mu} + \frac{i}{\hbar c} A_{\mu}$  in   $\C= \mathfrak{Lie}(\C^{\ast})$
and weight $W(\psi)= (-\frac{3}{2}, q)$. Then $D_{\mu} \psi= (\partial _{\mu} + \tilde{\Gamma}_{\mu} + W(\psi)z_{\mu})\psi$, presupposing an obvious convention for applying $ W(\psi) z$.  }

The two authors considered (\ref{Lagrangian Dirac field}) 
as Lagrangian of a ``massless'' theory, because the masslike factor of the spinor field $\gamma  |\Phi|$ was  scale invariant,\footnote{This argument is possible, but not compelling  $\gamma  |\Phi|$ has the correct scaling weight of mass and may be considered as such.} 
and  proposed to  proceed to a theory with masses by   introducing a ``scale symmetry breaking'' Lagrange term 
\beq  \mathcal{L}_{B} \sim \frac{R}{6} + (\frac{mc^2}{\hbar})^2 |\Phi|^2  \label{Lagrangian Weyl symmetry break}  \eeq
with fixed (non-scaling) $m$ \citep[sec. 4]{Drechsler/Tann}.\footnote{ Similar already in Tann's PhD dissertation.}
But they did not associate such a transition from a (seemingly) ``massless'' theory to a massive one with any kind of hypothetical  ``phase transition''.

At the end of the paper they even commented:
 \begin{quote}
It is clear from the role the modulus of the scalar field plays in this theory (\ldots) that the scalar field with nonlinear selfcoupling is not a true matter field describing scalar particles. It is a universal field necessary to establish  a scale of length in a theory and should probably not be interpreted as a field having a particle interpretation. \citep[1050]{Drechsler/Tann}
\end{quote}
Their  interpretation of the scalar field  $\Phi$ was rather geometric than that of an ordinary quantum field; but their term
 (\ref{Lagrangian Weyl symmetry break}) looked ad-hoc to the uninitiated.\footnote{Note that one could just as well do without  (\ref{Lagrangian Weyl symmetry break}) and proceed with fully scale covariant masses -- compare last footnote. }  

\subsubsection{Drechsler on mass acquirement of electroweak bosons}
Shortly after the joint article with Tann,  Drechsler  extended the investigation to a  gravitationally coupled electroweak theory \citep{Drechsler:Higgs}.
Covariant derivatives were   lifted  as $\tilde{D}$ to the electroweak bundle. It included the additional connection components and coupling coefficients $ {g}$ and $ {g}'$ with regard to $SU(2)$ and $U(1)_Y$ like in Cheng's work (\ref{ew covariant derivative}).
The Weyl geometric Lagrangian  could be generalized and transferred to the electroweak bundle  \citep[(2.29)]{Drechsler:Higgs},
\beq \mathcal{L}  =   \mathcal{L}_{\text{grav}}  + \mathcal{L}_{\Phi } + \mathcal{L}_{{\psi}} + \mathcal{L}_{\text{YM}} \, ,  \label{total Lagrangian ew-grav} \eeq
with contributions
like in (\ref{Drechsler/Tann Lagrangian}),  (\ref{phi term Weyl geometric}), (\ref{Lagrangian Dirac field}),  and (\ref{Cheng's ew interaction Lagrangian}) ({\em ew} terms only).
 Lagrangians for the fermion fields had to be rewritten similar to electromagnetic Dirac fields (\ref{Lagrangian Dirac field}) and were decomposed into the chiral left and right contributions. 

In principle, Drechsler's proposal coincided with Cheng's; but he proceeded with more care and with  more detailed explicit  constructions. 
He derived the  equations of motion with respect to all dynamical variables \citep[equs. (2.35) -- (2.41)]{Drechsler:Higgs} and calculated the energy-momentum tensors of all  fields ocurring in the Lagrangian. 

The symmetry reduction from the electroweak group $G_{\text{ew}}$ to the electromagnetic $U(1)_{\text{em}}$ could then be expressed    similar to the procedure in the  standard model. $SU(2)$ gauge freedom  allows  to chose a (local) trivialization of the electroweak bundle such that the $\Phi$  assumes the  form considered in the ordinary Higgs mechanism 
\beq \hat{ \Phi} \doteq  \left(  \begin{array} { c}   0  \\  \phi_o   \end{array} \right) \; ,
\eeq 
where $ \Phi_o$ denotes a real valued field, and ``$\doteq$'' equality in a specific gauge. $ \hat{\Phi} $ has the isotropy group $U(1)$ considered as  $U(1)_{\text{em}}$and was  called   the {\em electromagnetic gauge} of  $\Phi$.\footnote{In other parts of the literature  (e.g., the work of R\c{a}czka and Paw{\l}owski) it is  called ``unitary gauge'', cf. also  \citep{Flato/Raczka}. }

In two respects Drechsler  went beyond what had been done before. 
 He {\em  reconsidered the standard interpretation} of  symmetry breaking by the Higgs mechanism \citep[1345f.]{Drechsler:Higgs}. 
 And he   calculated the consequences of nonvanishing {\em electroweak curvature components} for  the 
{\em energy-momentum tensor } of the  scalar field    $\hat{\Phi}$ \citep[1353ff.]{Drechsler:Higgs}. 
With regard to the first point, he  made clear that he saw nothing compelling in the  interpretation of  symmetry reduction as ``spontaneous symmetry breaking due to a nonvanishing  vacuum expectation value of the scalar field''  \citep[1345]{Drechsler:Higgs}. He  analyzed the situation and came to the conclusion that 
the transition from our ${\Phi}$  to $\hat{\Phi}$  is to be regarded as a  ``choice of coordinates'' for the representation of the scalar field in the theory and has, in the first place, nothing to do with a ``vacuum expectation value'' of this field.\footnote{Mathematically spoken,  it is a change of trivialization of the $SU(2)\times U(1)$-bundle.}
\begin{quote}
 \ldots This choice is actually not a breaking of the orginal $\tilde{G}$ gauge symmetry [our $G_{\text{ew}}$, E.S.] but a different realization of it. (ibid.) \label{p Drechsler}
\end{quote}
He compared the stabilizer $U(1)_{\text{em}}$ of $\hat{\Phi}$ with the ``Wigner rotations'' in the study of the representations of the Poincar\'e group.
With regard to the second point, the energy-momentum tensor of the scalar field could be calculated roughly like 
 in the simpler case of a complex scalar field, (\ref{energy-momentum phi}).   Different to what one knew from the pseudo-Riemannian case,  the covariant derivatives $D_{\mu}\Phi$ etc. in (\ref{energy-momentum phi})  were then  dependent on  scale or $U(1)_{\text{em}}$ curvature. 

After breaking the Weyl symmetry by a Lagrangian of form (\ref{Lagrangian Weyl symmetry break})  (ibid. sec. 3), Drechsler 
 calculated the curvature contributions induced by the Yang-Mills potentials of the $ew$ group and its consequences for the energy-momentum tensor $ T_{\Phi}$ of the scalar  field.
 Typical contributions to  components of $ T_{\Phi}$ had the form of mass terms 
\beq m_W^2 W^{+  \, \ast}_{\mu} W^{- \,  \mu}, \; m_Z  Z_{\mu}^{\ast} Z^{\mu}  \, , \quad  \mbox{with} \quad  m_W^2 = \frac{1}{4} g^2 |\Phi_o|^2  , \;
m_Z^2 = \frac{1}{4} g_o^2 |\Phi_o|^2  \, , \label{Drechsler's mass terms}
\eeq
 $g_o^2 = g^2 +  g'^2 $, for the bosonic fields $W^{\pm}, Z$ corresponding to the generators $\tau _{\pm}, \tau _o$ of the electroweak group,  \citep[1353ff.]{Drechsler:Higgs}.\footnote{$W^{\pm}_{\mu} = \frac{1}{\sqrt{2}} (W^1_{\mu} \mp i W^2 _{\mu})$, $Z_{\mu}= \cos \Theta \,  W^3_{\mu} - \sin \Theta \,  B_ {\mu} $. } 
 They are identical with the mass expressions for the $W$ and $Z$ bosons in  conventional electroweak theory.
According to  Drechsler, the  terms (\ref{Drechsler's mass terms}) in $ T_{\Phi}$ indicate  that 
the`` boson and fermion mass terms appear in the total energy-momentum tensor'' through the energy tensor  of the scalar field after ``breaking the Weyl 
symmetry''.\footnote{One has to be careful, however.  Things become more complicated if one considers the trace. In fact,  
 $ tr \, T_{\Phi}$ contains a  mass terms of the Dirac field of form $\gamma | \Phi_o | \hat{\psi}^{\ast} \hat{\psi}$, with $\gamma$ coupling constant of the Yukawa term
 ($\hat{\psi  }\;$  indicating electromagnetic gauge).  One of the obstacles for  making quantum  matter fields compatible with classical gravity is  the vanishing of  $tr \, T_{\psi}$, in contrast to the (nonvanishing) trace of the energy momentum tensor of  classical matter.  
Might Drechsler's analysis  indicate   a way out of this impasse? -- Warning: The mass-like expressions for $W$ and $Z$  in (\ref{Drechsler's mass terms}) cancel in $tr\, T_{\Phi}$ \citep[equ. (3.55)]{Drechsler:Higgs} like in the energy-momentum tensor of the $W$ and $Z$ fields themselves.}

Drechsler and Tann studied their scalar fields (complex or Higgs-like) as possibilies for an  extension of the gravitational structure of spacetime.  In their scale covariant theory of mass acquirement they tried to understand  how {\em mass generation}  is linked to the gravitational structure. Drechsler added that in his view the scalar field ``\ldots should probably not be interpreted as a field having a particle interpretation'' \citep[1050]{Drechsler/Tann}. This was an interesting remark at a time when elementary particle physicists started to collect information on a possible scalar boson of the Higgs field. But the empirical confirmation of the existence of a Higgs-like boson was still far out of sight; it did not materialize before  the LHC started to operate at a sufficient level of energy and luminosity in 2012.\footnote{See, e.g., \citep{Franklin:Higgs}.} Even so,  a more indirect link between the Higgs field and gravity, in contrast to the  perspective  of our two authors compatible with a bosonic interpretation of the scalar field, would be an interesting point. Back in the  1990s Drechsler did not expect that the  search for a bosonic quantum of Higgs type might ever be  confirmed by experiment.

\subsection{\small The ``Higgs'' and  Weyl scaling after 2000\label{subsection SM recent}}
In the years following the onset of the new millennium, but still before  a Higgs-like boson would be observed at the LHC, different  authors continued to explore the near to scale invariance of the standard model  and attempted to bridge the gap between the SM and gravity, keeping as closely as possible to  the original Higgs ``mechanism'' developed in the special relativistic framework. They did not adhere to a  common 
 research program; the researchers  used   different geometric/conceptual frameworks and worked in differing perspectives. Weyl geometric methods did not always stand in the center of the investigations; some scientists worked with global scale invariance and unimodular gravity 
\citep{Shaposhnikov_ea:2009,Shaposhnikov/Zenh"ausern:2009}, 
others preferred a conformal approach without making use of Weyl geometric concepts \citep{Meissner/Nicolai,Bars/Steinhardt/Turok:2014}, and some started from conformal symmetry but were mainly interested in models with radiative breaking of scale symmetry \citep{Foot_ea:2007,Foot_ea:2008,Foot_ea:2013}. A small group of authors, however,  continued in the line of Weyl geometric studies (Nishino/Rajpoot, Quiros, Ohanian e.a) . They often were  not aware of the whole range of  studies made in the 1970s to 1990s and took up just one filament  of the latter. 
 With few exceptions,\footnote{For an exception  still standing under the spell of Drechsler/Tann, although with  a consistently scale covariant  approach without an explicit  scale symmetry breaking term, see , e.g., \citep{Scholz:2011Annalen}.}
the majority of the mentioned authors worked with two scalar fields, a {\em Higgs-like} one $\Phi$ with values in a spin $\frac{1}{2}$ representation of the electroweak group, and a real-valued one, here denoted by  $\phi$,  the {\em gravitational scalar field}. Keeping track with our main theme, we shall  concentrate on the last group of authors who worked in the framework of Weyl geometry. 

The modified Hilbert-Weyl term  and kinetic terms of the scalar fields of these authors were, up to notational conventions,  of the common form
\beqarr L_{HW} &=& -\frac{\epsilon_{sig}}{2}(\zeta_1 \phi^2 +\zeta_2 \Phi^{\dag}\Phi)R\, , \qquad \label{L HW}\\
L_{\phi}&=& \epsilon_{sig} \frac{\alpha_1}{2}  D_{\nu}\phi D^{\nu}\phi\, , \qquad 
L_{\Phi}= \epsilon_{sig} \frac{ \alpha_1}{2} (D_{\nu}\Phi) D^{\nu}\Phi^{\dag}\, , \label{L phi}\\
\mathfrak{L} &=& L \sqrt{|g|} \; ,  \qquad \quad \epsilon_{sig} = \left\{ {+1 \quad \mbox{ for}\;  sig\, g = (+ - - - ) } \atop  -1 \quad \mbox{ for} \; sig\, g = (- + + + )   \right.  
\eeqarr 
(in most cases  $\alpha_1=\alpha_2=1$),  with Weyl geometric scalar curvature $R$ and  electroweak and Weyl geometric covariant derivatives $D_{\mu}$.\footnote{In the high energy physical context, and accordingly in our section 3,  signature of $g$ 
$=(+ - - -) $. For sign  conventions regarding curvature see fn \ref{fn curvature conventions}.} A quadratic  curvature term $L_{R^2}$ was added by some, not by all, authors. 
Yang-Mills terms of  the electroweak connections (potentials) $W$ for the $SU_2$-component, $B$ for hypercharge $U(1)$, and $\varphi$  for the scale connection with field strength  $f=d\varphi$,  were added,
\beq L_{YM}=-\frac{1}{4}\left( tr(W_{\mu \nu}W^{\mu \nu}) + B_{\mu\nu}B^{\mu \nu} +f_{\mu\nu}f^{\mu\nu} \right) \, . \label{L YM}
\eeq 
Similarly Dirac kinetic terms  $L_{\Psi\, kin}$ and Yukawa mass terms $L_{\Psi\, Y}$ for the different  fermions $\Psi^{f g}_{i h}$, with indices taking care for the various types and properties ($f=q, l$ for quark or lepton, $g=1,2,3$ generation, $i= u, d$ (``up, down'') for the {3}-component of weak isospin, $h=R,L$ helicity) were added in a form adapted to the Weyl geometric framework. The Dirac terms could be written with or without the Weyl geometric scale connection term because, even if it is included, it finally cancels in the total expression. This had been noticed  already by Hayashi and Kugo  (see section \ref{subsection Dirac-Utiyama}).\footnote{See also \citep[p. 81]{Blagojevic:Gravitation}.} We need not reproduce the explicit form of the fermionic terms here, but have to keep in mind that the Yukawa terms contained the  matrices with relative mass coefficients (``mass matrix'') of the SM and a  {\em scale covariant} Higgs field.\footnote{For an explicit form of Dirac kinetic terms and Yukawa mass terms see,  e.g., \citep[equ. (1.2)]{Nishino/Rajpoot:2009}.}  

Breaking of  scale invariance without an explicit mass terms of the Higgs field became the crucial points for our authors. Because of its scaling behaviour ($w(\phi)=-1$) the gravitational scalar field $\phi$ already {\em specifies a preferred scale}  in which it assumes a constant value $\phi_o$ (scalar field gauge in the terminlology of section \ref{subsection Weyl geometry}):
\beq \phi(x) \doteq \phi_o = const \quad \label{scalar field gauge}
\eeq 
This was the reason behind  Utiyama  calling  $\phi$  a ``measuring field'' already in the 1970s.

But  the question still remains  how such an, at first sight only mathematical, specification may be incorporated in the material structures lying at the basis of measuring processes. In the context of the search of a connection between gravity and  the {\em ew} sector of fundamental fields it seemed natural to search for a relation between the two scalar fields $\phi$ and $\Phi$. For this 
 a biquadratic/quartic potential  in the two scalar fields,  and a corresponding Lagrange term,  plays a crucial role. Using the abbreviation $|\Phi|^2=\Phi^{\dag}\Phi$ it is:
\beqarr V(\Phi, \phi) &=&   \frac{\lambda_1}{4}|\Phi|^4-\frac{\mu}{2}|\Phi|^2\phi^2+ \frac{\lambda'}{4}\phi^4  \label{L V4} \\
 &=&\frac{\lambda_1}{4}\left(|\Phi|^2 - \frac{ \mu}{\lambda_1}\phi^2\right)^2+ \frac{\lambda}{4}\phi^4 \, , \qquad \lambda = \lambda'-\frac{\mu^2}{\lambda_1} >0 \, ,\nonumber \\
\mathcal{L}_{V} &=& -V(\Phi, \phi) \sqrt{|g|}
\eeqarr 
 Chromodynamics was usually not considered;   our group of authors concentrated on the electroweak sector of the SM and its possible link to gravity.

\subsubsection{Nishino/Rajpoot\label{subsection Nishino/Rajpoot}}
 In 2004 two theoretical high energy physicists at California State University, {\em Hitoshi Nishino } and {\em Subhash Rajpoot} posed the goal of ``extending the standard model with Weyl's scale invariance'', adding that the scale invariance is ``badly broken'' at the order of the Planck mass/energy. They made it clear that in the ``philosophy advocated in the present work \ldots the standard model Higgs is not eliminated, and is the sought for particle''  \citep[1]{Nishino/Rajpoot:2004}. 

For adapting the fermionic fields to the differential geometric setting, the authors outlined the usual spinor calculus in a Weyl geometric approach with scale dependent tetrads consisting of  point-dependent bases  $e_a =e_a^{\mu} \,\partial_{\mu} \; (a=0, \ldots 3)$ and their dual forms, here denoted by $\vartheta^a= \vartheta^{a}_{\mu}dx^{\mu}$, and the metric $g_{\mu \nu}= \vartheta^a_{\mu}e_{a\,\nu}$.
With  $g(x) \mapsto \tilde{g}=e^{2\Lambda(x)}g(x)$ the tetrads have to be rescaled like
\beq \vartheta^a_{\mu} \mapsto \tilde{ \vartheta^a_{\mu}} =  e^{\Lambda(x)} \vartheta^a_{\mu} \, \qquad  e_a^{\mu} \mapsto    \tilde{ e}_a^{\mu}  =  e^{-\Lambda(x)}  e_a^{\mu} \, ,
\eeq
that is $w(\vartheta^{a})=1, \; w(e_a)=-1$. The Weyl geometric affine connection,  the corresponding spin connection, Weyl geometric covariant derivatives, and curvature expressions were developed by the two authors, although not always completely reliable.\footnote{The expression for the scalar curvature is given in the paper (and also in the later papers by the same authors) as $R=\,  _g R{}\hspace{-0.2em} - 6 \nabla_{\mu}\varphi^{\mu}+ 6\varphi_{\mu}\varphi^{\mu} $, where a coupling constant $f$ introduced by the authors is here set to $f=1$ and transcribed into  our notation,  \citep[equ. (14)]{Nishino/Rajpoot:2004}. The correct Weyl geometric value (\ref{Weylian scalar curvature}) would be  $R=\,  _g{}\hspace{-0.2em} R - 6\, _g{}\hspace{-0.2em}\nabla_{\mu}\varphi^{\mu}- 6\varphi_{\mu}\varphi^{\mu}  $, cf. \citep[p. 21]{Weyl:InfGeo}, \citep[equ. (A 31)]{Drechsler/Tann} and others. Because of $ \nabla_{\mu}\varphi^{\mu}= \,  _g{}\hspace{-0.2em}\nabla_{\mu}\varphi^{\mu} + 4 \varphi_{\mu}\varphi^{\mu}$ this implies   $ R= \, _g{}\hspace{-0.2em} R - 6 \nabla_{\mu}\varphi^{\mu}+ 18\varphi_{\mu}\varphi^{\mu}$ ! \label{fn error R}}

On this background they described a {\em two stage process} of symmetry breaking. In the  first step they dealt with breaking  the scale symmetry,  formulated in terms of  the compactified scaling   group $\tilde{U}(1)$. The breaking was expressed   ``by setting'' the value of the gravitational scalar field to a constant $\phi_o$
\beq \phi(x)= \phi_o    \quad \mbox{with}\quad \zeta_1 \phi_o^2 = (8 \pi G)^{-1} \,,
\eeq 
in our terminology they introduced scalar field (Einstein) gauge.
In the second step the  {\em ew} symmetry was assumed to be broken ``spontaneously'' like in the special relativistic SM case ($SU_2\times U(1)_Y \mapsto U(1)_{em}$). For the  first step they gave a physical interpretation  which has some analogy with the Higgs ``mechanism'':
\begin{quote} At this stage the scalar field $\sigma$ [here denoted  $\phi$, E.S.] becomes the Goldstone boson \ldots. The vector particle associated with $\tilde{U}(1)$ breaking, the Weylon, absorbs the Goldstone field and becomes massive with mass  $M_S$ given by $M_S  = \sqrt{\frac{ 3f^2}{4 \pi G_N}} \approx 0.5\times f\, M_P$ [$f$ a coupling constant of the scale connection, E.S.]. \citep[4]{Nishino/Rajpoot:2004}
\end{quote}
Then the quartic potential (\ref{L V4}) is reduced to the Higgs potential like in the SM plus a cosmological term $\frac{\lambda'}{4}\phi_o^4$. In the ground state of the Higgs field only the cosmological term survives and the transition to scalar field gauge endows the Higgs field with  mass
\beq m_H  \doteq \sqrt{\mu} \phi_o \, . \label{Higgs mass}
\eeq

After  a short outline of how to adapt the parameters to the mass generation scheme of the SM the authors concluded
\begin{quote}
Our contention is that the present model presents a viable scheme in which gravity is
unified, albeit in a semi-satisfactory way, with the other interactions. (\ldots) When the
complete theory of all interactions is found, the model in its present form, it is hoped, will
serve as its low energy limit.

To conclude, we have accommodated Weyl’s scale invariance as a local symmetry in the
standard electroweak model. This inevitably leads to the introduction of general relativity. \citep[8]{Nishino/Rajpoot:2004}
\end{quote}

This paper  remained in a preprint stage. Although  its content seems to have been  presented at different  conferences it never was published in a scientific journal. The reason may have been that the authors considered it only  as a first, provisional step. In the following years they extended their approach to a $SU(5)$ grand unified theory (GUT) \citep{Nishino/Rajpoot:2007} and revised their presentation by taking up an idea  going back to Stueckelberg \citep{Nishino/Rajpoot:2009,Nishino/Rajpoot:2011}.  

 In the late 1930s {\em Ernst  Stueckelberg} had introduced a massive scalar field $B$ complementing an  $U(1)$ potential $A_{\mu}$ which expressed a field of electromagnetic type, but with mass (i.e. similar to a Proca field). $B$ was given a peculiar gauge behaviour involving a mass parameter $m$ under $U(1)$ gauge transformations $A_{\mu} \mapsto \tilde{A}_{\mu}= A_{\mu} + \partial_{\mu} \Lambda(x)$
\beq B(x) \mapsto \tilde{B} = B(x) + m\, \Lambda (x) \, . \label{Stueckelberg transformation}
\eeq 
Stueckelberg's context was the search for an  interaction of a scalar field with nucleons.  
Transformations of type (\ref{Stueckelberg transformation}) were taken up by Pauli and others. They  became to be known as {\em Stueckelberg transformations} and   $B$ as {\em Stueckelberg (compensating) field}. With an  appropriate $\Lambda$, the Stueckelberg field allowed to specify  a peculiar gauge with  $\tilde{B}=0$,  {\em without breaking} the $U(1)$ symmetry  which is only  given a ``different realization'' (in Drechsler's terms  quoted above, p. \pageref{p Drechsler}). This  turned out to be crucial for the  renormalizability of the theory and made the ``Stueckelberg trick''  attractive for  quantizing  the electromagnetic field or its relatives like Proca like fields.\footnote{The non-broken  $U(1)$ symmetry is important for the BRST relations, the quantum analogue of the Noether relations.  See \citep[75ff.]{Ruegg_ea:Stueckelberg}.} 

The careful reader may have noted the kinship between the Stueckelberg ``trick'' for $U(1)$ and the Higgs ``mechanism'' for the electroweak group. So did Nishino/Rajpoot. Moreover, they realized that,  just by taking the logarithm,  the transition to the Weylian scalar field gauge  can be given the  form of a Stueckelberg transformation. Transliterated to our notation they introduced an exponential expression of the form\footnote{The  factor $\zeta_1^{-\frac{1}{2}}$ in \citep[equ. (2.1)]{Nishino/Rajpoot:2009} was set by them to $\zeta_1=1$ while  transforming the Lagrangian  into their equ. (2.3).  The follow up paper \citep[ second paragraph of section 2]{Nishino/Rajpoot:2011} shows that this reduction was intended. Of course, a different factor $\zeta_1$  would heavily influence the  mass calculation in (\ref{mass contribution kinetic term}).} 
\beq   \phi(x) = \zeta_1^{-\frac{1}{2}}M_P\, e^{M_p^{-1} \beta(x) }\, .
\eeq
Then  the scale gauge transformation $\phi \mapsto \tilde{\phi} = e^{-\Lambda}\phi$ is expressed by 
\beq \beta \mapsto \tilde{\beta} + M_p \Lambda  \, , \label{Stuckelberg-Nishino transformation}
\eeq
and the transition to scalar field gauge corresponds to $\tilde{\beta}=0$, exactly like in the case of the the Stueckelberg ``trick''.

Nishino/Rapoot thus  rewrote their basic Lagrange density equivalent to our equations (\ref{L HW}, \ref{L phi}, \ref{L YM}, \ref{L V4}) in terms of the logarithmized scalar field \cite[equ. (2.3)]{Nishino/Rajpoot:2009} and  normed it to scalar field gauge \citep[equ. (2.6]{Nishino/Rajpoot:2009}. Then  the mass expression $m_{\varphi}$ for the scale connection field (``Weyl field'')  could be read off.\footnote{Warning: Nishino/Rajpoot used the notation $\varphi$ for the Stueckelberg ``compensator'', i.e. our $\beta$, and $S_{\mu}$ for the scale connection (the potential of the ``Weyl field''), our $\varphi_{\mu}$. In order to avoid confusion the notation in the present paper has been homogenized for  the authors discussed here.}
In scalar field gauge  the kinetic terms   (\ref{L phi}) of $\phi$ and $\Phi$ acquire  forms which makes them  contribute to $m_{\varphi}$.  For $\phi$ it is
\beq \frac{1}{2}D_{\nu}\phi D^{\nu}\phi \doteq \frac{1}{2} (f M_p)^2 \varphi_{\nu}\varphi^{\nu} \, , \label{mass contribution kinetic term}
\eeq 
while for $\Phi$ the   contribution to the mass of $\varphi$ is  $f(\Phi^{\dag}\Phi)$  (after {\em ew symmetry breaking} $f v^2$, with $v^2$ the vacuum expectation value of the operator $\Phi^{\dag}\Phi$).
In any case  the contribution  due to $\Phi$ is 
 much less than the one from $\phi$ and from the modified Hilbert term (\ref{L HW}), both of which are at the order of  the Planck scale. It  may safely be neglected at several orders of magnitude.\footnote{Nishino/Rajpoot did not consider the contribution of the modified Hilbert term, in contrast to  Smolin and Cheng  (see section \ref{subsection SM 1970s}).}

In the imaginative language of   the elementary particle community Nishino and Rajpoot  commented that the scalar field is ``now eaten up by the Weylon''. A little later they added,  more technically:
\begin{quote}  After all, the Weylon $\check{S}_{\mu}$ [our $\varphi_{\mu}$, E.S.] acquires the mass $f M_P$, the compensator $\varphi$ [our $\beta$, E.S.] is absorbed into the longitudinal component of $\check{S}_{\mu}$, and the potential terms are reduced to the Higgs potential in SM \ldots \citep[3]{Nishino/Rajpoot:2009}
\end{quote}
With this explanation they clad the mass derivation for the scale connection field in the mantle of a narrative which is   widely spread  in their community  and usually accepted as  scientifically explanatory.\footnote{Compare \citep{Stoeltzner:stories}.}

In a follow up paper, the two California State physicists came back to the topic and presented the results of some results concerning a quantized version of their theory. 
They started from their Lagrangian given in  terms of the logarithmized scalar field \cite[equ. (2.3)]{Nishino/Rajpoot:2009} and with  modified Hilbert term
\beq L_{HW} = -\frac{1}{2}\left(\zeta_1 M_P^2\, e^{2 M_p^{-1} \beta(x) }+\zeta_2 \Phi^{\dag}\Phi\right) \, R
\eeq
($R$ Weylian scalar curvature written there as $\tilde{R}$). 

At this point  Nishino and Rajpoot  left the  track of Weyl geometry and decided to switch to  the   JBD paradigm. They  considered the initial  Lagrangian a  ``Jordan frame'' and wanted to transform it to ``Einstein frame''.\footnote{Strictly speaking their framework does not contain any meaningful ``Jordan frame'', because their Weyl structure is not integrable, and thus the purely Riemannian representation of the affine connection presupposed in ordinary Jordan frame does not exist. Einstein frame, on the other hand, is meaningful in any Weyl geometric gravity approach with a scale covariant scalar field and corresponds to scalar field gauge (\ref{scalar field gauge}). }
For this goal they performed  a ``Weyl rescaling for the vierbein or metric   only'', i.e.  a ``{\em field re-definition}'' which did not include the corresponding transformations of the scale covariant fields and the scale connection.\footnote{``Note that the Weyl rescaling we made is a field re-definition, but it is not a part of any local scale transformation which is defined to act not only on $g_{\mu \nu}$ but also on $\Phi$ and $\varphi$
as in (2.2) [the equation for the full gauge transformation, E.S.]'' \cite[p. 4]{Nishino/Rajpoot:2011}. \label{fn field re-definition}}
Referring to calculations in the framework of JBD theory they arrived at a reduction of the Hilbert term to a form which depends only on the Riemannian component $_g R$ of the scalar curvature. According to their calculation,  the scale connection contributions drop out of the Lagrangian (but not the Yang-Mills term for the scale curvature ).\footnote{In the light of the  error for the scalar curvature indicated in fn. \ref{fn error R} one may be  inclined to doubt   the correctness of such a complete cancellation.}
 In other word,  a reduction to Einstein frame form of JBD with two scalar fields and an additional Yang-Mills field was achieved \cite[equ. (2.10)]{Nishino/Rajpoot:2011}. 

On this basis our authors performed a series of calculations at the quantum level. They determined   (Adler-Bell-Jackiw and trace) anomalies, studied the possibility for cancelling  the remaining (trace-) anomalies, considered quantum corrections to the cosmological constant, and studied the perturbative renormalizability of their model and the possible new divergences. All in all, these were remarkable results; but they were  arrived at in a hybrid  approach which  started  in a setting of Weyl geometric gravity and ended  in JBD gravity, after performing an  artificial and methodologically unconvincing transition by a  ``field re-definition'' type of rescaling. In spite of such  shortcomings the derivations  were a notable step towards connecting  the electroweak sector of elementary particle physics with gravitational structures, mainly formulated in a Weyl geometric framework. 

\subsubsection{Hao Wei, Rong-Gen Ca, Quiros\label{subsub Quiros ea}}
 H. Nishino and S. Rajpoot  were not the only  researchers who thought about the question how to establish a connection between gravity and the SM fields by exploiting Weyl geometric methods. Even though  we have to be selective here,  it  has to be clear that the Weyl geometric approach continues to be  alive in the era of the Higgs boson (or some close relative) being found in experimental observations. 
 A talk given   in July 2004   by  Hung Cheng at the Institute for Theoretical Physics of the Chinese Academy of Science, Beijing seems to have initiated interest in  Weyl geometric methods  by Chinese theoretical physicists {\em Hao Wei, Rong-Gen Cai} and others.\footnote{\cite[Acknowledgments]{Cai/Wei:2007}} 
It was natural for them to take the ``Cheng-Weyl vector field'' (i.e., the Weylian scale connection with massive boson studied by Cheng in the late 1980s) and Cheng's view as their starting point for a new look at   the standard model of elementary particle physics \cite{Wu:2004,Cai/Wei:2007}.

Another road was taken by  {\em Israel Quiros},  at the time we are interested in here, placed at Guanajuato, Mexico.  Coming from a background in Jordan-Brans-Dicke gravity and cosmology (see section \ref{subsection trivial}) he developed thoughts of his own about how ``scale invariance and broken electroweak symmetry may coexist together''  \citep{Quiros:2013}.  In this   conceptually clear paper he gave a nice introduction  to the basic  ideas  of  integrable Weyl geometry and showed that the scale covariance of the SM fields can not only  be imported into a general relativistic framework, if Weyl geometric gravity is used, but  can even be upheld  after  breaking the {\em ew} symmetry. One only need to accept, and to use,  mass parameters $m$  which scale with weight  $w(m)=-1$. 

For his presentation Quiros  used a simplified version of the Lagrangian (\ref{L HW}ff.)  similar to the one of Nishino/Rajpoot, whose papers he probably did not yet know. He encoded the gravitational scalar field in terms of a point-dependent scalar  exponent written by him  as $\varphi$  -- in order to avoid confusion we shall transliterate it like above  as $\beta$ -- of  the  factor in the  Hilbert-Weyl term. Compared with our notation above he wrote 
\beq \zeta_1\phi(x)^2 = M_p e^{\beta(x)} \label{simplification}
\eeq
 and considered Weylian scale connections exclusively of the form
\beq   \varphi =\varphi_{\mu} \,dx^{\mu} = d\beta \quad  \Longleftrightarrow \quad \varphi_{\mu} = \partial_{\mu}\beta
\eeq
\citep[equ. (8)]{Quiros:2013}. This implies the restriction
\beq \phi = const \Longleftrightarrow \beta=0  \Rightarrow  \varphi=0\; .
\eeq 
In our terminology (\ref{simplification})  implies an inbuilt {\em identification of   Riemann gauge and Einstein gauge}. That was probably unnoticed by the author, and is widely spread among scientists who  entered into  Weyl geometric methods from a JBD background. For the basically geometrical and conceptual, purposes of the paper this restriction  may have been of no particular disadvantage, but    the dynamical role of the scalar field was trivialized by this specialization. 

\subsection{\small  Towards Weyl scaling at the quantum level\label{subsection quantum level}}
\subsubsection{Scale invariant quantization procedures\label{subsection Codello et al.}} 
Problems on a more fundamental have been posed by a 
 group of  theoretical physicists    working at Trieste.  {\em  Alessandro Codello, Giulio D'Orico, Carlo Pagani} and {\em Roberto Percacci}  recently reconsidered the question of how scale invariance behaves under  quantization if one approaches it with the method of the so-called  ``renormalization group'' (RG)   and the use of  functional integral methods.  
 In   \citep{Codello_ea:2013} they gave a  report  on their work and  rebutted   the general view that quantization  necessarily leads to a breaking of (point-dependent) scale symmetry even if the classical Lagrangian  is scale invariant. In a step by step argumentation they show  how the 
 functional integrals can be given a scale invariant form by  using an integrable Weyl geometric background and a gravitational scalar field $\chi$  of weight $w(\chi)=-1$, called a ``dilaton'',   as external fields which are not quantized at the first stage. 

They started from the basic idea  that ``one can make any action Weyl-invariant by replacing all dimensionful couplings by dimensionless couplings multiplied by the powers of the dilaton'' \citep[p. 2]{Codello_ea:2013}. 
Then a dimensional  coupling coefficient of scaling dimension $k$, let aus say $\mu$,  is turned into a coupling parameter of the form $\chi^{-k} \check{\mu}$ with a ``dimensionless'', i.e. non-scaling, constant $\check{\mu}$. The authors achieve scale covariance/invariance  of the fields, respectively actions, by using  Weyl geometric expression with regard to an integrable scale connection with coefficients
\beq b_{\mu} = - \chi^{-1} \partial_{\mu} \chi \qquad \mbox{\citep[p. 3]{Codello_ea:2013}\footnote{Note that the differential form $b= b_{\mu}dx^{\mu}$ is sign inverted in comparison with our conventions of section \ref{subsection Weyl geometry}. }}
\eeq 
Like Nishino and Rajpoot they consider  this as a gravitational equivalent to the ``St\"uckelberg trick''. 
 Their main work then consisted in showing that the Weyl invariance which is easily achievable for the classical action  is left  intact, in their framework,  for the functional integrals, the differential equation governing the renormalization flow equation, and the UV and IR endpoints of the flow. 

 Classical quantum matter fields  (scalar or Dirac spinors) have a vanishing trace of the energy-momentum tensor, while  the expectation value of the quantized  trace no longer vanishes. This so-called {\em trace anomaly} of quantization has  puzzled theoretical physicists for a long time and is usually taken as a sign that scale invariance is broken at the quantum level. Our authors came to a different conclusion. They explained  that, although   the ``trace anomaly'' is still present in their approach,  it {\em no longer signifies  breaking of the local scale invariance}. The reason  lies in a  cancellation of the trace  terms of the quantized fields a  by corresponding counter-terms arising from the scalar field, the ``dilaton'' in the language of the paper. 

After some comments on the quantization of the metric field, and further discussions of the difference between strictly conformal theories and the Weyl geometrically ``conformalized'' ones, the authors finished with the remark:
\begin{quote}
The present work provides a general proof that with a suitable
quantization procedure, the equivalence between conformal frames can also be maintained in the quantum theory \citep[p. 21]{Codello_ea:2013}.
\end{quote}
But they also stated clearly that  their quantization procedure does not lead to new physical effects. In this sense their research shows a certain analogy to Kretschmann's view of diffeomorphism invariant re-formulations of physical theories which do not  per se lead to new physical insights.

 Even so, the authors have achieved to show that  the  extension of the mathematical automorphism group of the underlying theories (SM fields, implicitly also gravity theory) can be upheld under quantization. Whether a further enrichment of the theories delivers new insights at the quantum level will be a question for the future. Probably this can only be the case, if the scalar field and/or the scale connection acquires a dynamical role beyond its purely mathematical  ``compensatory'' character in the  scale transformation.

\subsubsection{Ohanian's retake of a ``spontaneous'' breaking of symmetry\label{subsection Ohanian}}
An   attempt at giving the scale connection a dynamical role has been made  
  by {\em Hans Ohanian} from the University of  Vermont. He proposed a model which  connects the standard model fields with general relativity in a Weyl geometric framework.   A complex scalar field  $\chi$ (``dilaton'')  acts as the crucial mediator. It undergoes spontaneous breaking of local scaling symmetry which the author preferred to call conformal symmetry,\footnote{Ohanian  preserved the label ``scale transformation'' for a global usage in Minkowski space, where, in addition to the rescaling of the fields $X\mapsto \tilde{X} = \Omega^k X$, a space dilation $x\mapsto \tilde{x}= \Omega x$ is applied \citep[p. 25]{Ohanian:2016}.}
 by a mechanism very similar to the breaking of electrodynamic $U(1)$ symmetry in a model studied by \citep{Coleman/Weinberg:1973}.    If gravitational effects can be neglected,  Ohanian's adaptation leads  to the SM field content in flat spacetime \citep{Ohanian:2016}. If, on the other hand, gravity is taken into account,  the transition from quantum to classical matter being leapfrogged,   it leads to  Einstein gravity as an ``effective field theory''.   Regarding the conformal expression of fields  Ohanian used a ``conformalization'' procedure with additional terms in the (Riemannian) scalar curvature  (in place of the more natural Weyl geometric expressions).
 Ohanian   proposed to assimilate the result  of Coleman/Weinberg by a simple substitution of coefficients and concluded:
\begin{quote}
After symmetry breaking, neither
the scalar field nor the vector field reveal themselves at the macroscopic level, and we
can ignore the effects of the Weyl gauge-vector on the transport of lengths \ldots . \citep[10f.]{Ohanian:2016}
\end{quote}
Because of the conformal coupling of the scalar fields to the {\em Riemannian} scalar curvature Ohanian  found that in his approach a  {\em modification of Riemannian geometry is excluded} in the long-range regime and comments:
\begin{quote}
This is in contrast to the standard Brans-Dicke theory, in which the massless scalar field makes a contribution to long-range gravitational effects, \ldots (ibid.)
\end{quote}

In the   high energy, short-range, regime Weyl geometric curvature does play a role in this model, as Ohanian discussed  in his section 4. Then the scale connection constitutes a ``vector'' field of its own, similar to the electromagnetic field, but with a mass term and with the dynamical current of the scale symmetry $\mathfrak{J}^{\mu}=\frac{\partial \mathfrak{L}}{\partial \varphi_{\mu}}$ as  right hand side of the dynamical equation.\footnote{ $\partial_{\nu}\left(\sqrt{|g|}f^{\mu \nu}  \right) =\mathfrak{J}^{\mu}$. In Ohanian's Lagrangian  $\varphi$ couples only to the ``dilaton'' scalar field $\chi$. This  leads to a  form for the variation of the Lagrangian under  scale transformations such that the dynamical current  coincides with the Noether current   \citep[equ. (13)]{Ohanian:2016}.} 

In his outlook Ohanian conjectured that certain problematic  features in the purely conformal approaches are essentially due to the lack of a a coherent metrical structure. In Weyl geometry the scale connection is the clue for making a Weylian metric consistent with conformal rescaling. 
Ohanian therefore  finished his paper with a remark which went right to the heart of the matter:
\begin{quote}
If the analysis of Ehlers et al. is correct, the absence of a Weyl vector and its
geometric paraphernalia is a fatal mistake -- if no Weyl vector, then no conformally
invariant theory with a geometric interpretation \citep[p. 16]{Ohanian:2016}.
\end{quote}

In this approach   Ohanian  proposed a   model which indicated why and how a Weyl field with curvature at the short-range, high energy level  looses its curvature in the  low energy regime and leads to  Einstein gravity in the long-range limit.

 Ohanian, like many other authors,   perceived the transition between the energy  regimes (high -- low) exclusively in the sense of hypothetical {\em successive temporal stages} in the cosmic development. This fits in  with the mainstream narrative connecting  cosmology and high energy physics shortly after the big bang. Philosophically inclined reader may  notice that one could interpret such kind of transition non-temporally, as a  {\em  structural passage  between different energy levels}, present at any time and any place of the world. This would be independent of the view regarding the reality content of the big bang  picture.\footnote{Physicists may well claim  that, e.g.,  the LHC experiments are important because they explore  how the world  has looked like a few ``nanoseconds after the big bang''. But  one need not take such stories at face value  in order  to appreciate the activities aiming at gaining knowledge about the respective energy levels and  the transitions between them.}

\section{\small Weyl geometric models  in astrophysics and cosmology since the 1990s\label{section cosmology}}

\subsection{\small The broader context: scalar fields in gravity,  conformal rescaling\label{subsection cosmology conformal}}
In the 1970s JBD theory underwent a  contradictory development: On the one hand, increasing precision of   radar tracking observations in the planetary system  showed    that Einstein gravity is an extremely   good description of gravity.\footnote{See C.Will's contribution to this volume.}
A tentative modification of the latter by a  Brans-Dicke type scalar field has at least  to be  suppressed  on this level, e.g. by an extremely high value of the coupling coefficient $\xi$ of the kinetic term in (\ref{Lagrangian JBD}), or may it be not adequate at all. 
On the other hand, the rise of   particle cosmology as a new subfield of theoretical physics opened ample space for  studying  models  in an assumed   very early phase of the universe. Here it appeared reasonable to think about modified gravity  and  elementary particle physics as an ensemble.  A fertile   environment for   studying  speculative models  emerged, some of which  were designed   for  combining the gravitational scalar field and a Higgs-type scalar field of elementary particle physics  \citep{Kaiser:mass,Kaiser:colliding}. This environment gave  new motivations and  incentives  for studying  scalar-tensor theories, completely different from those of the 1960/70s \citep[chaps. 3, 7]{Capozziello/Faraoni}.  One of the new roles rehearsed for the scalar field on this stage was that of  an agent, called {\em inflaton}, which drives a hypothetical phase of very early accelerated expansion of the spacetime.

Another role arose from string theory where a new type of scalar field, a so-called {\em dilaton},   entered the stage. Originally it coupled to the trace of the (2-dimensional) stress tensor of the string. But in  the form of a  constraint for restoring conformal symmetry, after its breaking under quantization, the dilaton re-appeared  as a source term in a classical Einstein-like equation.   That gave rise to speculate about deriving Einstein gravity as an effective theory arising from string theory, with the dilaton scalar field and conformal symmetry as mediators \citep[p. 14f.]{Brans:roots}.

 All in all, a vast field for studying scalar field theories in generalized theories of gravity arose.\footnote{For extensive surveys of this field see \citep{Fujii/Maeda,Capozziello/Faraoni}.} Only few  authors of this field remembered Weyl geometry and took it up for their purpose. This was the case, e.g., in string models; but  they  remain outside the scope of this survey. They  would  need a  study of their own; here we look at more mundane manifestations of Weyl geometry in cosmology and astrophysics during the last two decades. Because of the close kinship between Weyl geometric rescaling and conformal invariance of field theories in a Riemannian environment  I here  bring  only a few  examples for recent conformal approaches in cosmology to the mind. They are far from exhaustive and have been selected because they connect in specific ways to our core topic.

\subsubsection*{Conformal approaches in cosmology}
An unusual analysis of  the ``dark'' sectors of recent cosmology  was given by {\em Philip Mannheim} and {\em Demosthenes Kazanas}.  They  argued that  the flat rotation curves of galaxies can be explained on the basis of a  conformal approach to gravity  \citep{Mannheim/Kazanas:1989}. In their  conformal theory, a static spherically symmetric matter distribution was described by the solution of a fourth order Poisson equation
\beq  \nabla ^4 B(r) = f(r) 
\eeq 
  with a typical coefficient  $B(r)$ proportional to $  - g_{oo}=  g_{rr}^{-1} $
of a metric $ds^2  = g_{oo}dt^2 -  g_{rr}dr^2 - r^2 d\Omega ^2$ (up to a conformal factor).    The r.h.s. of the Poisson equation, $f(r)$, depended on the mass distribution, e.g., in a spiral galaxy. The result of a comparison of their theory with data  for 11 galaxies with different behaviour of rotation curves led to a  good fit and encouraged the authors to present their approach as a possible candidate for a modified gravity explanation of dark matter phenomena  \citep{Mannheim/Kazanas:1989,Mannheim:1994}.

During the following years the approach was  extended to the question of dark energy in a peculiar perspective. In the special case of conformally flat models, like Robertson-Walker geometries, Mannheim proposed to consider the  Hilbert-Einstein Lagrangian term $-\frac{1}{12}|\phi |^2 R \sqrt{|det\, g|}  $ of a conformally coupled scalar field $\phi$ as part of the {\em matter} Lagrangian. Due to this sign choice, he arrived at a version of the Einstein  equation with {\em inverted } sign. He interpreted this as a  kind of ``repulsive gravity'' which supposedly  operates on cosmic scales in addition to the ``attractive gravity'' on smaller scales,  indicated by the conformally modified Schwarzschild solution. In his eyes, such a repulsive gravity might step into the place of the dark energy of the cosmological constant term of standard gravity \citep[729]{Mannheim:2000}.

In spite of such a grave difference to Einstein gravity, Mannheim did not consider his conformal view to   disagree  with the standard model of cosmology and its accelerated expansion. He rather argued that his approach may lead to a more satisfying explanation of the expansion dynamics. In his view, ``repulsive gravity'' would take over the role of dark energy. Moreover he expected that a conformal approach  with quadratic curvature terms may shed new light on the  initial singularity and, perhaps, also on the black hole singularities inside galaxies.

A completely different  approach  using local conformal symmetry in particle physics and cosmology is due to {\em Izhak Bars, Paul Steinhardt} and {\em Neil Turok}. A silent background for their interest in this question seems to have been the idea of a cyclic, respectively oscillating, model of the universe,  proposed a decade earlier by two of them  as an alternative to the  ``inflationary'' paradigm \citep{Steinhardt/Turok:2002}. 
In the latter proposal the minima of the oscillation were related to some kind of speculative physics of the string and brane type.\footnote{For a historical discussion  of oscillating models see \citep{Kragh:2009fascination} and, in an even wider perspective,  H. Kragh's contribution to this book; C. Smeenk (this volume) nicely describes the rise of the inflationary paradigm. } In  \citep{Bars/Steinhardt/Turok:2014}   the three authors  explored the possibility  that a conformal theory of gravity and the standard model fields might suffice for understanding   the  bridging process  between two cycles without necessarily much new speculative physics. They worked with  a locally scale  invariant version of the standard model, combined with gravity, similar to Nishino/Rajpoot (section \ref{subsection SM recent}), but  in the framework of purely  conformal geometry rather than Weyl geometry. 
They considered a complex valued gravitational scalar field $\phi$, called a {\em dilaton}, in addition to the Higgs field $\Phi$, both scaling with the same weight (in our notation $w=-1$). The dilaton couples only to the Higgs field by a common biquadratic potential like in (\ref{L V4}) and to the right-handed singlet neutrinos by Yukawa terms of its own \citep[p. 6]{Bars/Steinhardt/Turok:2014}.  All other masses are ``generated'' by coupling to the Higgs field like in the standard model. 

The authors investigated possible general  forms for locally scale invariant gravitational Lagrangians including a kinetic term for the  dilaton (equ. (10), loc. cit.). They were heading towards ``a fully scale-invariant approach to all physics'' (p. 5, loc. cit.) by several reasons. At first, the ``dimensionless constants in a conformally invariant theory are logarithmically divergent as opposed to the quadratic divergence of a bare Higgs mass term''  and the recent studies of \citep{Codello_ea:2013} have shown that ``the local scale invariance survives even though there is a trace anomaly'' \citep[p. 2]{Bars/Steinhardt/Turok:2014}. Moreover, so they claimed, the conformal freedom of chosing different scale gauges makes their cosmological models geodesically complete. That was a bit cavalier, but it is not the aim of this paper to evaluate such claims critically.\footnote{The authors declared geodesic
incompleteness as ``an artifact of an unsuitable frame choice:
geodesically incomplete solutions in Einstein frame may be
completed in other frames, even though the theories are
entirely equivalent away from the singularity'' \citep[p. 13]{Bars/Steinhardt/Turok:2014}.} 
More important, in our context, is to  recognize the similarity in outlook between the Weyl geometric proposals for combining gravity with  standard model fields in  a consequently scale invariant approach and the concern of our three authors. Here a perspective on a putative ``geodesic completeness'' of cosmological models came in sight, although in a rather peculiar way, not taking into account the problem of an invariant characterization of the  proper time along timelike geodesics. 

This question was discussed in  more detail  by  {\em R. Penrose} in his recent proposal for embedding the standard model of cosmology in a long cycle of iterations connected by conformal bridges between Riemannian phases of cosmic evolution  \citep{Penrose:2006CCC}. 
 He argues that  for  very high energy states in the past   timelike trajectories  lose their physical meaning anyhow  and the whole physically relevant information can be described by the structure on the lightcone. By some not yet understood processes a similar argument is imputed for  states in the asymptotic future. This  idea developed  a purely conformal perspective of how to  extend  Riemann-Einstein gravity beyond the conformally compactified past and future infinities and would probably fit well to the Bars/Steinhardt/Turok approach. Both proposals assumed that it is possible to develop a meaningful physics of the bridging process between to cycles under abstraction  from all those geometrical features which distinguish Weyl geometry  from a purely conformal structure.  

\newpage
\subsection{\small Diverse  views of Weyl geometry in cosmology\label{subsection diverse views}}
\subsubsection{Continuation of Rosen's work\label{subsection Israelit}}
M. Israelit  investigated Weyl geometric methods in cosmology in the first half of the 1990s  together with his mentor N. Rosen. After  Rosen's death in 1995 he continued publishing on his own for nearly two decades.\footnote{Isrealit died  in 2015 at the age of 87. His last paper known to me is  \citep{Israelit:2012},  a slightly changed version of \citep{Israelit:2010}. }
 In this work the question of dark matter  was studied from different perspectives,  always based on geometrical fields.  \citep{Israelit/Rosen:1992dm} explored  the  neutral massive boson interpretation of the Weylian scale connection, hinted at by Rosen already in his 1982 paper  (cf. section \ref{subsection Dirac followers}).    The authors assumed a ``chaotic Weylian microstructure'', constituted physically by a ``Weylon gas''. 
 On  large distances the scale curvature effects were negligible and  a Riemannian space structure arose  in their approach. On this basis Rosen and Israelit started to study  a hypothetical  Bose-Einstein {\em Weylon gas}   satisfying the equation of state  $\rho = 3 p$ and its consequences  for different cosmological models  in Einstein gravity \citep{Israelit/Rosen:1993dm}.

In one of their next papers they turned toward the {\em scalar field}  and tried to find out in which way it may contribute to dark matter ``pervading all of cosmic space'', i.e., on the largest scales,  not in the sense of local inhomogeneities in galaxies (like in theories of the MOND family) and in galaxy clusters. Although in their approach the Einstein gauge ($\beta = 1$) leads to ``the usual formalism of general relativity'' \citep[p. 764]{Israelit/Rosen:1995dm} our two authors  believed that different gauges with non-constant $\beta$-field might lead to new physical insight.   They  declared:
\begin{quote}
Although the gauge function is arbitrary, it leads to the presence of
dark matter which, in principle, can be observed. \cite[pp. 777]{Israelit/Rosen:1995dm}
\end{quote}
This was not particularly convincing, and it remained open how such observation ``in principle'' could be made

In some papers  of the late 1990s  and several at the beginning of the new millennium Israelit  continued this research line using {\em integrable Weyl geometric} gravity. In this context he realized that even under the assumption of an integrable Weylian scale connection the resulting modification of Einstein gravity can be {\em non-trivial},  if the potential of the scale connection $w$ is different from the scalar field, respectively its logarithm   \citep[chap. 7]{Israelit:1999Book}, \citep[equs. (17)f.]{Israelit:1999matter_creation} (compare our equs. (\ref{integrable scale connection}),  (\ref{triviality condition})). Israelit derived the dynamical equations with regard to $g_{\mu \nu}, \varphi_{\mu}$ and $\beta$, and also the Noether relations due to diffeomorphism invariance and to the scale invariance of the Lagrange density. The latter showed that on shell of the Einstein equation the dynamical equations of the scale connection $\varphi_{\mu}$ and of the scalar field $\beta$ are equivalent.\footnote{Because of scale invariance  there is, in fact, only one true scalar field degree of freedom (compare subsection \ref{subsection Weyl geometry}). \label{fn one degree of freedom}}

Israelit's aim was to explain not only dark matter but also the  accelerated expansion of standard cosmology by the gravitational scalar field which he called the ``Dirac gauge function''.\footnote{\citep{Israelit:1996,Israelit:1999matter_creation,Israelit:2002Matter_creation,Israelit:2002Quintessence}; chapters 6 and 7 in his book  \citep{Israelit:1999Book}. } 
In these papers  Israelit tried to convince his colleagues that ``cosmic matter was created by geometry'', viz. out of the energy of the gravitational scalar field \citep[p. 295]{Israelit:2002Matter_creation}. According to him, his  scalar field was able to generate dark matter and the magical substrate {\em quintessence} flourishing in the mainstream narratives on the early ``history'' of the universe. These were imaginative proposals.\footnote{I  doubt that they stand on a solid base,    although I am unable to check them in detail (E.S.). In any case, they are   too multifarious for being discussed in this survey.} 

 In one of his latest papers Israelit came back to considering a non-integrable Weylian scale connection,  now no longer as a representative of the electromagnetic potential  but again as a  field with massive bosons, ``Weylons'', of spin $-1$ and mass $> 10\, MeV$. On a microlevel, so his argument, the Weyl geometric structure appears  non-integrable, ``chaotic'', while on  larger scale there remains an effective gauge ``vector field'' with vanishing curvature. The author concluded with the remark: 
\begin{quote} `` \ldots the purpose of the present work was to show that on the basis of
the Weyl-Dirac theory one can build up a model, where conventionally matter, DM
and DE are created by geometry. This aim is achieved. \citep[sec. 8]{Israelit:2010}
\end{quote}

The ``creation'' described by Israelit did not even claim to   establish a connection between  geometry and the standard model fields. His discussion appeared as a reflex from afar on the cosmological mainstream in which elementary particle physicists  had been so successful in occupying the debate on the ``early history'' of the universe \citep{Kaiser:mass}.
Perhaps this is one of the reasons why in none  of the investigations of our last section,  nor in the ones discussed in the next subsection, we find much overlap with those  of  Rosen and Israelit.

\subsubsection{Weyl geometric extensions of gravity: trivial or provocative?\label{subsection trivial}}
Coming from Jordan-Brans-Dicke theory, {\em Israel Quiros}   got interested in  Weyl geometry  while still working at Santa Clara, Cuba, several years before his work mentioned in section  \ref{subsub Quiros ea}. He was one of those in the JBD community who took Dicke's proposal seriously, which postulated to   state natural laws in a form  that does not depend on (localized)  choices of measurement units \citep{Quiros_ea:JBD,Quiros:dual_geometries}.\footnote{Compare on this point \citep[pp. 86ff.]{Capozziello/Faraoni}.}  
At first he developed the formulation of  ``dual'' views for the interchange from Jordan to Einstein frame  \citep{Quiros:unit_transformations}. A decade later, after he had  moved to L\'eon, Mexico, he wrote a joint paper with three other Mexican authors, {\em Jose E. Madriz, Ricardo Garc\'{i}a-Salcedo, Tonatiuh Matos} in which the authors explained  how the different frames of JBD theory may be interpreted as ``complementary geometrical descriptions of a same phenomenon''  \citep{Quiros_ea:2013}. 

From there it was only a small step to  entering   Weyl geometric gravity. As  we have seen in the last section, Quiros  looked, and still looks,    for a common perspective  on gravity and a scale invariant formulation of  SM model fields, \citep{Quiros:2014a}. In a recent paper he investigated the  purely conformal approach to scale invariant Lagrangian field theories and criticized them  for lacking a well defined metrical structure with a uniquely determined affine connection. He concluded
\begin{quote}
\ldots that there will be problems with a theory
which pretends to be Weyl-invariant only because the action -- and the derived field equations --  is invariant under (2) [point-dependent  scale transformations, E.S.], but which is sustained by spacetimes whose geometrical structure does not share the gauge symmetry of the
action. \citep[p. 3]{Quiros:2014b}
\end{quote}
Quiros therefore pleaded for the use of  Weyl geometry as an  appropriate framework for his research goal.  

But alas, simplifying the Lagrangian used in \citep{Quiros:2014a} he only foresaw a kinetic term for the Higgs (or a Higgs-like) field $\Phi$, not for the gravitational scalar field $\phi$ coupling to the Hilbert term. With a gravitational Lagrangian including quartic potential
\beq L_{grav} = \frac{1}{12}\phi^2 R + \lambda \phi^4 \qquad \mbox{\citep[equ. (20)]{Quiros:2014b}} \label{L grav Quiros}
\eeq 
he found that his scalar field equation for $\phi$ reduced to the trace of the Einstein equation like in the case of conformal coupling in Riemannian geometry. 
After pondering about  the possibility of having ``an infinity of feasible – fully equivalent
– geometrical descriptions'' and the resulting paradoxical picture of an ``infinity of possible patterns of cosmological evolution'' he passed over to Einstein scalar-field gauge as ``simplest gauge one may choose''.

 For the choice of (\ref{L grav Quiros}) as the gravitational Lagrangian  this resulted in
 the  Hilbert action of Einstein gravity ``minimally coupled to the standard model of particles with no new physics beyond the standard model at low energies'' \citep[p. 9]{Quiros:2014b}.  His following discussion  reduced to the simple observation that conformal rescaling allows to scale singularities away. All this remained without new physical insights or effects; in this sense the  Weyl geometric extension of gravity  considered by Quiros up to 2014  remained physically trivial. But  it  was characterized by a conceptually clear   exposition of ideas and methods,  so   we may hope that in the  further  development of Quiros' research program  he will go beyond these  limitations.

{\em  Carlos Castro},  after the turn of the millennium  working  at  the Centre for Theoretical Studies of Physical Systems in Atlanta, USA, had become  acquainted with Weyl geometry already in the early 1990s (seee section \ref{Santamato two phases}). At that time  Santamato's  proposal for using Weyl's scale connection for  geometrizing   the quantum potential stood at the center of his interest \cite{Castro:1992}. When he became aware of the new attempts at using Weyl geometric methods  in gravity and in high energy physics,  he took up the Weylian thread again. His guiding questions were now  how Weyl's scale geometry may be  used  for  understanding dark energy and, perhaps, the Pioneer anomaly which at that time could still appear as a challenge for gravity theories \cite{Castro:FoP2007,Castro:2009}.\footnote{A few years later high precision numerical modelling showed that  thermal effects can completely account  for the  observations known as the flyby anomaly of the Pioneer spacecrafts \citep{Rievers/Laemmerzahl:Pioneer}.}
 Castro speculated with grand visions for his newly detected interest in  Weyl geometric methods, in contrast to Quiros' more sober perspective. An even  sharper  contrast comes forward with regard to  style. His  papers lack a clear conceptual  exposition; this is the reason why they are not discussed here in more detail.

Another unconventional view was put forward by the present author ({\em Erhard Scholz}, Wuppertal). His historical  studies on the work of H. Weyl lead him to the impression that already the comparatively simple modification of Riemannian geometry  by integrable Weyl geometry, 
combined with  a non-trivial scalar field extension of Einstein gravity (in the sense of our section \ref{subsection Weyl geometry} with $v+w\neq 0$),  may shed new light on certain points of present day cosmology. He was glad to find some recent activities in Weyl geometric gravity among the Munich ``group'', although it went in a different direction (section \ref{subsection Drechsler/Tann}).  

He found it most intriguing to see that  in a  Weyl geometric approach to gravity  the  cosmological   redshift    need  no longer be due  to an  expansion of the spacelike folia of Friedmann-Robertson-Walker manifolds. This was clear because, in the transition from the Riemann gauge to Einstein gauge, the warp function may be  scaled away  partially or completely 
 \citep{Scholz:2005model_building,Scholz:2005BerlinEnglish}.\footnote{A similar argument was already given by Rosen (see section \ref{subsection Dirac followers}) and more recently by  \citep[sec. 7]{Romero_ea:2012}; compare \citep[chap. 5]{Perlick:Diss}.} 
Thus   a part of the cosmological redshift $z$  may be  due to  the time component  $\varphi_o$ of the scale connection, rather than to a spatial expansion of the ``universe''. The reason for this observation is the scale invariance of  $z$, if  scale covariant geodesics of weight $w=-1$  are used.
 With regard to cosmological observers defined by a timelike geodesic flow $X$, the redshift  is given by  the quotient of energies $E_o, E_1$ of light signals (idealized ``photons'') at the event of emission $p_o$ and of observation $p_1$, 
\beq  z+1 = \frac{E_o}{E_1} = \frac{g(\gamma'(\tau_o), X(p_o))}{g(\gamma'(\tau_1), X(p_1))} \, , 
\eeq
where $\gamma(\tau)$ denotes the null-geodesic representing the trajectory of the signal. Because of the parametrization of the geodesics with weight $w=-1$  the quotient is scale invariant. 
 Although for Robertson-Walker models with non-trivial  (i.e.,  not constant) scalar field and warp function $a(t)$ in Riemann gauge, the {\em redshift} seems to result from an  expanding warp function,   in Einstein scalar-field gauge it  {\em may}  at least   partially  be due  to the scale connection, i.e., to a {\em field effect} of the additional component of the gravitational structure. 

This effect is particular striking in certain models which  appear expanding in Riemann gauge but 
have a static metric in Einstein gauge  \citep{Scholz:2005BerlinEnglish,Scholz:FoP}. Here the cosmological redshift in Einstein gauge turned out to be completely due to the time component of the Weylian scale connection, $H = \varphi_o$, with $H$  the Hubble parameter. Although Scholz initially overestimated the physical import of this example,  it can probably be  a fruitful  epistemic  provocation.\footnote{It was the topic of the author's talk at the Mainz conference.}  
As a toy model it may continue to serve as an incentive for critically  rethinking the  foundations of our   standard picture of  the universe.  

\newpage

\subsection{\small Attempts at dark matter, MOND-like\label{subsection dark matter}}
The success of {\em modified Newtonian dynamics}, MOND, since the 1980s  for explaining the rotation curves of galaxies and  the Tully-Fisher relation between the luminosity of spiral galaxies and their angular velocity  led to diverse attempts for general relativistic generalizations \citep{Sanders:DarkMatter}. Some of them introduced  non-geometrical structures, like an additional vector field in the so-called tensor-vector-scalar field theory, TeVeS; but the earliest attempt at a relativistic MOND-like theory was formulated by {\em Mordechai Milgrom} and {\em Jacob Bekenstein} in the framework of JBD gravity \citep{Bekenstein/Milgrom:1984}. 

This approach worked with a non-quadratic kinetic Lagrangian for the scalar field with MOND-typical transition function; therefore its name  ``relativistic a-quadratic'' rAQUAL theory. Bekenstein and Milgrom were torn between the Jordan frame  and Einstein frame.  In a later review  paper  by one of the authors the Jordan frame was declared to be  the ``physical metric'', while Einstein frame  was considered  as the ``primitive metric''  \citep[p. 6]{Bekenstein:2004}.\footnote{Bekenstein considered it still in 2004 as ``evident'' that  measurements with `` clocks and rods'' are expressed by the Jordan metric. Moreover the latter's  Levi-Civita connection governs   the free fall of test particles. But the dynamics does not satisfy the ``usual Einstein equation'' in Jordan frame (because of explicit terms in  the scalar field). The Einstein frame represented for him the  ``primitive metric'' because here the gravitational action reduces to  the classical form of the Hilbert term, and the dynamics is given by the Einstein equation \citep[p. 5f.]{Bekenstein:2004}. In their common paper, Milgrom and Bekenstein  used the terminology of ``dual descriptions'' working in ``gravitational units'' (Einstein frame) respectively ``atomic units'' (Jordan frame), which sounded a bit like Dirac's distinction \citep[p. 14]{Bekenstein/Milgrom:1984}.}
In this framework the MOND-like free fall of particles in an extremely weak gravitational field could be derived.  This approach was not free of  shortcomings, as the  authors themselves remarked: Lensing effects seemed unexplainable by the approach, because the conformal change between the two ``dual'' frames seemed not to affect light-like geodesics.\footnote{In his later review paper Bekenstein qualified this point by stating that only as long as the scalar field ``\ldots  contributes comparatively little to the
energy-momentum tensor, it cannot affect light deflection,
which will thus be due to the visible matter alone'' \citep[p. 6]{Bekenstein:2004}. One can read this observation the other way round: If the scalar field carries  a considerable contribution to the energy-momentum it influences  light deflection. }
 Moreover,  the scalar field  allows   perturbations which propagate with superluminal velocity. The authors relativized this problem, however, by adding that such perturbations  probably ``cannot induce acausal effects in the behavior of particles and electromagnetic fields'', because they only relate to the conformal factor of the metric \citep[p. 14]{Bekenstein/Milgrom:1984}.
An additional critical point, not only for rAQUAL but for all theories of the original MOND family, was their inability to explain the anomalous dynamics  in  galaxy clusters, without assuming some additional unseen matter. 

Because of the close relation between integrable Weyl geometric gravity  and JBD theory, Bekenstein's and Milgrom's rAQUAL may be an interesting  challenge for testing  what happens if it is transformed into a  scale invariant framework.  In two recent papers the present author   investigated this problem   \citep{Scholz:2016MONDlike,Scholz:2016Clusters}. The first paper contained a Weyl geometrical  reformulation of  rAQUAL, at least for the so-called ``deep MOND'' regime and the upper transitional regime. Different from Bekenstein/Milgrom's view,   in the Weyl geometric approach  observable quantities   are most directly expressed in Einstein gauge. Spacelike components of the Weylian scale connection, $\varphi_j, \, j=1,2,3$  express additional accelerations in comparison with those induced by  the Riemannian part of the metric (corresponding to Newtonian ones in the weak gravity limit).  In the extremely weak gravity regime two different components of additional accelerations, $ a_ {\phi},\, a_{\varphi}$,  can   be distinguished. The first one is part of the Riemannian acceleration, in  Einstein gauge, and due to   the energy density   of the scalar field $\phi$;  the second one results from the the Weylian scale connection $\varphi$ (in Einstein gauge). In extremely weak static gravitational constellations (i.e., order of magnitude of Newton acceleration $a_N$ close to the MOND acceleration $a_o \approx 1.2\, 10^{-10}m s^{-2}$), MOND-like phenomenology is reproduced similar to rAQUAL. But here half of the additional acceleration is due to the {\em scalar field's  energy}. It thus influences the light trajectories.\footnote{Cf. last footnote.} Whether this suffices for explaining the observed lensing effects remains to be seen. 

Moreover, contributions on different length scales to local inhomogeneities of  the scalar field's energy density can add up to a common effect. This seems to  have striking consequences  for  the dynamics of galaxy clusters. In a  heuristic investigation   of data  from 17(+2) clusters\footnote{The data for 2 clusters are outliers, already from the phenomenological point of view.} our author found an encouraging agreement of accelerations  predicted  by the Weyl geometric scalar tensor theory with the corresponding empirical values    \citep{Scholz:2016Clusters}. This was done on the basis  of the observed baryonic masses alone, without  assuming  additional unseen, ``dark'', matter.\footnote{The famous Coma cluster  which led  Zwicky introduce the hypothesis of dark matter is among the galaxy clusters for which the Weyl geometric model is consistent with most recent empirical data on mass distributions and accelerations.}

Calculations with ordinary  MOND, or even its relativistic generalization  TeVeS,  reduces the need of assuming additional hypothetical dark matter, but cannot do  completely without it. R. Sanders argued that sterile neutrinos could do the job. In this context it is   interesting to see how,    in the  Weyl geometric  framework,  an {\em outlandish kinetic term} of  the scalar field  seems to {\em  suffice for explaining} the otherwise anomalous {\em dynamics of galaxy clusters}. 

Of course, there remained problems: For safeguarding  the dynamics on the solar system level the author invented an  (ad-hoc) hypothesis postulating that  scalar field inhomogeneities are suppressed in regions where the value of at least one  sectional curvature of the Riemannian component exceeds a certain threshold \citep[p. 6]{Scholz:2016Clusters}. That saves the dynamics, but the origin of such a hypothetical suppression remained  unclear.  Moreover, the cosmological consequences of this approach are far from clear.\footnote{Unpublished calculations indicate scenarios of a cosmic evolution in agreement with many features of standard cosmology: initial singularity, large parts of the cosmological redshift due to the  expansi on of spatial folia in Einstein gauge, accelerated ``late time'' expansion etc.} In spite of such shortcomings this model may be of some value for exploring the possibilities of Weyl geometric scalar fields in the realm of dark matter phenomena.

\subsection{\small  The Brazilian approach\label{subsection Brazilian approach}}
A   challenge to the standard big bang picture, drawing upon Weyl geometric methods, came from Brazil.   Interest in Weyl geometrical approaches to  cosmology have been present  in the Brazilian theoretical physics community since the 1990s.    The central person for this development,  {\em  M\'ario Novello}, acquired his doctorate  in 1972 at Geneva under the supervision of J.M. Jauch.  Already as a young PhD student he published a paper on Dirac spinors expressed in quaternionic calculus in a Weyl space \citep{Novello:1969}.\footnote{In this paper Novello still thought in terms of Weyl's first interpretation of the scale connection, the ``em dogma'' in the terminology above. }
Back in Brazil and working at the {\em Centro Brasileiro de Pesquisas Fisicas}, Rio de Janeiro, he   cooperated with many international guests. In  2003 he   became the founding director of  the {\em Instituto de Cosmologia Relatividade e Astrofisica (ICRA)}. Due to his influence  Weyl geometric ideas were introduced in  the Brazilian community of theoretical physicists  in the course of  the 1990s. They flourished and  turned into an  research tradition of its own, the {\em Brazilian approach} to Weyl geometric gravity as I want to call it.  

\subsubsection{A Palatini-type path to integrable Weyl geometry\label{section Palatini}}
In the early 1980s Novello and a  co-author from Cologne, {\em  H. Heintzmann},   reflected   possible consequences  for cosmology if one allows to model it in a slightly more general framework than   Riemannian  geometry \citep{Novello/Heintzmann:1983}.  Like  other authors before them, they used a {\em metric-affine } approach to  gravity, presupposing   a  metric $g$ and an independent  affine connection $\Gamma$. This allows to define curvature tensors like in Riemannian geometry, including the scalar curvature $R$.  Starting from a  gravitational Lagrangian which included a term of the form 
\beq \mathfrak{L}_R=-e^{\omega}R\sqrt{|g|} \label{mod Hilbert term}
\eeq
 with point-dependent function $\omega(x)$,\footnote{In this paper $\omega(x)$ was not yet introduced as a scalar field of its own, but via  the square of the electromagnetic potential $A_{\mu}$, i.e., $\omega=\log{A_{\mu}A^{\mu}}$.}
 metric and affine connection were varied  independently according to the so-called {\em Palatini approach}. They found that the  variation with regard to the connection implies 
\beq \nabla_{\lambda} \,g_{\mu \nu} = - \partial_{\lambda}\omega \, g_{\mu \nu} \;  . \label{Palatini relation}
\eeq
Our authors immediately  realized that this relation can be identified with  the Weyl geometrical compatibility condition of our equation (\ref{metric compatibility})  for the  integrable scale connection 
\beq \varphi = \frac{1}{2}d\omega \label{varphi Novello}
\eeq 
  (in our notation). This  approach was not without limitations: it identified the scale gauges in which the coefficient $e^{\omega}$ of the Hilbert term in (\ref{mod Hilbert term})  becomes constant  and the one in which the scale connection (\ref{varphi Novello}) vanishes. 	In our terminology above   {\em no difference} between {\em  Riemann gauge} and scalar field -- {\em Einstein gauge} could be conceived! This structural identification of the two gauges made the   	Palatini  approach to Weyl geometric gravity a   trivial extension of Einstein gravity, if the full scale invariance of the Lagrange density is observed. But such a comparison was not in the mind of our authors. 	
	Referring  to   Canuto et al. and in the wake of Dirac  (cf. section \ref{subsection Dirac-Utiyama}),  they pondered about the possibility that atomic clocks and gravitational clocks at  different places  might be  related by a  variable factor $\omega(x)$. If $\omega(x)$ is asymptotically constant, different ``Riemannian domains'' would arise, possibly  connected by ``Weyl integrable regions of space''. Moreover, the ``age'' of the universe might become ``arbitrarily large'' \citep{Novello/Heintzmann:1983}.

In the following years Novello developed broad activities in gravitation theory, elementary particle physics, and cosmology;  in particular  he was interested in understanding how the  initial singularity of standard Riemann-Einstein cosmology can be avoided. In a joint paper with {\em Edgar Elbaz}, a colleague from France,  {\em Jose M. Salim} and {\em L.A.R.  Oliveira} from his group,  he and his co-authors proposed an imaginative model for what they called  the ``creation of the universe'' (a clause from the title of the paper) \citep{Novello/Oliveira_ea:1992}. Using some Weyl geometric features and  a scalar field $\omega$, the authors  were able to display a ``cosmic'' development from a flat vacuum state (described by Minkowski space) via a contracting phase, ``bouncing'' at a minimum of a scale function, to an expanding  ``inflationary''  phase. Without  going into details of this study we want to see how and why this paper became a classical point of reference for the Brazilian tradition in Weyl geometric methods.\footnote{In many papers of the Brazilian tradition \citep{Novello/Oliveira_ea:1992} is   quoted as a starting point \citep{Salim/Sautu:1996,Oliveira/Salim/Sautu:1997,Fonseca/Romero_ea:GR_Weyl_frames,Romero_ea:2012}, to cite just a few. Sometimes it is even called the ``first approach to scalar-tensor theory in WIST'' [Weyl integrable space-time] \cite{Pucheu_ea:2016}. }

 Weyl geometry was introduced,  like in \citep{Novello/Heintzmann:1983}, i.e., by the  Palatini method of variation (\ref{Palatini relation}). This led to an integrable  Weyl geometry   characterized by a scalar function $\omega (x)$, the potential of the scale connection  $\varphi= \frac{1}{2} d\omega$.   For Novello, Elbaz et al.  the  above mentioned identification of Riemann gauge and Einstein  did not appear detrimental, because their goal was not a modification of Einstein gravity. They  rather set out  modelling semi-classical quantum ``perturbations of the system of measurement units'' described by  $\delta \omega_{\lambda}$, such that
\beq \delta(\nabla_{\lambda}g_{\mu \nu})= (\delta \omega_{\lambda}) g_{\mu \nu} \; .
\eeq 
Perturbations of such a kind are inconsistent with Riemannian geometry but consistent with Weyl geometry, as the authors noted with references to \citep{EPS,Audretsch:1983,Perlick:Observerfields}.
Like many physicists in the last third of the 20th century, they thought in terms of a time-evolution of the cosmos, here even in the sense of a {\em temporal evolution of its geometrical structure}.  

They hoped to find ``a definite conceptual context \ldots for the description of such  structural transitions'' during the cosmic evolution in the Weyl geometric approach \citep[p. 650]{Novello/Oliveira_ea:1992}. For this goal they  considered a process 
governed by the  Lagrangian
\beq \mathfrak{L}_{vac} = (R+ \xi \nabla_{\nu}\partial^{\nu} \omega)\sqrt{|g|} \qquad \mbox{\citep[equ. (4.2)]{Novello/Oliveira_ea:1992},} \,  \label{Lagrangian Novello ea 1992}
\eeq 
where $\nabla_{\nu}$ is (our) notation for the Weyl geometric derivative and $R$ denotes the Weyl geometric scalar curvature.\footnote{In a side remark the authors reminded that $-2\xi \omega_{\lambda}\omega^{\lambda}$ is  a  variationally equivalent  kinetic term, because the difference to the kinetic term in (\ref{Lagrangian Novello ea 1992}) is a  total divergence \citep[p. 654]{Novello/Oliveira_ea:1992}. \label{fn equivalent kinetic term Novello ea}} 
From the point of view of Weyl geometry this  was a {\em hybrid approach}; the Lagrange density  was not scale invariant, although Weyl geometric concepts and expressions were used.  The authors considered this an advantage, because  a difference between ``gravitational'' units (expressed by a point dependent gravitational ``constant'') and  atomic units, originally assumed by Dirac, Canuto et al. (see section \ref{subsection Dirac}),  appeared   unacceptable to them. They rather assumed a broken (active) scale symmetry \citep[p. 653]{Novello/Oliveira_ea:1992}; a mere transformation of units in Dicke's sense, i.e. a passive conception of scale covariance was not their case. Understandably, they decided for Einstein gauge as the expression for the broken symmetry state. They thus 
understood (\ref{Lagrangian Novello ea 1992}) as an ``effective canonical action'' of  a broken underlying scale symmetric dynamics with some surviving residual Weylian terms. Guided by such kind of physical intuition the authors avoided a reduction of their approach to Einstein gravity, which would have become necessary, had they assumed full scale invariance. 

The resulting dynamical equation could  be expressed, without loss of content, in Riemannian terms,
\beq _g\hspace{-0.15em}R_{\mu \nu} - \frac{_g\hspace{-0.15em}R}{2}g_{\mu \nu}= \lambda^2 \omega_{\mu}\omega_{\nu} - \frac{\lambda^2}{2} \omega_{\alpha}\omega^{\alpha} g_{\mu \nu} \, , \qquad  \label{Einstein equ Novello ea}
\eeq 
with $ \lambda^2= \frac{1}{2}(4\xi - 3)$. 
Then it was ``equivalent to an Einstein equation in which the WIST\footnote{WIST was (and is) the abbreviation, preferred by  the Brazilian authors, for ``Weyl integrable scalar tensor theory''.}
  field $\omega$ provides the source of the Riemannian curvature'' \citep[p. 655]{Novello/Oliveira_ea:1992}.
As $\omega$ was  the integral of the scale connection, it  had a ``purely geometrical origin'' and appeared acceptable to them, although it had strange physical properties: negative energy density, positive pressure of the same value (``stiff'' matter).

The scalar field equation derived from (\ref{Lagrangian Novello ea 1992}) and  the Einstein equation (\ref{Einstein equ Novello ea}) evaluated for a homogeneous, isotropic spacetime led to a model without initial singularity. In the far past it looks like a contracting  Minkowski space with a  non-trivial and in this sense ``excited''  scalar field $\omega$. After a first phase of an accelerated contraction, their warp function $a(t)$ reaches a minimum value $a_o$, after which it turns into an expansionary phase. The authors interpreted the first, contracting phase as a vacuum with a geometrical  scalar field excitation. Near the minimum they sketched quantum processes of photon and baryon genesis ``driven'' by the scalar field. Then an expansionary phase follows, ending in a state which, so they argued, could connect to the radiation dominated phase of the standard model of cosmology. All in all the calculations were embedded in an imaginative narrative which claimed to solve several pressing problems inherent in the standard picture of the ``hot big bang'' (no initial singularity, causal horizon and flatness problems, matter anti-matter asymmetry).

\subsubsection{Cosmological models with fluid matter}
Several  follow up papers appeared, among them \citep{Salim/Sautu:1996,Oliveira/Salim/Sautu:1997}. In the  first one Salim and {\em S.L. Saut\'u}  added different types of ``external fields'' representing matter and its interaction  to the  vacuum Lagrangian  (\ref{Lagrangian Novello ea 1992}). At first they  dealt with an electromagnetic field and an external scalar field. Their   matter terms of the  Lagrangian had a  scale invariant form
 \citep[equ. (12)]{Salim/Sautu:1996}.\footnote{This adds flavour to the hybrid approach mentioned above.} 
 More important  for cosmology, in the next step they adapted  the Lagrangian of a perfect fluid following the trajectories of a timelike vector field to their framework.\footnote{The fluid Lagrangian was taken over  from \citep{Ray:1972}.}
Here the hybrid form of the approach with the specified scale gauge was of great advantage, because it facilitated the adaptation of the fluid Lagrangian. The authors derived the dynamical equations and constraints in their framework, first in terms of the Weyl geometric derivative and  curvature expressions, with  particular taking care for the interaction with the geometrical scalar field $\omega$ \citep[equs. (34)-- (40)]{Salim/Sautu:1996}. 
After that they rewrote the Einstein equation and the scalar field equation in Riemannian terms  and derived the corresponding generalized Raychaudhuri equation for the  homogeneous isotropic case (equ. (47)). Rewriting the coupling constant $\xi$ of (\ref{Lagrangian Novello ea 1992}) by $\lambda=\frac{1}{2}(4\xi -3)$  they  concluded that 
``\ldots depending on the sign of $\lambda$, the cosmological solution under
consideration can be non-singular and inflationary'' \citep[p. 359]{Salim/Sautu:1996}.

That was a considerable step forward and  generalized the effect observed for the case of  the special vacuum solution in \citep{Novello/Oliveira_ea:1992}.  The authors rightly concluded:
\begin{quote}
We have shown that the Weyl integrable
geometry can be used in a natural way to geometrize a long-range scalar field. Using
a general principle to prescribe the interaction of the geometric scalar field with other
physical systems, we can describe in WIST all the classical situations studied by EGR [Einstein gravity, E.S.]. \citep[p. 359]{Salim/Sautu:1996}
\end{quote}

In their next paper the two authors, now supported by {\em Henrique P.  de Oliveira}, studied ``non-singular inflationary cosmologies in Weyl integrable spacetime'' \citep{Oliveira/Salim/Sautu:1997}. To the  gravitational Lagrangian (\ref{Lagrangian Novello ea 1992}) of  Novello et al.  they added  a self-interaction potential of the scalar field $V(\omega)$ and the fluid Lagrangian of \citep{Salim/Sautu:1996}. Referring to the same parameter $\lambda$ as above, they came to the conclusion that for $\lambda >0$ the Friedmann-like solutions had strong similarities to those of Einstein gravity, while for $\lambda<0$ interesting ``novelties appear in WIST'' (p. 2835).

\begin{figure}[t]
	\centering 
		\hspace*{-4em}
		\includegraphics[scale=0.8,trim=0mm 70mm 0mm 110mm,clip]{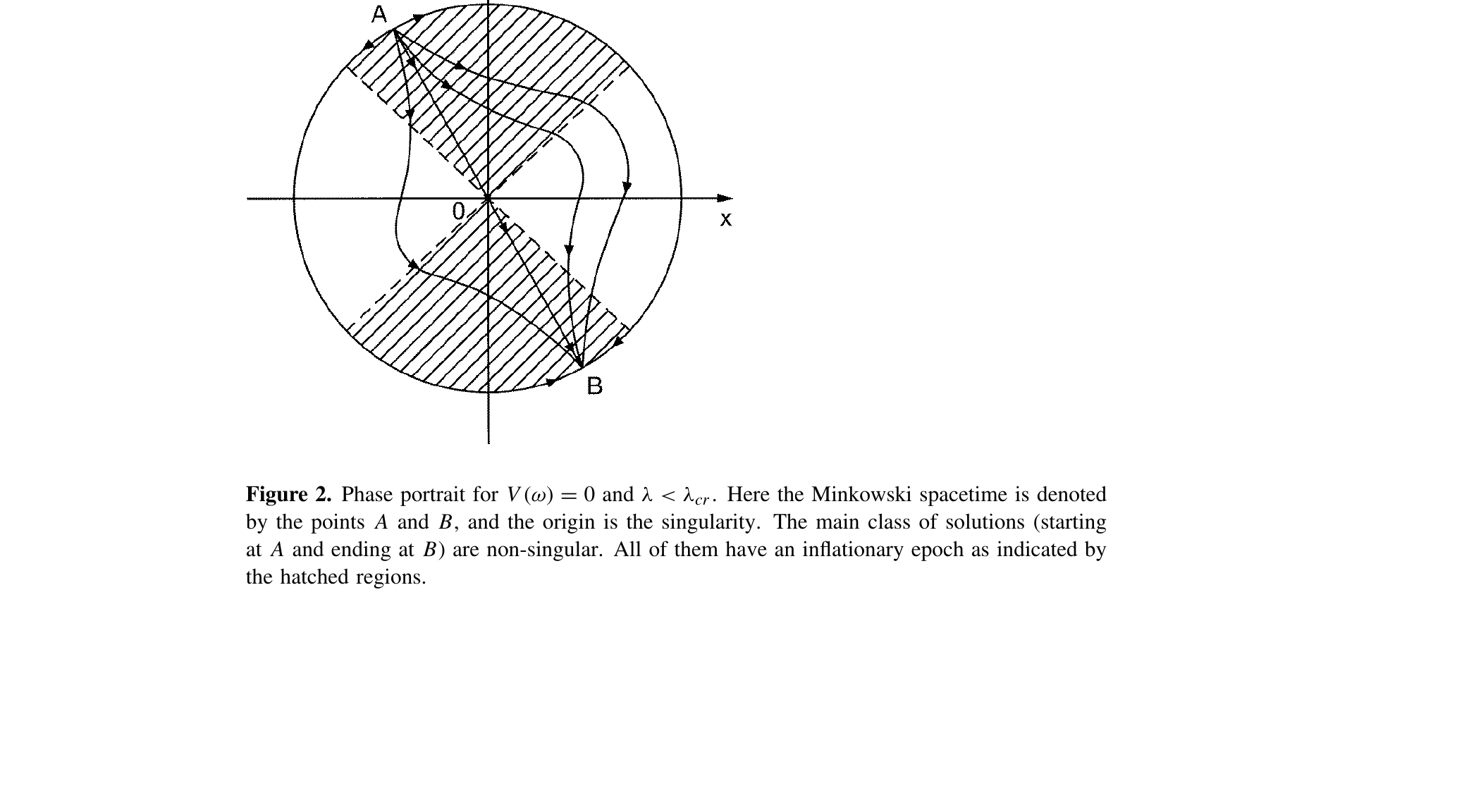}
				\vspace*{-30mm}
	\caption{Phase portrait (de Oliveira e.a. 1997, fig 2)}
	\label{fig:phase_portrait_1}
\end{figure}

The three authors studied the qualitative behaviour of the modified Friedmann and scalar field equations of their model  in the parameter plane $(x,y)$ with
\beq x= \frac{\dot{a}}{a}\, , \qquad \qquad y = \dot{\omega} \, .
\eeq    
For  vanishing potential, $V(\omega)=0$, and for an exponential potential $V(\omega)= V_o e^{\beta \omega}$ ($V_o, \beta$ constants) they found that the solutions of the Friedmann equation are generically singularity free, while the solutions with initial or final singularity are instable (see our fig. \ref{fig:phase_portrait_1}). This was a  striking result. The authors commented:
\begin{quote}
Depending
on the parameter $\lambda$, we obtained non-singular models as a general feature. (\ldots)
The non-singular behaviour is explained by the violation of the strong energy condition
provided by the geometric scalar field. \citep[p. 2842f.]{Oliveira/Salim/Sautu:1997}
\end{quote}
The explanation indicates that   the ``strong energy condition'' was understood in the geometrical sense ($R_{\mu\nu} V^{\mu}V^{\nu}>0$ for any timelike vector field $V^{\mu}$).\footnote{The ``physical sense'' of the strong energy condition is $T_{\mu \nu} - \frac{1}{2} tr\, T\, g_{\mu \nu}V^{\mu} V^{\nu}$. Geometrical and physical conditions are equivalent {\em in Einstein gravity}; see, e.g.,  \citep[p. 49]{Curiel:Energy_cond}.}
Later investigations of the Brazilian school would show that  the geometrical energy condition may be violated, while the  physical energy  condition may still be satisfied (see below).

\subsubsection{A tension between  the Palatini approach and scale invariance}
Many more papers on Weyl geometric gravity were published by the Brazilian group.   Some young researchers joined the network and started to publish with colleagues from the older generation in different constellations; among them and  not yet mentioned before   {\em  Tony S. Almeida,  F.A.P. Alves, Aadriano B. Barreto, J.B. Fonseca-Neto, F.P. Poulis} and {\em  Carlos Romero } (in alphabetical order). They dealt with the  relationship of the Palatini variant of Weyl geometric gravity to  Einstein gravity    and to JBD theory, and continued to study singularity behaviour of cosmological models in their slightly extended framework. 

Several of these papers   \citep{Fonseca/Romero_ea:GR_Weyl_frames,Romero_ea:2012}  used  a   Lagrangian of the form:
\beq \mathfrak{L}= e^{(1-\frac{n}{2})\omega}(R + 2 \Lambda e^{-\omega} + \kappa  e^{-\omega} L_m ) \sqrt{|g|} \; , \label{Lagrangian Romero ea 2011}
\eeq 
where  the scalar curvature is to be understood  in the metric-affine sense ($R=R(g,\Gamma)$).  After a Palatini type variation like in the transition from (\ref{mod Hilbert term}) to (\ref{Palatini relation}) $R$ turns into the scalar curvature of a Weylian metric given by the pair $(g, \varphi =\frac{1}{2}d\omega) $.\footnote{\citep[eq. (7)]{Fonseca/Romero_ea:GR_Weyl_frames},\citep[equ. (12)]{Romero_ea:2012} } 
The authors emphasized the importance of what  appeared  to them a ``new kind of invariance, namely with respect to Weyl transformations''  \citep[8]{Romero_ea:2012}  without, however, keeping coherently to scale invariance as a guiding principle of their investigation. 
Rewritten in Riemannian terms this Lagrangian acquires the form of a Brans-Dicke Lagrangian with a conformally coupled scalar field.
Not very convincingly, this was  presented as a ``geometrization'' of  JBD scalar fields in general and, in addition, as an argument for a compatibilty of Einstein gravity and scale invariance in the sense of   Weyl geometry.\footnote{``An important conclusion (\ldots) is that general relativity can perfectly `survive' in a non-Riemannian environment'' \citep{Fonseca/Romero_ea:GR_Weyl_frames} etc.}

Such a  generalization of Einstein gravity is  clearly too weak to lead to interesting new features (see section \ref{subsection trivial});  it even does not allow to recuperate the Lagrangian (\ref{Lagrangian Novello ea 1992}), so important for the Brasilian tradition.   Other papers thus start from  a metric-affine generalization of  JBD-type  Lagrangian with general coupling coefficient  $\omega$,   written in the form 
\beq \mathfrak{L}_{JBD} = e^{-\phi} (R + \omega\, \partial_{\mu}\phi  \partial^{\mu}\phi) \sqrt{|g|} \; ,
\eeq 
 with $R=R(g,\Gamma)$ as above. Again the Palatini variation implied the relation (\ref{Palatini relation}) and, in this way, a motivation for specifying the metric and connection in the sense of  integrable Weyl geometry \citep[equ. (2)]{Almeida/Pucheu/Romero:2014}.

 But then the kinetic term,  taken over without change from usual, Riemannian,   Brans-Dicke theory breaks the scale invariance for general $\omega$. Accordingly the authors used a restricted Weylian scale transformation only. For the transition to the  ``Einstein frame'' they   transformed the quantities $ e^{-\phi}$ and $R$ only, while leaving the core expression of the  kinetic term unaffected (a ``field substitution'' rather than a gauge transformation), with the result
\beq \mathfrak{L}= (\tilde{R}  + \omega\, \partial_{\mu}\phi  \partial^{\mu}\phi) \sqrt{|g|}  ) \label{recent Brazilian Lagrangian}
\eeq 
plus a matter action $\mathfrak{L}_m$. They were thus led back to the archetypical form (\ref{Lagrangian Novello ea 1992}) of the gravitational Lagrangian in the Brazilian tradition.\footnote{\citep[equ. (3.24)]{Almeida/Pucheu/Romero:2014}, \citep[equ. (16)]{Almeida/Pucheu:2014}; compare (\ref{Lagrangian Novello ea 1992}) in the light of fn. \ref{fn equivalent kinetic term Novello ea}.} 

The  tension between the Weyl geometric frame and the general methodology did not pass unnoticed by the authors. But it  seems that they  were prevented from resolving it, because their adherence to the Palatini method of variation and the difficulty  to express  the  kinetic term of the scalar field in  a scale invariant form without using Weyl geometric scale covariant derivatives (\ref{L phi}). In the conclusion of one of the papers they wrote  `that neither the action
nor the field equations of the proposed theory are invariant
under Weyl transformations'',  admitting that ``it would perhaps be desirable, at least
from the aesthetic viewpoint, that the whole theory should
exhibit Weyl invariance`'' \citep[p. 8]{Almeida/Pucheu:2014}. In another paper three of them even gave a twisted explanation why this seemingly must be so \citep[p. 39]{Almeida/Pucheu/Romero:2014}.
 This is surprising, because  two years earlier  two of them, in this case added by Fonseca-Neto, had already  noticed  that the Lagrangian (\ref{recent Brazilian Lagrangian}) can be scale-transformed into the form of the general JBD Lagrangian (\ref{Lagrangian JBD Fujii/Maeda}). So they were quite close to bringing the Brazilian approach into a coherently scale  invariant form;\footnote{One only needed to put the  JBD Lagrangians in a Weyl geometric framework. Alternatively, if  one wants to  start from the Brazilian point of view, one may read the constant coefficient of the Hilbert term in (\ref{recent Brazilian Lagrangian}) as  the  value of a scale covariant scalar field $\chi$ in (Einstein-) scalar field gauge, $\chi_o \doteq 1$, and $\partial_{\mu}\phi  \partial^{\mu}\phi$ as the scalar field gauged expression of the scale covariant  kinetic term $D_{\mu}\chi D^{\mu}\chi$ with scale covariant derivative (\ref{scale covariant derivative}). }
 but by some reason or other the members of the Brazilian group did not dare, or felt unable, to  transcribe their approach in a scale invariant mode.

Gravitational Lagrangians of the form (\ref{recent Brazilian Lagrangian})  also played a  role in  recent  qualitative studies of ``isotropic cosmologies in Weyl geometry'' by {\em John Miritzis} from Athens  and  by authors from  the Brazilian network itself  \citep{Miritzis:2004,Pucheu_ea:2016}. We find    qualitative studies of cosmological models  with or without (initial or final) singularities. The questions and results extended those of  the 1990s.    New glances at global singularities  in the slightly extended framework described above were added   \citep{Lobo/Romero_ea:singularities}.  After a  detailed investigation of the  Raychaudhuri equation, the authors of the last mentioned paper  showed that  the geometrical version of the strong energy condition can  be violated in the  Brazilian approach, while the physical one may be maintained due to contributions of the  energy-momentum of the  scalar field.  This was a sharp observation and may be of wider import.


\section{\small Discussion\label{section discussion}}
\subsection{\small A rich  history aside the mainstream}
Our  survey over the reappearance of Weyl geometry  has encountered {\em  four different entrance channels} through which central concepts of Weyl's scaling invariant but still fully metrical, geometry of 1918 were reintroduced into  late 20th century physics. They differed in motivation and systematics;  three of them  were opened even twice by essentially independent research initiatives and  slightly differing  systematic ideas (these are characterized by an ``and'' in the following list):
 \begin{itemize}
\item[1.] Axiomatic foundations of gravity: Ehlers/Pirani/Schild (section \ref{subsection EPS})
\item[2.] Scale co/invariant scalar tensor theory of gravity: Om\-ote/Utiy\-ama and Dirac (section \ref{subsection Dirac-Utiyama})
\item[3.] Cartan geometric approach: Bregman and Charap/Tait (section \ref{subsection Cartan-Weyl})
\item[4.] de Broglie-Bohm-Madelung (dBMB) approach: Santamato and Sho\-jai/Sho\-jai/Golsh\-ani (section \ref{section geometrical quantum mechanics})
\end{itemize}

The {\em first three} were initiated in the  short time interval  1971--1974, and two of the openings were taken twice. Dirac's motivations for his multiple gauge approach to gravity was quite idiosyncratic and played a minor role for its reception. The broader scenario indicates an intellectual environment which  let it appear natural to come back to Weyl's proposal of   generalizing Riemannian geometry in a new field theoretic context.  Bjorken scaling had attracted attention in the late 1960s, but was known to be only approximatively valid already at that time. So it could not have been a major driving force.\footnote{Only for the authors of \citep{Hehl_ea:1988Kiel} this appeared to be different, see section \ref{subsection Cartan-Weyl}. }
 On the other hand, the field structures of elementary particle physics  were just acquiring the form and status of a new, gauge theoretic   standard model, due to to the renormalizability results of 't Hooft and Veltman (1972) and the experimental detection of quark binding states, called  ``J-$\Psi$'' (1974). Their Lagrangians  were basically (globally) scale invariant in Minkowski space, with only the mass term of the  hypothetical Higgs field as a scale breaking term. This context may have strongly motivated researches which explored possible connections to gravity in an enlarged scale covariant framework. In such a wider perspective a new look at   Weyl geometric generalizations of Einstein gravity must have appeared a promising perspective. In this respect it was important that the  Jordan-Brans-Dicke research program of scalar tensor theories  had shown already a decade earlier  how one could model  gravity without taking  recourse to a quadratic curvature term. It was natural to do so also in a renewed Weyl geometric setting. This brought it much closer to Einstein gravity than the quadratic gravity theories studied since the time of Weyl. That was important because during this time Einstein gravity lived through a vivid phase of new empirical confirmations.\footnote{Cf. C. Will's contribution to this volume.}

The {\em fourth opening} was anchored in the completely different intellectual context of the  de Broglie-Bohm program for reconsidering the foundations of quantum mechanics. Its two research lines  started a decade later  than the first three (Santamato), 
or even two decades later (Shojai/Shojai/Golshani) which in the following will be referred to as the ``Tehran approach''. The two lines differed among each other more strongly than the respective double starts in items 2 and 3.  In the 1990s the Bohmian approaches entered a latency phase (Santamato) or were  just heading towards a new beginning,  still developed in a JBD framework (Iranian approach).  The authors of the latter started  to use  Weyl geometric concepts explicitly only after the turn to the new millennium. 

All in all, the time  until roughly 2000 was  a {\em first phase of exploration} for all the approaches. For several years the immediate continuators of the Dirac line explored astrophysical consequences of Dirac's distinction of an ``atomic gauge'' and ``Einstein gauge'' (Bouvier, Maeder, Canuto et al.) or refined and extended the theory (Rosen). At  first they sticked  to  Dirac's interpretation of the scale connection  as an electromagnetic potential, the  {\em em} dogma.   In the 1980s such a literal allegiance  of Dirac's ideas faded out. Those who continued to appeal to Dirac's approach, like Rosen and Israelit, enriched the perspective by considering the scale connection as a representative of a Proca-like massive gauge field, or saw it in a different context anyhow like Smolin. This  boiled down to a merging of the modified Dirac line  with  the research  following Omote/Utiyama's  initiative (Hayashi, Kugo et al.), which intensified with the attempts to bring Weyl's scale geometry in contact with the field content of the rising standard model (see below). 

On the other hand, the later researches of Rosen and those of Israelit explored a vast terrain of theoretical possibilities, many of them  quite speculative, of how  the energy-momentum of a Weyl geometric gravitational scalar field or  a hypothetical ``Weylon gas'' might contribute to dark matter phenomena and/or to the accelerated expansion diagnosed in the  usual Riemannian approach to gravity. But these studies remained on a relatively general level and remained without closer links to astrophysical or astronomical observations (section \ref{subsection Israelit}). 

 The {\em Cartan-Weyl geometric approach} was soon relegated to a very special case in the broader Cartan geometric metric-affine theories (Hehl et al.) or was studied in relation to Kaluza-Klein theory.\footnote{E.g. \citep{Drechsler/Hartley}.} 
 As the latter are not included in this survey, it disappear more or less from the range of our panorama. The foundational studies of Ehlers/Pirani/Schild, on the other hand, found a broad and continued reception and development in the philosophy of physics and remained a point of orientation for  foundational studies of gravity.

For some authors  (Englert, Smolin, Cheng, later Drechsler, Tann) the {\em rise of the standard model}  suggested to  connect Weyl geometric gravity, or at least scale covariant gravity in the case of Englert et al., with standard model fields, in particular the Higgs field. 
 Cheng's seminal paper of 1988 was  the first relatively detailed account of the electroweak sector of the SM assimilated to a Weyl geometric context, although only on  the pseudo-classical level of the theory (section \ref{subsection SM 1970s}).\footnote{Studies of QFT on Weylian manifolds, comparable to the corresponding researches for Lorentzian manifolds,  discussed in R. Wald's contribution to this volume, are still a desideratum.}
 Here we also find  explicit references to the papers of Dirac and the Utiyama research tradition, indicating the merging of these lines mentioned above.  Nearly a decade later Drechsler and Tann  found much of the electroweak structure in their own development of Weyl geometry,\footnote{The two authors neither referred to the Dirac tradition in Weyl geometric gravity  nor to Utiyama's; their Weyl geometric starting point was ``self-made'' \citep{Drechsler/Hartley} aside from Weyl's original papers.}
but with the peculiar idea of  considering the Higgs field as a part of the gravitational structure.

In the {\em new century} this peculiar idea was superseded by the studies of Nishina and Rajpoot who continued the research opened up by Cheng and stayed closer to the mainstream expectations of a massive Higgs field which at that time  was still hypothetical (section \ref{subsection Nishino/Rajpoot}). With the empirical  detection of the Higgs quantum excitation (``particle'') this  line was accentuated as the most realistic among the Weyl geometric approaches to SM fields. But the question how scale symmetry  is related to the quantum level remains still open. Ohanian's attempt for convincing us that scale symmetry is ``spontaneously broken'' near the Planck scale and leads back to Einstein gravity  developed a nice toy model (section \ref{subsection Ohanian}), but the investigations of Codello et al. show that the last word has not yet been spoken (section \ref{subsection Codello et al.})

With regard to astrophysics and cosmology the first exploratory phase of investigations, as it was called above, was superseded in the  Brazilian research tradition of Weyl geometric gravity, initiated by the work of Novello et al. Although this research line has been confined to a geometrically ``hybrid'' approach  which would imply a dynamically inert scalar field, if the Weyl geometric scaling symmetry would be taken seriously, this group of authors followed the physical intuition of their founding ``father'', or at least of the founding paper of the tradition \citep{Novello/Oliveira_ea:1992} which  assumed  an effective action  of a broken underlying scale symmetry  with some surviving residual Weylian terms  (section \ref{subsection Brazilian approach}). This allowed to investigate  concrete cosmological models which give some impression of what possibilities a Weyl geometric extension of present Riemann-Einstein cosmological models might offer. 
From a different side a bridge to the family of MOND-like theories of dark matter has been established;  it has widened the research horizon for the Weyl geometric extension of gravity  theories  even further (section \ref{subsection dark matter}). 

All in all the scale covariant, and often explicitly Weyl geometric approaches to gravity, elementary particle fields, foundations of quantum mechanics, astrophysics and cosmology have developed a rich panorama of models since the 1970s. In many cases less known scientists contributed to this research. Once in while  it  attracted the attention of internationally  renown physicists.   Although the Weyl geometric perspective remained  a side-stream in all of the mentioned fields up to now, it may well offer  interesting challenges and openings for the future. 

\subsection{\small \ldots and an open research horizon}
Our panorama has shown a variety of approaches which do yet not form a coherent research program as a whole. The Bohmian research lines, e.g.,  still stand separate from the other approaches, although some formal connections to the scalar fields of Dirac/Omote/Utiyama  type have been established in the later phase. It is not clear, however, whether the  ``Tehran''  perspective stands on solid grounds, and if so whether it can be integrated with the  ``Italian'' approach into a consistent common picture.   Other filaments of the whole field indicate perspectives which may reinforce each other.\footnote{The following perspectives of an ``open research horizon'' are necessarily subjective, but may be of help for orientation.}

Although cosmology  has increased its observational basis in the past few decades so tremendeously, it continues to call for alternative approaches to its many  conundrums. Several contributions to this volume 
deal with such alternatives.  The Weyl geometric approach joins this challenge,  although for the time being with minor strength.  The Brazilian work has made the most concrete contribution to this subject, but it is still hampered by the constraints resulting from the Palatini approach of variation  (section \ref{subsection Brazilian approach}). Moreover it  has not yet started to fully explore the consequences of the rescaling  freedom in Weyl geometric Robertson-Walker models and the possibility that part of the cosmological redshift may be  due to the scale connection rather than to a ``real''  expansion   (section \ref{subsection trivial}). Such a turn towards a more field theoretic explanation of the cosmological redshift would open new vistas for the geometry of cosmological model building.
 The Weyl geometric approach is clearly well suited for such investigations. 

Among the most recent papers we have come across first steps towards a field quantization scheme in a Weyl geometric environment,    which preserves scale symmetry at the quantum level (Codello et al., section \ref{subsection quantum level}). If this quantization procedure, or another one with the  property of preserving the scale symmetry, can be extended to the complete set of standard model fields plus the Weyl geometric scale connection and the gravitational scalar  field, we may arrive at a {\em modest} integration of gravity and the SM, in which only  the scale degree of freedom of the metric is quantized. Bars/Steinhardt/Turok have already argued that a theory with scale symmetry at the quantum level may   lead to a cancelling of the quadratically divergent terms in the radiative corrections to the Higgs mass, which constitute   the hard core of the  {\em naturalness problem} in present elementary particle physics.\footnote{\citep[p. 2]{Bars/Steinhardt/Turok:2014} }

This is a highly  interesting observation, although still  an unproven expectation. Together with the long standing speculations of the scalar field and/or the ``Weylon'' (scale connection) field as candidates for dark matter (sections \ref{subsection Dirac followers} and \ref{subsection Israelit}) the Weyl geometric approach seems to offer chances for attacking the naturalness problem of the SM  and the dark matter problem jointly,  essentially by  extending  the underlying autormorphism group of gravity and field theory. This complex of expectations has fed much of the research dynamics of the supersymmetry program; here we seem to be approaching  a similar thematic complex in a  more modest  form. We  also have seen that  a  classical, ``effective'' view of the gravitational field can lead to MOND-like phenomenology if also unusual kinematical terms for the scalar field are taken into account (section \ref{subsection dark matter}).   We may thus look forward with interest and curiosity to see  what the future research will lead to. \\[1em]


\small
\subsection*{\small Acknowledgements}
This paper owes its existence to {\em David Rowe}'s initiative  in several respects. He  encouraged me to present heterodox ideas on  Weyl geometric methods in cosmology at the Mainz conference and invited  me to rethink the case after a cool reception of the talk by the other participants. That gave me the chance  to place my views in the wider range of the recent attempts for using Weyl geometric methods in physics.    After an  interruption of several years,  an earlier first  draft of this paper \citep{Scholz:2011MainzarXiv} had be to be rewritten completely for the final version of this book. The new version overlaps nicely with  the wider ambit of the investigations of the interdisciplinary  group {\em Epistemology of the LHC} with center at Wuppertal, generously supported   by the  DFG/FWF.\footnote{Preprint number ELHC\underline{\hspace{0.3em}}2017-002. } This group offers the chance for a close communication between historians and philosophers of science and colleagues from the elementary particle community. 
 David  generously accepted the resulting oversize of the paper.
 \vspace{2em}

\small
\bibliographystyle{apsr}
  \bibliography{a_lit_hist,a_lit_mathsci,a_lit_scholz}

\end{document}